\newcommand{\Mg}{({\mathcal M},{\mathfrak  g}) } 
\newcommand{\M}{{\mathcal M}}               

\newcommand{\R}{I\!\!R}                     
\newcommand{\Z}{{\mathbb Z}}                
\newcommand{\C}{{\mathbb C}}                
\newcommand{\N}{{I\!\!N}}                   

\newcommand{\Jts}{{\mathbb J}}   
\newcommand{\Jt}{J}              
\newcommand{\Jms}{{\mathfrak  J}}     
\newcommand{\Jm}{{\mathcal J}}   

\newcommand{\Iota}[2]{{\mathfrak  D}_{#1,#2}}

\newcommand{\sgn}{\mbox{\sl sgn}}
\newcommand{\dgn}{\mbox{\sl dgn}}
\newcommand{\rk}{\mbox{\sl rk}}
\newcommand{\mul}{\mbox{\sl mul}\,}
\newcommand{\Bge}[2]{{\rm B}_{\rm sym}^{\vphantom{a}^{\ge #1}}(#2,\R)}

\newcommand{\ddt}{\frac{{\rm d}^2}{{\rm d}t^2}}

\newcommand{\mycontentsline}[4]{%
 \ifnum#1>1\medskip\else\bigskip\fi
 \hbox to\hsize{%
 \ifnum#1>1\hspace{1cm}\fi%
 \hbox to1cm{#2\hfil}%
 \vbox{\advance\hsize by-#1cm\parindent0pt\parfillskip0pt #3\dotfill#4}%
 }}

\newcommand{\mycontentsbigline}[4]{%
 \ifnum#1>1\medskip\else\bigskip\fi%
 \newdimen\auxiliar
 \hbox to\hsize{%
 \ifnum#1>1\hspace{1cm}\fi%
 \hbox to1cm{#2\hfil}%
 \vtop{\advance\hsize by-#1cm\parindent0pt\parfillskip0pt
 \auxiliar=\hsize\advance\auxiliar by-1cm
 \parshape=2 0pt \auxiliar 0pt \hsize #3\dotfill#4}%
 }}

\makeatletter
\newcommand{\mysubsection}{\@startsection
{subsection}%
{2}%
{0mm}%
{-\baselineskip}%
{0.5\baselineskip}%
{\normalfont\normalsize\itshape}
}%
\makeatother

\documentclass[oneside,titlepage,letterpaper,draft,12pt]{amsart}

\usepackage{amsmath}               
\usepackage{amsthm}
\usepackage{amscd}
\usepackage{times}

\numberwithin{equation}{subsection}   

\title[{\tiny Stability of the Focal and Geometric Index}]%
{Stability of the Focal and Geometric Index in semi-Riemannian
Geometry via the Maslov Index}
\author[{\tiny F.\ Mercuri}]{Francesco Mercuri}
\address{Departamento de Matem\'atica,\hfill\break\indent IMECC,
Universidade Estadual de Campinas\hfill\break\indent  SP, Brazil}
\email{mercuri@ime.unicamp.br}
\author[{\tiny P.\ Piccione}]{Paolo Piccione}
\address{Departamento de Matem\'atica,\hfill\break\indent IME, Universidade de S\~ao
Paulo\hfill\break\indent SP, Brazil}
\email{piccione@ime.usp.br}
\urladdr{http://www.ime.usp.br/\~{}piccione}
\author[{\tiny D.\ V.\ Tausk}]{Daniel V.\ Tausk}
\address{Departamento de Matem\'atica,\hfill\break\indent IME,  Universidade de S\~ao
Paulo\hfill\break\indent SP, Brazil}
\email{tausk@ime.usp.br}
\urladdr{http://www.ime.usp.br/\~{}tausk}
\thanks{{\it 1991 Mathematics Subject Classification.} 34B24, 34C10, 34L05, 53C22, 53C50, 53C80}
\date{August 1999}

\theoremstyle{definition}\newtheorem{defsym}{Definition}[subsection]
\theoremstyle{definition}\newtheorem{defsiggeom}[defsym]{Definition} 
\theoremstyle{definition}\newtheorem{deffocal}{Definition}[subsection]
\theoremstyle{definition}\newtheorem{defsigdiff}[deffocal]{Definition} 
\theoremstyle{plain}\newtheorem{correspondence}{Proposition}[subsection]
\theoremstyle{plain}\newtheorem{slightly}[correspondence]{Lemma}
\theoremstyle{remark}\newtheorem{remlight}[correspondence]{Remark}
\theoremstyle{plain}\newtheorem{matatudo}[correspondence]{Proposition} 
\theoremstyle{remark}\newtheorem{unifying}[correspondence]{Remark}
\theoremstyle{definition}\newtheorem{defomega}{Definition}[subsection] 
\theoremstyle{plain}\newtheorem{isolated}[defomega]{Proposition} 
\theoremstyle{plain}\newtheorem{riemlorcausal}[defomega]{Corollary}
\theoremstyle{remark}\newtheorem{analytic}[defomega]{Remark}

\theoremstyle{remark}\newtheorem{remcomplLagr}{Remark}[subsection]
\theoremstyle{definition}\newtheorem{cartas}{Definition}[subsection]
\theoremstyle{plain}\newtheorem{geometriaLambda}[cartas]{Proposition}
\theoremstyle{plain}\newtheorem{LieAl}[cartas]{Proposition}
\theoremstyle{definition}\newtheorem{defH*}[cartas]{Definition} 
\theoremstyle{plain}\newtheorem{UnOn}[cartas]{Proposition}
\theoremstyle{definition}\newtheorem{defLambdak}[cartas]{Definition}
\theoremstyle{remark}\newtheorem{abertofechado}[cartas]{Remark}
\theoremstyle{plain}\newtheorem{transitividade}[cartas]{Proposition}
\theoremstyle{plain}\newtheorem{complcomun}[cartas]{Corollary}
\theoremstyle{plain}\newtheorem{geometriaLambdak}[cartas]{Proposition}
\theoremstyle{remark}\newtheorem{remtangent}[cartas]{Remark}

\theoremstyle{plain}\newtheorem{quotient}{Lemma}[subsection]
\theoremstyle{plain}\newtheorem{pi1}[quotient]{Corollary}
\theoremstyle{remark}\newtheorem{induced}{Remark}[subsection]
\theoremstyle{plain}\newtheorem{preliminary}[induced]{Lemma}
\theoremstyle{plain}\newtheorem{elementary}[induced]{Lemma}
\theoremstyle{plain}\newtheorem{homologous}[induced]{Proposition}
\theoremstyle{definition}\newtheorem{defmaslov}[induced]{Definition}
\theoremstyle{plain}\newtheorem{method}{Proposition}[subsection]
\theoremstyle{plain}\newtheorem{estimate}[method]{Corollary}
\theoremstyle{plain}\newtheorem{signature}[method]{Corollary}

\theoremstyle{definition}\newtheorem{MaslovDiff}{Definition}[subsection]
\theoremstyle{plain}\newtheorem{MaslovFocal}[MaslovDiff]{Theorem}
\theoremstyle{plain}\newtheorem{MaslovCausal}[MaslovDiff]{Corollary}
\theoremstyle{plain}\newtheorem{LimitMaslov}{Theorem}[subsection]
\theoremstyle{plain}\newtheorem{LimitCausal}[LimitMaslov]{Corollary}
\theoremstyle{remark}\newtheorem{uniformity}[LimitMaslov]{Remark}

\theoremstyle{plain}\newtheorem{technical}{Lemma}[subsection]
\theoremstyle{plain}\newtheorem{eigenvalues}[technical]{Proposition}
\theoremstyle{plain}\newtheorem{finiteeigen}[technical]{Corollary}
\theoremstyle{remark}\newtheorem{remuniform}[technical]{Remark}
\theoremstyle{definition}\newtheorem{defindspectral}[technical]{Definition}
\theoremstyle{plain}\newtheorem{restrictions}{Lemma}[subsection]
\theoremstyle{remark}\newtheorem{remerror}[restrictions]{Remark}
\theoremstyle{plain}\newtheorem{indlambda}[restrictions]{Corollary}
\theoremstyle{plain}\newtheorem{Morse}[restrictions]{Theorem}
\theoremstyle{remark}\newtheorem{remMorse}[restrictions]{Remark}

\theoremstyle{plain}\newtheorem{Ck}{Lemma}[subsection]
\theoremstyle{plain}\newtheorem{dense}[Ck]{Lemma}

\begin{document}
\begin{abstract}
We investigate the problem of the stability of the number of conjugate or focal
points (counted with multiplicity) along a semi-Riemannian geodesic $\gamma$.
For a Riemannian or a non spacelike Lorentzian geodesic, such number is equal
to the intersection number (Maslov index) of a continuous curve with a subvariety 
of codimension one of the Lagrangian Grassmannian of a symplectic space.
Such intersection number is proven to be stable in a large variety of circumstances.
In the general semi-Riemannian case, under suitable hypotheses this number 
is equal to an algebraic count of the multiplicities of the conjugate points, and
it is related to the spectral properties of a non self-adjoint differential
operator. This last relation gives a weak extension of the classical Morse Index Theorem
in Riemannian and Lorentzian geometry. In this paper we reprove some results
that were incorrectly stated by Helfer in \cite{Hel1}; in particular, a counterexample
to one of Helfer's results, which is essential for the theory, is given. 
In the last part of the paper we discuss a general technique for the construction
of examples and counterexamples in the index theory for semi-Riemannian
metrics, in which some new phenomena appear.
\end{abstract}
\maketitle
\section*{Contents}
\bigskip
\mycontentsline{1}{\ref{sec:intro}}{Introduction}{\pageref{sec:intro}}
\mycontentsline{1}{\ref{sec:preliminaries}}{Preliminaries}{\pageref{sec:preliminaries}}
\mycontentsline{2}{\ref{sub:geometrical}}{The geometrical problem}{\pageref{sub:geometrical}}
\mycontentsline{2}{\ref{sub:differential}}{The differential problem}{\pageref{sub:differential}}
\mycontentsline{2}{\ref{sub:geomdiff}}{Relations between the geometrical and the differential
problem}{\pageref{sub:geomdiff}}
\mycontentsline{2}{\ref{sub:symplstruct}}{The symplectic structure associated to a differential
problem}{\pageref{sub:symplstruct}}
\mycontentsline{2}{\ref{sub:isolatedness}}{On the discreteness of the set of $(P,S)$-focal
instants}{\pageref{sub:isolatedness}}
\mycontentsline{1}{\ref{sec:geometry}}{%
Geometry of the Lagrangian Grassmannian}{\pageref{sec:geometry}}
\mycontentsline{2}{\ref{sub:generalities}}{Generalities on symplectic
spaces}{\pageref{sub:generalities}}
\mycontentsline{2}{\ref{sub:lagrangian}}{The Lagrangian Grassmannian}{\pageref{sub:lagrangian}}
\mycontentsline{1}{\ref{sec:intersection}}{%
Intersection Theory: the Maslov Index}{\pageref{sec:intersection}}
\mycontentsline{2}{\ref{sub:fundamental}}{The fundamental group of the Lagrangian
Grassmannian}{\pageref{sub:fundamental}}
\mycontentsline{2}{\ref{sub:construction}}{The intersection theory and the construction of the
Maslov index}{\pageref{sub:construction}}
\mycontentsline{2}{\ref{sub:computation}}{Computation of the Maslov
index}{\pageref{sub:computation}}
\mycontentsline{1}{\ref{sec:applications}}{%
Applications of the Maslov Index: Stability of the Geometric Index}{\pageref{sec:applications}}
\mycontentsline{2}{\ref{sub:maslovdiff}}{The Maslov index of a differential
problem}{\pageref{sub:maslovdiff}}
\mycontentsline{2}{\ref{sub:stability}}{Stability of the indexes}{\pageref{sub:stability}}
\mycontentsbigline{1}{\ref{sec:spectral}}{%
The Spectral Index. Some Remarks on a Possible Extension
of the Morse Index Theorem}{\pageref{sec:spectral}}
\mycontentsline{2}{\ref{sub:spectral}}{Eigenvalues of the differential problem and the spectral
index}{\pageref{sub:spectral}}
\mycontentsline{2}{\ref{sub:morse}}{A generalized Morse Index Theorem}{\pageref{sub:morse}}
\mycontentsline{1}{\ref{sec:final}}{%
Curves of Lagrangians originating
from differential problems}{\pageref{sec:final}}
\mycontentsbigline{2}{\ref{sub:distributions}}{Differential problems determine curves of
Lagrangians that are tangent to distributions of affine spaces}{\pageref{sub:distributions}}
\mycontentsline{2}{\ref{sub:localcoordinates}}{A study of curves of Lagrangians in local
coordinates}{\pageref{sub:localcoordinates}}
\mycontentsline{2}{\ref{sub:nondegenerate}}{The case where $g$ is nondegenerate on ${\rm
Ker}(\beta(t_0))$}{\pageref{sub:nondegenerate}}
\mycontentsline{2}{\ref{sub:counterexample}}{A counterexample to the equality $\mu(g,R,P,S)={\rm
i}_{\rm foc}$}{\pageref{sub:counterexample}}
\mycontentsline{2}{\ref{sub:instabilityfocal}}{Instability of the focal
index}{\pageref{sub:instabilityfocal}}
\mycontentsline{2}{\ref{sub:instabilitypoints}}{Instability of focal points with signature
zero}{\pageref{sub:instabilitypoints}}
\newpage
\begin{section}{Introduction}\label{sec:intro}
The original motivation for writing this paper was given by the following
problem. Given a sequence $\gamma_j$ of geodesics in a semi-Riemannian
manifold $\Mg$ converging to a geodesic $\gamma_\infty$, what can be said 
about the convergence of the {\em geometric index\/} of $\gamma_j$ to that of
$\gamma_\infty$? The question arose in the context of
Lorentzian geometry, where the problem originated in an attempt to develop
a Morse theory for lightlike geodesics as {\em limit\/} of the theory
for timelike geodesics.

Recall that the geometric index of a geodesic
$\gamma:[a,b]\mapsto\M$ is the number of points that are 
conjugate to $\gamma(a)$ along $\gamma$, 
counted with multiplicity. Conjugate points along a geodesic 
correspond to the zeroes of (non
trivial) Jacobi fields along $\gamma$, which are vector fields 
that annihilate the 
{\em index form\/} $I_\gamma$. The index form $I_\gamma$ is a
symmetric bilinear form defined in the space of vector fields
along $\gamma$ that is written in terms of the Levi-Civita connection and
the curvature tensor of $\mathfrak g$; the celebrated {\em Morse Index Theorem\/}
(see~\cite{BE, BEE, dC, EK, M, ON, Sm} for different versions of this theorem)
states that the geometric index of a Riemannian or causal Lorentzian geodesic
$\gamma$ is equal to the number (with multiplicity) of negative eigenvalues of 
$I_\gamma$, provided that the final point
$\gamma(b)$ is not considered  in the count of conjugate points.
The number (with multiplicity) of negative eigenvalues of a symmetric bilinear
form on a vector space is called the {\em index\/} of the form; 
a simple argument shows that if a symmetric bilinear form 
is continuous with respect to some norm in the vector space $V$, 
then its index does not change
when one extends $I$ to the Banach space completion of $V$.

It is not too hard to prove that the convergence of a sequence of geodesics
to a geodesic implies a {\em strong\/} convergence of the corresponding
index forms. Hence, by considering suitable Hilbert space completions
of the set of vector fields along $\gamma$ and representing $I_\gamma$
as a self-adjoint operator on this Hilbert space, the problem of convergence
of the geometric index can be reformulated in terms of convergence
of the index of a sequence of self-adjoint operators converging
in the operator norm. However,  it is very easy to give examples of sequences
of (real) symmetric matrices with constant index converging to a symmetric matrix
having a different index. In finite dimension, this phenomenon arises
only when the limit matrix is non invertible, and in the case that, when
passing to the limit, some
negative eigenspace of the matrices of sequence {\em falls\/} into the kernel
of the limit, causing a drop in the index. In the infinite dimensional case
the situation is even worse, and one can have a sequence of self-adjoint
operators converging to an invertible self-adjoint operator whose index
is strictly less that the infimum of the indices of the approximating
family.

Some questions concerning the continuity of the conjugate points in Riemannian 
geometry are studied in reference~\cite{Pat}.

If $\Mg$ is Riemannian, i.e.,  $\mathfrak g$ is positive definite,
then, considering a suitable $H^1$-Sobolev completion of the space of
vector fields along the geodesic, the self-adjoint operator associated
to the index form is a {\em compact perturbation\/} of the identity.
In this case, if the limit is an  invertible operator, i.e., if the
point $\gamma_\infty(b)$ is not conjugate to $\gamma_\infty(a)$ along
$\gamma_\infty$, the geometric index of $\gamma_j$ is eventually constant,
and equal to the geometric index of $\gamma_\infty$. 
The same conclusion holds
for timelike Lorentzian geodesics, provided that the limit be also timelike.
In this case, we consider  the index of the index form
defined only in the space of vector fields which are
everywhere orthogonal to the geodesic, and the Lorentzian metric is positive
on such fields. 

When one considers the case of lightlike geodesics, though, the situation is 
complicated by the fact that
the index form always has a non trivial kernel, 
even when the final endpoint
is non conjugate to the initial one. Namely, any vector 
field which is a multiple of
the tangent field to the geodesic is in the kernel of the index form.
The presence of the kernel in the lightlike case could be avoided
by considering suitable quotients of the normal bundle, but then
one loses the relation with the non lightlike geodesics,
for which an index form cannot be defined on such quotient.

Thus, using abstract spectral arguments one only proves the semi-continuity
of the index for Lorentzian causal geodesics.
\smallskip

A different technique to study the stability of the index is suggested
by an analogy with the classical {\em Sturm problem\/} in the theory
of ordinary differential equations. The Sturm oscillation theorem deals
with second order differential systems of the form $-(px')'+rx=\lambda x$
where $p$ and $r$ are functions with $p>0$, and $\lambda$ is a real parameter.
The theorem states that, denoting by $C^1_o[\alpha,\beta]$ the space of
$C^1$-functions on $[\alpha,\beta]$ vanishing at $\alpha$ and $\beta$, 
the index of the symmetric bilinear form 
$B(x,y)=\int_a^b[px'y'+r xy]\;{\rm d}t$ in $C^1_o[a,b]$ is equal  to the  
sum over $t\in\,]a,b[$ of the dimension of the kernel of the bilinear form 
$\int_a^t[px'y'+r xy]\;{\rm d}t$ in $C^1_o[a,t]$.

The classical proof of the Sturm oscillation theorem (\cite[Chapter~8]{CL}) is obtained by
showing that the two quantities involved in the thesis can be obtained
as the {\em winding number\/} of two homotopic closed curves in the  real projective 
line. As a {\em side effect\/} of this theory, one obtains immediately that,
since the winding number is stable by homotopies, and in particular by small
$C^0$-perturbations, the index of $B$ is stable by small perturbations.

So, the stability of the index for a Sturm system is proven by 
relating the index form $B$ to some homotopical invariants of
the problem.
In this paper we exploit this method to obtain the stability
of the geometrical index, or of suitable substitutes of it, for
semi-Riemannian geodesics. More precisely, we investigate the
notion of {\em Maslov index\/} for a Jacobi type system of ordinary
differential equations, which is obtained as the {\em intersection number\/}
between a curve and a subvariety of codimension one of
a smooth manifold. The Maslov index
of a system arising from the Jacobi equation of a Riemannian or
a causal Lorentzian geodesic is then proven to be equal to the geometric
index of the geodesic.

The idea and some of the results presented are not new.

In \cite{Bott} and \cite{Ed} the authors develop an approach to the index
problem through topological methods.
For instance, in~\cite{Ed}, it is employed a similar technique to prove 
a generalization of the Sturm's Theorems in the case of an arbitrary 
self-adjoint system
of even order and subject to arbitrary self-adjoint boundary conditions.
The main technical tool used in the proof is  the notion of {\em $U$-manifold},
which is obtained from an even dimensional vector space $E$, endowed with
a non degenerate Hermitian form $\psi$, as the set of all maximal subspaces
of $E$ on which $\psi$ vanishes. Non trivial solutions of the eigenvalue
boundary value problem determine continuous curves in a $U$-manifold, and
the proof of the result is based in studying the number of intersections
of such curves with a subvariety of codimension one. The intersection theory
developed by the author for $U$-manifolds is based on the {\em relative
homotopy\/} theory.

In~\cite{Hel1} (see also \cite{Hel2}), the author carries out a similar analysis for
Morse--Sturm--Liouville systems, which are symmetric with respect
to an indefinite inner product. In this context, the environment
for the intersection theory is given by the set of all Lagrangians
of a symplectic space which is naturally associated to the differential
system. This set, the {\em Lagrangian Grassmannian\/} of the symplectic
space, has a natural manifold structure, which is in general non orientable.
A non trivial solution of the differential system gives a continuous curve
in the Lagrangian Grassmannian, and the zeroes of the solution correspond
to intersections of this curve with the subvariety of all Lagrangians $L$
which are not complementary to a given one $L_0$. The intersection
theory employed in this situation is based on the notion of {\em Maslov index\/}
of a curve, obtained using duality in homology theory (Arnol'd--Maslov cycle).

This approach does not seem to deal properly with the lack of orientability
of the Lagrangian Grassmannian.  Furthermore, several results of Helfer's
paper are incorrectly stated, due mainly to the lack of an essential assumption
of {\em nondegeneracy\/} for the restrictions of certain symmetric bilinear forms.
More precisely, in \cite[Proposition~5.1 (b)]{Hel1} it is claimed the equality
between the Maslov index and the sum of the signatures of the conjugate
points; we give a counterexample to such equality in subsection~\ref{sub:counterexample}.
In \cite[Proposition~6.2]{Hel1}, the proof is incomplete, because
the case of non simple eigenvalues is not treated properly (see Remark~\ref{thm:remerror}). 
Consequently, also the proof of \cite[Proposition~7.1]{Hel1} is affected by these problems;
a more restrictive statement of the Index Theorem 
is proven in Section~\ref{sec:spectral}.

It is important to observe that the possibility of such degeneracies,
which do not occur in Riemannian or causal Lorentzian geodesics,
is responsible for many of the new phenomena which arise in the general
semi-Riemannian case, like for instance, the possibility of accumulation
of the conjugate points along a geodesic. Curiously enough, also in the
book by O'Neill \cite[Exercise~8, page~299]{ON}, the non degeneracy
assumption is missing, and the author claims incorrectly that the
set of focal points along a geodesic is discrete.

For this reasons, we have opted to provide an alternative, self-contained, 
presentation of some of Helfer's results needed for our proof 
of the stability of the geometrical index. For the sake of completeness,
the geometrical results of \cite{Hel1} that are being considered in 
this paper have been reproven in the slightly more general context of 
geodesics starting orthogonally to a given non degenerate submanifold 
$\mathcal P$ of $\mathcal M$. With such a generalization we are able to
prove our stability results also for the {\em focal index\/} of a geodesic
relative to a submanifold.

In order to keep our notation as uniform as possible and to make the results
accessible to mathematicians and physicists from different areas,  in the paper 
we make an effort to give a formal proof of almost everything we claim, 
even though this attitude has the disadvantage of not distinguishing between
new and old results. For instance, the statement and the proof of
Proposition~\ref{thm:correspondence} in this paper, apart from the result 
of Lemma~\ref{thm:slightly}, is almost identical to \cite[Section~3]{Hel1};
the proof of Lemma~\ref{thm:restrictions} is essentially contained
in the proof of \cite[Proposition~6.2]{Hel1}, and
also the proof of Theorem~\ref{thm:Morse} is identical to the proof of a somewhat
similar result proven in \cite[Theorem~7.1]{Hel1}. It should also be remarked
that many of the results concerning the geometry of the Lagrangian Grassmannian
of a symplectic space  presented in Section~\ref{sec:geometry} may appear in
similar forms on other references, like for instance \cite{Ar,GS,Tre}.

We outline briefly the contents of this paper.
In Section~\ref{sec:preliminaries} we introduce the notations and we give a few
preliminary results on the geometrical and the differential framework.
In Section~\ref{sec:geometry} we present a detailed description of the analytical
structure of the Lagrangian Grassmannian of a symplectic space. All the results
are given in an intrinsic, i.e., coordinate independent, form.

Section~\ref{sec:intersection} is devoted to the intersection theory used for
the definition and the properties of the Maslov index. The treatment presented
is inspired by the intersection theory for {\em $U$-manifolds\/} of \cite{Ed};
this approach has the advantage of avoiding the  homology/cohomology duality
issues  in non orientable spaces. We remark that the Lagrangian Grassmannian
is in general a non orientable manifold (see~\cite{Fu}).
Some of the results presented in Section~\ref{sec:intersection} are
already proven in \cite{Ar}, although only in the case of closed curves.
For the computation of the Maslov index of a loop, Arnol'd uses the isomorphism of 
the fundamental groups of  the Lagrangian Grassmannian and of the circle
$S^1$ induced by the square of the determinant function. 
For our purposes, the techniques
developed in Section~\ref{sec:intersection} using suitable coordinate charts,
are more appropriate for the computation of the Maslov index (see
subsection~\ref{sub:computation}). Namely, we describe the Maslov index in terms
of change of signatures of symmetric bilinear forms, obtaining a natural
relation with the geometric index (see Theorem~\ref{thm:MaslovFocal}
and the following corollary).

In Section~\ref{sec:applications} we apply the results of the previous sections
to the problem of the stability of the geometrical and focal indexes in semi-Riemannian
geometry. A special attention is given to the case of the approximation of a lightlike
geodesic by timelike geodesics in a Lorentzian manifold. If $\gamma_\infty$ is
a lightlike geodesic in a Lorentzian manifold whose endpoints are not conjugate and
$\gamma_n$ is a sequence of timelike geodesics converging to $\gamma_\infty$, then
the geometric index of $\gamma_n$ is eventually constant and equal to the geometric
index of $\gamma_\infty$. 

In Section~\ref{sec:spectral} we present a spectral theorem
(Theorem~\ref{thm:Morse}), which is a weak generalization
of the classical Morse Index Theorem for Riemannian or causal Lorentzian geodesics
to the case of geodesics in an arbitrary semi-Riemannian manifold. The proof of
the spectral theorem is obtained by a simple homotopy argument in the Lagrangian 
Grassmannian. It is interesting to observe that the proof of Theorem~\ref{thm:Morse}
gives an alternative  and unifying  proof of all the previous versions of
the Morse Index Theorem in Riemannian and Lorentzian geometry (see \cite{M}). 

In Section~\ref{sec:final}, we study the problem of determining
which curves of Lagrangians are associated to solutions of
Jacobi systems. We give some necessary and sufficient 
conditions for this occurrence, and we use these conditions to find
examples of situations described generically in the rest of the paper.
Remarkably, we give an example in which the equality between the 
Maslov index and the focal index of a Lorentzian spacelike geodesic $\gamma$
fails, due to the degeneracy of the metric on the space $\Jms[t] $ consisting of
the evaluations at $t$ of the Jacobi fields along $\gamma$ that vanish at
the initial instant.
\smallskip

\noindent{\bf Acknowledgments.}\enspace
Several people have given support to the authors during the writing of this
paper. Particularly, the authors wish to express their gratitude to Professor
Daniel Baumann Henry, Professor Daciberg Lima Gon\c calves and Professor Volker Perlick
who have given many fruitful suggestions. 
\end{section}
\newpage
\begin{section}{Preliminaries}\label{sec:preliminaries}
Let $\Mg$ be a smooth semi-Riemannian manifold, i.e.,
$\M$ is a finite dimensional real  $C^\infty$ manifold 
whose  topology satisfies the second countability axiom and 
the  Hausdorff separation axiom, and  $\mathfrak g$ is a 
a smooth $(2,0)$-tensor field on $\M$ which defines a non degenerate
symmetric bilinear form on each tangent space of $\M$. We set
$m={\rm dim}(\M)$; a non zero vector $v\in T\M$ will be called {\em spacelike},
{\em timelike} or {\em lightlike\/} according to $\mathfrak  g(v,v)$ being
positive, negative or null.

We give some general definitions 
concerning symmetric bilinear forms for later use.
\begin{defsym}\label{thm:defsym}
Let $V$ be any real vector space and $B:V\times V\mapsto\R$ a symmetric
bilinear form. The {\em negative type number\/} (or {\em index}) $n_-(B)$ of $B$
is the possibly infinite number defined by
\begin{equation}\label{eq:def-}
n_-(B)=\sup\Big\{{\rm dim}(W):W\ \text{subspace of}\ V \ \text{on which}\ B\ 
\text{is negative definite}\Big\}.
\end{equation}
The {\em positive type number\/} $n_+(B)$ is given by
$n_+(B)=n_-(-B)$; if at least one of these two numbers is finite,
the {\em signature\/} $\sgn(B)$ is defined by:
\[\sgn(B)=n_+(B)-n_-(B).\]
The {\em kernel\/} of $B$, ${\rm Ker}(B)$, is the set $V^\perp$ of vectors $v\in V$
such that $B(v,w)=0$ for all $w\in V$; the {\em degeneracy\/} $\dgn(B)$
of $B$ is the (possibly infinite) dimension of ${\rm Ker}(B)$. 
\end{defsym}
If $V$ is finite dimensional, then the numbers $n_+(B)$, $n_-(B)$ and
$\dgn(B)$ are respectively  the number of $1$'s, $-1$'s and $0$'s
in the canonical form of $B$ as given by the Sylvester's Inertia Theorem.
In this case, $n_+(B)+n_-(B)$ is equal to the codimension of ${\rm Ker}(B)$,
and it is also called the {\em rank} of $B$, $\rk(B)$.
\smallskip

The semi-Riemannian manifold $\Mg$ is said to be {\em Lorentzian\/}
if the index of $\mathfrak  g$ at every point of $\M$ is equal to $1$.
A four dimensional Lorentzian manifold is the mathematical model
for a general relativistic spacetime; in this case, the
timelike and lightlike geodesics in $\M$, i.e., geodesics $\gamma$ with
$\gamma'(t)$  timelike or lightlike  for all $t$ respectively, represent
the trajectories of massive and massless objects freely falling
under the action of the gravitational field.
\smallskip

Let $\nabla$ denote the covariant derivative of the Levi--Civita connection
of $\mathfrak  g$ and let ${\mathcal R}(X,Y)=\nabla_X\nabla_Y-\nabla_Y\nabla_X-\nabla_{[X,Y]}$ 
be the curvature tensor of $\mathfrak  g$. If $V:[a,b]\mapsto T\M$ is a vector field
along a given curve $\gamma:[a,b]\mapsto\M$, we will write $V'$ for the
covariant derivative of $V$ along $\gamma$.

\mysubsection{The geometrical problem}
\label{sub:geometrical}

Let $\mathcal P$ be a smooth submanifold of $\M$, $p\in P$ and $n\in
T_p{\mathcal P}^\perp$, i.e., $n\in T_p\M$ and $\mathfrak  g(n,v)=0$ for all $v\in T_p{\mathcal P}$.
The {\em second fundamental form\/} of
$\mathcal P$ at 
$p$ in the direction $n$ is the symmetric bilinear form ${\mathcal S}_n:
T_p{\mathcal P}\times T_p{\mathcal P}\mapsto\R$ given by:
\[{\mathcal S}_n(v_1,v_2)={\mathfrak  g}(\nabla_{v_1}V_2,n),\]
where $V_2$ is any extension of $v_2$ to a vector field on $\mathcal P$.
If $T_p{\mathcal P}$ is {\em non degenerate}, i.e., if the restriction
of $\mathfrak  g$ to $T_p{\mathcal P}$ is non degenerate, then there exists 
a linear operator, still denoted by ${\mathcal S}_n$, on $T_p{\mathcal P}$,
such that ${\mathcal S}_n(v_1,v_2)={\mathfrak  g}({\mathcal S}_n(v_1),v_2)$
for all $v_1,v_2\in T_p{\mathcal P}$. Also, if $T_p{\mathcal P}$ is non degenerate,
then the second fundamental form can be viewed as a $T_p{\mathcal P}^\perp$-valued
symmetric bilinear form $\mathcal S$ on $T_p{\mathcal P}$ defined by:
\[{\mathcal S}(v_1,v_2)=\text{orthogonal projection of}\ \nabla_{v_1}V_2\ \text{onto}\ 
T_p{\mathcal P}^\perp,\]
so that
\[{\mathcal S}_n(v_1,v_2)=\mathfrak  g({\mathcal S}(v_1,v_2),n),\quad\forall\, v_1,v_2\in
T_p{\mathcal P},\ n\in T_p{\mathcal P}^\perp.\]
Let $\gamma:[a,b]\mapsto\M$ be a non constant geodesic orthogonal to $\mathcal P$ at the initial
point, i.e.,  $\gamma(a)\in{\mathcal P}$ and $\gamma'(a)\in T_{\gamma(a)}{\mathcal P}^\perp$.
Let's assume  that $T_{\gamma(a)}{\mathcal P}$ is non degenerate;
we will say that the family of objects $(\M,\mathfrak  g,\gamma, {\mathcal P})$
is an {\em admissible quadruple for the geometrical problem}.

A {\em Jacobi field\/} along $\gamma$ is  a smooth vector field $\Jm$
along $\gamma$ that satisfies the Jacobi equation:
\begin{equation}\label{eq:eqJacobi}
\Jm''={\mathcal R}(\gamma',\Jm)\,\gamma';
\end{equation}
we say that $\Jm$ is a ${\mathcal P}$-Jacobi field if, in addition, $\Jm(a)$ satisfies:
\begin{equation}\label{eq:PJacobi}
\Jm(a)\in T_{\gamma(a)} {\mathcal P}\quad\text{and}\quad
\Jm'(a)+{\mathcal S}_{\gamma'(a)}[\Jm(a)]\in T_{\gamma'(a)}{\mathcal P}^\perp. 
\end{equation}
Observe that, if $\mathcal P$ is a single point of $\M$, then \eqref{eq:PJacobi}
reduces to $\Jm(a)=0$. Geometrically, equation \eqref{eq:eqJacobi}
means that $\Jm$ is the variational vector field corresponding to a variation
of $\gamma$ by geodesics; condition \eqref{eq:PJacobi} says that these geodesics
are orthogonal to $\mathcal P$ at their initial points.

We define the following vector spaces:
\begin{eqnarray}
\label{eq:defJms}\Jms&=&\Big\{\Jm: \Jm\ \text{is a $\mathcal P$-Jacobi field along}\
\gamma\Big\};
\\
\label{eq:defJmsperp}
\Jms^\perp&=&\Big\{\Jm\in\Jms:\Jm(t)\in\gamma'(t)^\perp,\quad\forall\,t\in[a,b]\Big\}.
\end{eqnarray}
Observe that if $\Jm\in\Jms$ is such that $\Jm(t)\in\gamma'(t)^\perp$ for
some $t\in\,]a,b]$, then $\Jm\in\Jms^\perp$. Namely, for any Jacobi field
$\Jm$, the function ${\mathfrak  g}(\Jm(t),\gamma'(t))$ is affine on $[a,b]$,
and if $\Jm$ is $\mathcal P$-Jacobi, then ${\mathfrak  g}(\Jm(a),\gamma'(a))=0$;
similarly, for a $\mathcal P$-Jacobi field $\Jm$, the condition
$\Jm'(a)\in\gamma'(a)^\perp$ is equivalent to $\Jm\in\Jms^\perp$.
We conclude  that $\Jms^\perp$ can be described alternatively as the space
of Jacobi fields along $\gamma$ satisfying the initial conditions:
\begin{equation}\label{eq:incondalt}
\Jm(a)\in T_{\gamma(a)}{\mathcal P}\quad\text{and}\quad \Jm'(a)
+{\mathcal S}_{\gamma'(a)}[\Jm(a)]\in\big[T_{\gamma(a)}{\mathcal P}
\oplus\gamma'(a)\big]^\perp.
\end{equation}
It is easy to see that ${\rm dim}(\Jms)=m$; namely, the dimension of the subspace
of $T_{\gamma(a)}\M\oplus T_{\gamma(a)}\M$ defined by the initial conditions
\eqref{eq:PJacobi} is equal to $m$. Similarly, ${\rm dim}(\Jms^\perp)=m-1$,
since the dimension of the subspace
of $T_{\gamma(a)}\M\oplus T_{\gamma(a)}\M$ defined by the initial conditions
\eqref{eq:incondalt} is equal to $m-1$.

For all $t\in[a,b]$, we define the subspaces $\Jms[t]$ and $\Jms^\perp[t]$ 
of $T_{\gamma(t)}\M$ by
\[\Jms[t]=\Big\{\Jm(t):\Jm\in\Jms\Big\}\quad\text{and}
\quad\Jms^\perp[t]=\Big\{\Jm(t):\Jm\in\Jms^\perp\Big\};\]
observe that $\Jms^\perp[t]$ should not be confused with $\Jms[t]^\perp$,
which is the orthogonal complement of $\Jms[t]$ in $T_{\gamma(t)}\M$.
Precisely, the following relation holds:
\begin{equation}\label{eq:relJmsperpt}
\Jms^\perp[t]=\Jms[t]\cap\gamma'(t)^\perp,\quad\forall\,t\in\,]a,b];
\end{equation}
this follows immediately from the observation that if $t\in\,]a,b]$ and 
${\mathcal J}\in\Jms$ is such that
${\mathcal J}(t)\in\gamma'(t)^\perp$, then ${\mathcal J}\in \Jms^\perp$.
The vector field $(t-a)\,\gamma'(t)$ is always a $\mathcal P$-Jacobi field,
thus $\gamma'(t)\in \Jms[t]$ for $t\in\,]a,b]$; observe that $\gamma'(t)
\in\Jms^\perp[t]$ if and only if $\gamma$ is lightlike. We also remark
that, for $t\in\,]a,b]$, the following formula holds:
\begin{equation}\label{eq:codim1}
{\rm dim}( \Jms[t])={\rm dim}(\Jms^\perp[t])+1,\quad t\in\,]a,b].
\end{equation}
Indeed, the codimension of $\Jms^\perp$ in $\Jms$ is $1$; moreover, for $t\in\,]a,b]$,
the linear operator $\Jm\mapsto \Jm(t)$ of evaluation at $t$ has the same kernel
in $\Jms$ and in $\Jms^\perp$.

\smallskip

The point $\gamma(t_0)$, $t_0\in\,]a,b]$ is said to be a
{\em $\mathcal P$-focal point\/} along $\gamma$ if there exists
a non zero $\mathcal P$-Jacobi field $\mathcal J$ along $\gamma$
with ${\mathcal J}(t_0)=0$. The {\em multiplicity\/} $\mul(t_0)$
of $\gamma(t_0)$ is the dimension of the space of $\mathcal P$-Jacobi 
fields that vanish at $t_0$; we set $\mul(t_0)=0$ if $\gamma(t_0)$ is
not $\mathcal P$-focal:
\begin{equation}\label{eq:defmul}
\mul(t_0)={\rm dim}\, \big\{\Jm\in\Jms:\Jm(t_0)=0\big\} .
\end{equation}

Observe that if $\Jm\in\Jms$ vanishes at $t_0$, then $\Jm\in\Jms^\perp$;
since ${\rm dim}(\Jms)\!\!=\!{\rm dim}(T_{\gamma(t_0)}\M)$ and ${\rm dim}(\Jms^\perp)=
{\rm dim}(\gamma'(t_0)^\perp)$, we conclude that $\mul(t_0)$
is equal to the codimension of $\Jms[t_0]$ in $T_{\gamma(t_0)}\M$, and
also equal to the codimension of $\Jms^\perp[t_0]$ in $\gamma'(t_0)^\perp$.

We remark that \eqref{eq:codim1} implies that $\Jms^\perp[t]$ is a {\em proper\/}
subspace of $\Jms[t]$, for $t\in\,]a,b]$; 
this fact is trivial if $\gamma$ is not lightlike,
but it has interesting consequences otherwise. For instance, if ${\rm dim}(\M)=2$,
this implies that there can be no focal points along any lightlike geodesics.
Namely, if $t_0\in\,]a,b]$ is $\mathcal P$-focal, then $\Jms[t_0]$ is at the
most one dimensional, and therefore $\Jms^\perp[t_0]=\{0\}$, contradicting
the fact that $\gamma'(t_0)\in \Jms^\perp[t_0]$.

If the number of $\mathcal P$-focal points along $\gamma$ is finite,
one defines the {\em geometric index}, ${\rm i}_{\rm geom}(\gamma)$,
of $\gamma$ relative to the initial submanifold $\mathcal P$ to be the sum
of the multiplicities of the $\mathcal P$-focal points:
\[{\rm i}_{\rm geom}(\gamma)=\sum_{t\in\,]a,b]}\mul(t).\]

In order to extend to the semi-Riemannian case the classical
Morse Theory, we need to introduce the concept of {\em signature\/}
for a ${\mathcal P}$-focal point. 
\begin{defsiggeom}\label{thm:defsiggeom}
If $\gamma(t_0)$ is a ${\mathcal P}$-focal point, its {\em signature\/}
$\sgn(t_0)$ is defined to be the signature of the restriction
of the metric $\mathfrak  g$ to the space $\Jms[t_0]^\perp$.
If $\gamma(t_0)$ is not a ${\mathcal P}$-focal point, we set $\sgn(t_0)=0$.
The {\em focal index\/} ${\rm i}_{\rm foc}(\gamma)$
of the geodesic $\gamma$ relative to the initial submanifold $\mathcal P$
is defined by the sum:
\[{\rm i}_{\rm foc}(\gamma)=\sum_{t\in\,]a,b]}\sgn(t),\]
provided that the number of ${\mathcal P}$-focal points along $\gamma$ is finite.
\end{defsiggeom}
Sufficient conditions for the finiteness of the number of ${\mathcal P}$-focal points
will be discussed at the end of this Section (see Proposition~\ref{thm:isolated},
Corollary~\ref{thm:riemlorcausal}   and Remark~\ref{thm:analytic}).


\mysubsection{The differential problem}
\label{sub:differential}

Using a trivialization of the normal bundle along $\gamma$ by means
of a parallel moving frame, we now reformulate the Jacobi problem given by
\eqref{eq:eqJacobi} and \eqref{eq:PJacobi} in terms of a second order linear
differential equation in $\R^n$ with suitable initial conditions.

To this aim, we consider the following objects.
Let $g$ be a non degenerate symmetric bilinear form in $\R^n$, 
and let $R(t)$, $t\in[a,b]$, be a continuous curve in the space of linear
operators in $\R^n$ such that $R(t)$ is $g$-symmetric for all $t\in[a,b]$, i.e.,
$g(R(t)[v],w)=g(v,R(t)[w])$ for all $v,w\in\R^n$.

Let $P\subset\R^n$ be a subspace such that the restriction of $g$ to $P$
is non degenerate, and let $S$ be a symmetric bilinear form on $P$;
then, there exists a $g$-symmetric linear operator on $P$,
which we also denote by $S$, satisfying $S(v,w)=g(S[v],w)$ for all
$v,w\in P$. We will say that the family $(g,R,P,S)$ is an
{\em admissible quadruple for the differential problem\/} in $\R^n$.

We consider the following linear differential equation in $\R^n$:
\begin{equation}\label{eq:DE}
J''(t)=R(t)[J(t)],\quad t\in[a,b];
\end{equation}
we will consider solutions $J$  of \eqref{eq:DE} that satisfy
in addition the following initial conditions:
\begin{equation}\label{eq:IC}
J(a)\in P,\quad J'(a)+S[J(a)]\in P^\perp,
\end{equation}
where $P^\perp$ is the $g$-orthogonal complement of $P$ in $\R^n$;
such vector fields will be called {\em $(P,S)$-solutions}.
Note that, if $P=\{0\}$ (and thus $S=0$), a $(P,S)$-solution is simply a solution
of \eqref{eq:DE} vanishing at $t=a$.

Let $\Jts$ be the space of all $(P,S)$-solutions:
\begin{equation}\label{eq:defJts}
\Jts=\Big\{J: J\ \text{satisfies \eqref{eq:DE} and \eqref{eq:IC}}\Big\};
\end{equation}
and, for $t\in\,[a,b]$, we set $\Jts[t]=\big\{J(t):J\in\Jts\big\}$.\

Observe that ${\rm dim}(\Jts)=n$ since the subspace of $\R^n\oplus\R^n$ determined
by \eqref{eq:IC} is $n$-dimensional.
\begin{deffocal}\label{thm:deffocal}
An instant $t_0\in\,]a,b]$ is {\em $(P,S)$-focal\/} if there exists a non null
$(P,S)$-solution $J$ such that $J(t_0)=0$. The {\em multiplicity\/}
$\mul(t_0)$ of $t_0$ is the dimension of the subspace of $\Jts$ consisting
of such solutions; if $t_0$ is not $(P,S)$-focal we set $\mul(t_0)=0$:
\begin{equation}\label{eq:defmulbis}
\mul(t_0)={\rm dim}\, \big\{\Jt\in\Jts:\Jt(t_0)=0\big\} .
\end{equation} 
\end{deffocal}
Since ${\rm dim}(\Jts)$ is equal to $n$, then the multiplicity $\mul(t_0)$
is the codimension of $\Jts[t_0]$ in $\R^n$:
\begin{equation}\label{eq:dimcodim}
\mul(t_0)={\rm codim}(\Jts[t_0])={\rm dim}(\Jts[t_0]^\perp).
\end{equation}
In analogy with Definition~\ref{thm:defsiggeom}, we now give the following:
\begin{defsigdiff}\label{thm:defsigdiff}
The {\em signature\/} $\sgn(t_0)$ of the $(P,S)$-focal instant
$t_0$ is defined to be the signature of the restriction of $g$ to the
space $\Jts[t_0]^\perp$. If $t_0$ is not $(P,S)$-focal, we set
$\sgn(t_0)=0$; if the set of $(P,S)$-focal instants is finite,
we define the {\em focal index\/}  ${\rm i}_{\rm foc}$ of the
quadruple $(g,R,P,S)$ to be the sum of the signatures of the
$(P,S)$-focal instants:
\[{\rm i}_{\rm foc}=\sum_{t\in\,]a,b]}\sgn(t).\]
\end{defsigdiff}

\mysubsection{Relations between the geometrical and the differential problem}
\label{sub:geomdiff}

Suppose that an admissible quadruple  $(\M,\mathfrak  g, {\mathcal P},\gamma)$ for the
geometrical problem is given. For all
$t\in[a,b]$, the linear operator $v\mapsto {\mathcal R}(\gamma'(t),v)\,\gamma'(t)$
in $T_{\gamma(t)}\M$ is $\mathfrak  g$-symmetric, and, by the
usual symmetry properties of the curvature tensor, it takes values in 
$\gamma'(t)^\perp$; we consider its restriction to $\gamma'(t)^\perp$
and we denote it by ${\mathcal R}(t)$.

If we choose an arbitrary parallel moving frame that trivializes the normal
bundle $(\gamma')^\perp$ along $\gamma$, so that we have an isomorphism
between $\gamma'(t)^\perp$ and $\R^{m-1}$, we get a linear
operator $R(t)$ on $\R^{m-1}$ corresponding to ${\mathcal R}(t)$,
and a symmetric bilinear form $g$ corresponding to (the restriction to 
$\gamma'(t)^\perp$ of) $\mathfrak  g$. Since $\mathfrak  g$ is parallel,
then $g$ is constant; obviously, $R(t)$ is $g$-symmetric. Similarly,
the subspace $T_{\gamma(a)}{\mathcal P}\subset\gamma'(a)^\perp$ corresponds
to a subspace $P$ of $\R^{m-1}$, and the second fundamental
form ${\mathcal S}_{\gamma'(a)}$ corresponds to a symmetric bilinear
form $S$ on $P$; moreover, the restriction of $g$ to $P$ is non
degenerate.

If $\gamma$ is not lightlike, then the restriction of $\mathfrak  g $ to $\gamma'(t)^\perp$
is non degenerate for all $t$, which implies that $g$ is non degenerate in
$\R^{m-1}$ so that $(g,R,P,S)$ is an admissible quadruple for the differential
problem in $\R^{m-1}$. We will say that $(g,R,P,S)$ is  {\em associated\/} to the quadruple
$(\M,\mathfrak  g,\gamma, {\mathcal P})$ by the choice 
of a parallel trivialization of the normal bundle along $\gamma$.

If $(g,R,P,S)$ is associated to an admissible quadruple for the geometrical
problem $(\M,\mathfrak  g,\gamma, {\mathcal P})$, then $R$ is indeed a smooth
map. Conversely, every quadruple $(g,R,P,S)$  with $R$ smooth
arises in this way:

\begin{correspondence}\label{thm:correspondence}
If $(g,R,P,S)$ is an admissible quadruple for the differential problem
in $\R^n$,
with $R$ smooth, then there exists an admissible quadruple
$(\M,\mathfrak  g,\gamma, {\mathcal P})$ for the geometrical problem 
such that $(g,R,P,S)$ is associated to $(\M,\mathfrak  g,\gamma, {\mathcal P})$
by some choice of a parallel trivialization of the normal bundle along $\gamma$.
Moreover, the quadruple $(\M,\mathfrak  g,\gamma, {\mathcal P})$ can be chosen
with $\gamma$ timelike as well as spacelike, and $\Mg$ can be chosen to be
conformally flat. If $g$ is positive definite, then $\Mg$ is
Riemannian  if $\gamma$ is spacelike and Lorentzian if $\gamma$ is timelike.
\end{correspondence}
\begin{proof}
Consider $\M=\R^{n+1}$ with coordinates $(x_1,x_2,\ldots,x_{n+1})$
and canonical basis $\{e_1,e_2,\ldots,e_{n+1}\}$; 
let $\gamma:~[a,b]\mapsto\M$ be the curve
$\gamma(t)=t\cdot e_{n+1}$. We consider the non degenerate symmetric bilinear form $\mathfrak  g_0$ on
$\M$ given by $\mathfrak  g_0(e_i,e_j)=g(e_i,e_j)$ for $i,j=1,\ldots,n$,
$\mathfrak  g_0(e_{n+1},e_{n+1})=\pm1$, and $\mathfrak  g_0(e_i,e_j)=0$ otherwise.

The choice of the sign of $\mathfrak  g_0(e_{n+1},e_{n+1})$ is done according
to whether $\gamma$ should be timelike or spacelike, as desired.

Let $\M$ be endowed with the conformally flat metric $\mathfrak  g=e^\Omega\cdot\mathfrak  g_0$,
where $\Omega$ is a smooth function in $\R^{n+1}$ that vanishes together  with
its partial derivatives on the $e_{n+1}$-axis. 
The factor $\Omega$ will be chosen so that the
corresponding metric $\mathfrak  g$ will satisfy the required properties.

To this goal, we recall some formulas about the covariant derivative and
the geodesic equation in general conformal metrics. Let $\nabla^{(0)}$
and $\nabla$ denote the covariant derivative or the gradient operators
in the metrics $\mathfrak  g_0$ and $\mathfrak  g$ respectively;
note that the covariant derivative $\nabla^{(0)}$ is the usual directional 
derivative in $\R^{n+1}$, although the gradient $\nabla^{(0)}$ is
{\em not\/} the usual gradient in $\R^{n+1}$. 

For smooth vector fields
$X,Y$ in $\M$, we have:
\begin{equation}\label{eq:dercovconf}
\nabla_XY=\frac12\Big[\mathfrak  g_0(\nabla^{(0)}\Omega,X)\,Y+
\mathfrak  g_0(\nabla^{(0)}\Omega,Y)\,X-\mathfrak  g_0(X,Y)\,\nabla^{(0)}\Omega\Big]+\nabla^{(0)}_XY;
\end{equation} 
moreover, the geodesic equation in $\Mg$ is:
\begin{equation}\label{eq:geoeqconf}
\nabla^{(0)}_{ c'} c'=\frac12\, \mathfrak  g_0(  c',c')\,\nabla^{(0)}\Omega
-\mathfrak  g_0(\nabla^{(0)}\Omega,c')\,c'.
\end{equation}
Since $\nabla^{(0)}\Omega\equiv0$ on $\gamma$, then $\gamma$ is a geodesic in $\Mg$;
moreover, by \eqref{eq:dercovconf}, the parallel vector fields along $\gamma$
in $\Mg$ are just the constant vector fields. Hence, we trivialize the normal
bundle along $\gamma$ in $\Mg$ by choosing the first $n$ vectors of the canonical
basis as a parallel moving frame.

To compute the Jacobi equation along $\gamma$ in $\Mg$, we linearize the geodesic 
equation \eqref{eq:geoeqconf}, obtaining:
\begin{equation}\label{eq:Jacconf}
J''=\frac12\,\mathfrak  g_0(\gamma',\gamma')\,{\rm Hess}_\Omega^{(0)}(J),
\end{equation}
where $J''$ is the ordinary second derivative in $\R^{n+1}$ and ${\rm Hess}^{(0)}_\Omega$
is the $\mathfrak  g_0$-sym\-metric linear operator given by ${\rm Hess}^{(0)}_\Omega (v)=
\nabla_v^{(0)}\,\nabla^{(0)}\Omega$. In the deduction of \eqref{eq:Jacconf}
we have used the fact that $\nabla^{(0)}\Omega$ and ${\rm Hess}_\Omega^{(0)}(\gamma')$
vanish on $\gamma$.

Observing that the covariant derivative along $\gamma$ in $\Mg$ equals ordinary derivative
in $\R^{n+1}$ and comparing equation \eqref{eq:Jacconf} with the general Jacobi equation
\eqref{eq:eqJacobi} we see that the curvature tensor $\mathcal R$ of $\Mg$ along $\gamma$ is
given by:
\[\mathcal R(\gamma',v)\,\gamma'=
\frac12\,\mathfrak  g_0(\gamma',\gamma')\,{\rm Hess}_\Omega^{(0)}(v).\]
It is easily checked that:
\[\mathfrak  g_0({\rm Hess}_\Omega^{(0)}(e_i),e_j)=\frac{\partial^2\Omega}{\partial x_i\partial x_j},
\quad\forall\,i,j=1,\ldots, n+1;\]
if we set:
\[a_{ij}(t)=\frac 2{\mathfrak  g_0(\gamma',\gamma')}\,\mathfrak  g_0(R(t)\,e_i,e_j),
\quad i,j=1,\ldots,n,\ t\in[a,b]\]
and consider an arbitrary smooth extension of $a_{ij}$ on $\R$, then the assignment
\[\Omega(x_1,\ldots,x_{n+1})=\frac12\sum_{i,j=1}^n a_{ij}(x_{n+1})\,x_i\,x_j\]
gives the required function.

To conclude the proof, we now need to exhibit a submanifold $\mathcal P$
of $\M$, passing through $\gamma(a)$ with tangent space $T_{\gamma(a)}{\mathcal P}
=P\oplus\{0\}$, and whose second fundamental form in the normal
direction $\gamma'(a)$ equals $S$. This will follow immediately from the next Lemma,
in which we prove something slightly more general.

The last assertion in the statement of the proposition is totally obvious.
\end{proof}
\begin{slightly}\label{thm:slightly}
Let $\Mg$ be a semi-Riemannian manifold, $p\in\M$,  $P$ a non degenerate subspace
of $T_p\M$ and  $S:P\times P\mapsto P^\perp$ be a symmetric bilinear map.
Then, there exists a smooth submanifold $\mathcal P$ of $\M$, with $p\in{\mathcal P}$,
such that $T_p{\mathcal P}=P$ and such that the second fundamental form
${\mathcal S}$ of ${\mathcal P}$ at $p$ equals $S$.
\end{slightly}
\begin{proof}
Let $U_0\subset T_p\M$ be an open neighborhood of the origin such that
the exponential map $\exp_p$ of $\Mg$ maps $U_0$ diffeomorphically onto
an open neighborhood of $p$ in $\M$. Regarding $\exp_p$ as a coordinate
map around $p$, it is well known that the Christoffel symbols of the
Levi--Civita connection vanish at the point $0$. Hence, the covariant derivative
at this point coincide with the usual directional derivative in $T_p\M$.
If ${\mathcal P}_0$ is a submanifold of $U_0$ passing through $0$
and ${\mathcal P}=\exp_p({\mathcal P}_0)$, then, since ${\rm d}\exp_p(0)$ is the
identity map, the tangent space $T_p{\mathcal P}$ is $T_{\scriptscriptstyle0}{\mathcal P}_0$;
moreover, by the above observation about the covariant derivative,
the second fundamental form of $\mathcal P$ at $p$
equals the second fundamental form of ${\mathcal P}_0$ at $0$ in the flat
space $T_p\M$.

We define ${\mathcal P}_0$ to be the smooth submanifold of $T_p\M$ 
given by the graph of the map $v\mapsto \frac12\, S(v,v)$ in the decomposition
$P\oplus P^\perp$, namely:
\[{\mathcal P}_0=\Big\{v+{\textstyle{\frac12}}\, S(v,v):v\in P\Big\}\cap U_0.\]
The conclusion follows from an elementary calculation of the second fundamental
form of ${\mathcal P}_0$.
\end{proof}
So far, we have associated an admissible quadruple for the differential
problem only to quadruples $(\M,\mathfrak  g,\gamma, {\mathcal P})$ with
$\gamma$ spacelike or timelike. Indeed, if $\gamma$ is lightlike, then the
symmetric bilinear form $g$ previously defined is degenerate on $\R^{m-1}$.
One way to avoid this problem, following a customary procedure in Morse
Theory (see~\cite{BEE} for the Lorentzian case), is to consider a suitable
quotient of the normal bundle along the lightlike geodesic $\gamma$.

More precisely, given an admissible quadruple $(\M,\mathfrak  g,\gamma, {\mathcal P})$
with $\gamma$ lightlike, for all  $t\in[a,b]$ we consider the quotient space
${\mathcal N}(t)=\gamma'(t)^\perp/\big[\R\,\gamma'(t)\big]$, where $\R\,\gamma'(t)$
is the one dimensional subspace generated by $\gamma'(t)$.
It is easy to see that ${\mathcal N}=\bigcup_{t\in[a,b]}{\mathcal N}(t)$ 
is a vector bundle along $\gamma$.

Since the kernel of the restriction of the metric $\mathfrak  g$ to $\gamma'(t)^\perp$
is precisely $\R\,\gamma'(t)$, then $\mathfrak  g$ gives a well defined non degenerate
symmetric bilinear form $\overline{\mathfrak  g}$ on the quotient space ${\mathcal N}(t)$.
Similarly, the linear operator ${\mathcal R}(t)={\mathcal R}(\gamma',\cdot)\,\gamma'$
annihilates $\R\,\gamma'(t)$, and therefore it defines a linear operator
$\overline{\mathcal R}(t)$ on ${\mathcal N}(t)$. Obviously, $\overline{\mathcal R}$
is $\overline{\mathfrak  g}$-symmetric.

The subspace $T_{\gamma(a)}{\mathcal P}\subset\gamma'(a)^\perp$ does not
contain $\gamma'(a)$, because of our nondegeneracy assumption on 
$T_{\gamma(a)}{\mathcal P}$, hence it may be identified with a subspace
of ${\mathcal N}(a)$, which will be denoted by the same symbol.

Let us now consider a trivialization of the normal bundle 
along $\gamma$ by a parallel moving frame in such a way that
the last vector field of the frame is the tangent vector $\gamma'$
itself. The remaining $m-2$ vector fields define a moving frame
in the bundle $\mathcal N$, and they induce a trivialization of $\mathcal N$.
We therefore get tensors
$R(t)$ and $g$ on $\R^{m-2}$ corresponding to the tensors $\overline{\mathcal R}(t)$
and $\overline{\mathfrak  g}$, as well as a subspace $P\subset \R^{m-2}$ and
a symmetric bilinear form $S:P\times P\mapsto\R$ corresponding to the subspace
$T_{\gamma(a)}{\mathcal P}$ of ${\mathcal N}(a)$ and the second fundamental
form ${\mathcal S}_{\gamma'(a)}$ of ${\mathcal P}$, respectively.

We have thus constructed an admissible quadruple $(g,R,P,S)$ for the
differential problem in $\R^{m-2}$ which we call {\em the associated
quadruple\/} to $(\M,\mathfrak  g,\gamma, {\mathcal P})$ in the case of
a lightlike geodesic $\gamma$.
\begin{remlight}\label{thm:remlight}
Given an admissible quadruple $(\M,\mathfrak  g,\gamma, {\mathcal P})$
for the geometric problem  
and an associated quadruple $(g,R,P,S)$ corresponding 
to some parallel trivialization
of the normal bundle along $\gamma$, we introduce a linear map $\Phi$ 
that carries vector fields orthogonal to $\gamma$ into vector fields
in $\R^n$, as follows. If $\gamma$ is non lightlike, $\Phi(v)$ is simply
the set of $m-1$ coordinates  of $v$ with respect to the chosen parallel
moving frame. When $\gamma$ is lightlike, $\Phi(v)$ is the set of $m-2$
coordinates of the projection of $v$ in ${\mathcal N}$ with respect
to the chosen parallel basis of $\mathcal N$. 

If $\gamma$ is non lightlike, such a map $\Phi$ gives an isomorphism between
$\Jms^\perp$ and $\Jts$; if $\gamma$ is lightlike, $\Phi$ maps $\Jms^\perp$
{\em onto\/} $\Jts$, and its kernel consists of affine multiples of $\gamma'$.
The surjectivity of $\Phi$ in the lightlike case follows by observing that
if $\Jm$ is a solution  of 
\[\Jm''={\mathcal R}(\gamma',\Jm)\,\gamma'+f\,\gamma'\]
for some fixed smooth map $f:[a,b]\mapsto\R$, then 
$\Jm-F\,\gamma'$ is a Jacobi field along $\gamma$, where $F''=f$.
\end{remlight}

The relation between the focal indexes of the geometric and differential
problems is clarified by the following:

\begin{matatudo}\label{thm:matatudo}
Let $(\M,\mathfrak  g,\gamma, {\mathcal P})$ be an admissible quadruple
for the geometric problem and $(g,R,P,S)$ be an associated quadruple
corresponding to some parallel trivialization of the normal bundle along $\gamma$.
Then, for all $t_0\in\,]a,b]$ there exists an isomorphism between
$\Jms[t_0]^\perp$ and $\Jts[t_0]^\perp$ which carries the restriction
of $\mathfrak  g$ to the restriction of $g$. In particular, for $t_0\in\,]a,b]$,
$\gamma(t_0)$ is a $\mathcal P$-focal point if and only if $t_0$ is a 
$(P,S)$-focal instant.
In this case, its multiplicity and signature in the geometric and in the differential
problem coincide, from which it follows that the focal indexes of the problems
are equal. 
\end{matatudo}
\begin{proof}
If $\gamma$ is not lightlike, let $\phi:\gamma'(t_0)^\perp\mapsto\R^{m-1}$ 
be the isomorphism given by the chosen trivialization of the normal bundle to $\gamma$.
For a lightlike $\gamma$, let's denote by $\phi:{\mathcal N}(t_0)\mapsto\R^{m-2}$
the isomorphism determined by the choice of the trivialization of the quotient bundle,
as described above. By construction, $\phi$ carries $\mathfrak  g$
(or $\overline{\mathfrak  g}$ for $\gamma$ lightlike) to $g$.

\noindent\ 
For $\gamma$ non lightlike, it is easily checked that $\phi$ carries $\Jms^\perp[t_0]$
onto $\Jts[t_0]$ by observing the correspondence between $\mathcal P$-Jacobi
fields  orthogonal to $\gamma$ and $(P,S)$-solutions of \eqref{eq:DE}.
Similarly, if $\gamma$ is lightlike, $\phi$ carries the quotient
$\Jms^\perp[t_0]/\big[\R\,\gamma'(t_0)\big]$ onto $\Jts[t_0]$
(see Remark~\ref{thm:remlight}). 

For $\gamma$ non lightlike, taking the orthogonal complements of
$\Jms^\perp[t_0]$ in $\gamma'(t_0)^\perp$ and of
$\Jts[t_0]$ in $\R^{m-1}$, using \eqref{eq:relJmsperpt} we conclude that 
$\phi$ induces the desired isomorphism between $\Jms[t_0]^\perp$ and
$\Jts[t_0]^\perp$.

If $\gamma$ is lightlike, we take the orthogonal complements of $\Jms^\perp[t_0]
/\big[\R\,\gamma'(t_0)\big]$ in ${\mathcal N}(t_0)$ and of
$\Jts[t_0]$ in $\R^{m-2}$. Again, using \eqref{eq:relJmsperpt}
we get that $\phi$ induces an isomorphism between the image
of $\Jms[t_0]^\perp$ in the quotient space ${\mathcal N}(t_0)$ and
$\Jts[t_0]^\perp$. To conclude the proof, we observe that \eqref{eq:codim1}
implies that $\gamma'(t_0)$ does not belong to $\Jms[t_0]^\perp$, which
implies that it maps isomorphically into ${\mathcal N}(t_0)$.
\end{proof}
\begin{unifying}\label{thm:unifying}
The main feature of the method we have described for associating
quadruples $(g,R,P,S)$ to quadruples $(\M,\mathfrak  g,\gamma, {\mathcal P})$
consists in the fact that, in the Lorentzian case, 
if $\gamma$ is non spacelike, then the bilinear form $g$ is positive definite.
In particular, from Proposition~\ref{thm:matatudo} it follows that,
if $\Mg$ is Lorentzian and $\gamma$ is non spacelike,
then $\mathfrak  g$ is positive definite in $\Jms[t_0]^\perp$.
The positivity of $g$ in $\R^n$ will be used later (see Remark~\ref{thm:remMorse})
to derive an alternative proof of the Morse index theorem for non spacelike Lorentzian 
geodesics.

\noindent
A different way of associating a quadruple $(g,R,P,S)$
to a quadruple  $(\M,\mathfrak  g,\gamma, {\mathcal P})$ is 
to consider a trivialization of the {\em entire\/} tangent bundle
along $\gamma$. Using this approach, one unifies the construction
for the lightlike and the non lightlike case; we will need this construction
in Section~\ref{sec:applications}, where we will discuss a problem of
approximation of lightlike geodesics by timelike geodesics. 
This construction will not introduce substantial modifications
of the solution space for the differential problem. Namely,
the statement of Proposition~\ref{thm:matatudo} is trivial
if the association of quadruple is understood in this sense.
It follows that, by using the two different associations of quadruples,
we get the same $(P,S)$-focal instants, with the same multiplicities and 
signatures, and thus the same focal index.

We remark also that the statement and the proof of 
Proposition~\ref{thm:correspondence} can be adapted to the case that the 
association of quadruples is made by trivializing the whole
tangent bundle. 
\end{unifying}

\mysubsection{The symplectic structure associated to a differential problem}
\label{sub:symplstruct}
Given the perfect analogy between the geometrical and the differential
problem, as given by Proposition~\ref{thm:correspondence} and 
Proposition~\ref{thm:matatudo}, we will henceforth concentrate our
attention on an admissible quadruple for the differential
problem $(g,R,P,S)$ in $\R^n$. 

Given two solutions $J_1$ and $J_2$ of the differential equation
\eqref{eq:DE}, the quantity
\begin{equation}\label{eq:const}
\sigma(t)=g(J_1(t),J_2'(t))-g(J_1'(t),J_2(t))
\end{equation}
is constant in $[a,b]$. Namely, a straightforward calculation
using equation \eqref{eq:DE} shows that $\sigma'$ 
vanishes identically. This motivates the following definition:
\begin{defomega}\label{thm:defomega}
The symplectic form $\omega$ on $\R^{2n}$ associated to $g$ is
given by:
\[\omega[(x_1,x_2),(y_1,y_2)]=g(x_1,y_2)-g(x_2,y_1).\]
The nondegeneracy of $\omega$ follows easily from the nondegeneracy of $g$.
\end{defomega}
The initial conditions $(J(a),J'(a))\in\R^{2n}$ determine
uniquely a solution of \eqref{eq:DE}, therefore
the space of solutions of \eqref{eq:DE} can be identified with
$\R^{2n}$. For all $t\in[a,b]$, we have a linear automorphism $\Psi(t)$
of $\R^{2n}$ satisfying
\begin{equation}\label{eq:defPsit}
\Psi(t)[(J(a),J'(a))]=(J(t),J'(t)),
\end{equation}
for every solution $J$ of \eqref{eq:DE}. This automorphisms are implemented
by what is usually called the {\em fundamental matrix\/} of the first order
linear differential system associated to \eqref{eq:DE}. Observe that
$t\mapsto\Psi(t)$ is a curve of class $C^1$ in the general linear group
${\rm GL}(2n,\R)$ which satisfies $\Psi(0)={\rm Id}$.

Using the fact that the quantity \eqref{eq:const} is constant,
it is also easy to observe that $\Psi(t)$ preserves the symplectic form
$\omega$ for all $t$:
\[\omega[\Psi(t)\,x,\Psi(t)\,y]=\omega[x,y],\quad\forall\,x,y\in\R^{2n},\]
hence, $\Psi(t)$ is a curve in the symplectic group of $\R^{2n}$ corresponding
to $\omega$.

The important observation here is that $\omega$ vanishes on the $n$-dimensional
subspace of $\R^{2n}$ determined by the initial conditions \eqref{eq:IC}.
Namely, if $J_1,J_2\in\Jts$, then $J_i(a)\in P$ and $J_i'(a)+S[J_i(a)]\in P^\perp$
for $i=1,2$, and:
\begin{equation}\label{eq:ell0Lagr}
\begin{split}
\omega[(J_1(a),J_1'(a))&,(J_2(a),J_2'(a))]=\\&= g(J_1(a),J_2'(a))-g(J_1'(a),J_2(a))=\\
&= g(J_1(a),-S[J_2(a)])-g(-S[J_1(a)],J_2(a))=0,
\end{split}
\end{equation}
where the last equality follows from the $g$-symmetry of $S$.

Summarizing the facts that \eqref{eq:const} is constant and that $\omega$ vanishes
on the space of initial conditions of $(P,S)$-solutions, we have the following
identity:
\begin{equation}\label{eq:identity}
g(J_1'(t),J_2(t))=g(J_1(t),J_2'(t)),\quad\forall\, J_1,J_2\in\Jts,
\end{equation}
for all $t\in[a,b]$.
\medskip

\mysubsection{On the discreteness of the set of $(P,S)$-focal instants}
\label{sub:isolatedness}
We give some conditions that guarantee the discreteness of the set of $(P,S)$-focal instants.
\begin{isolated}\label{thm:isolated}
Let $(g,R,P,S)$ be an admissible quadruple for the differential problem
in $\R^n$, 
and let $t_0$ be a $(P,S)$-focal instant. If $g$ is non degenerate  on $\Jts[t_0]$,
then there are no $(P,S)$-focal instants other than $t_0$ in some neighborhood
of $t_0$. Moreover, there are no $(P,S)$-focal instants in some neighborhood
of the initial instant $a$.
\end{isolated}
\begin{proof}
Let $\mul(t_0)=n-k>0$ be the multiplicity of the focal instant $t_0$.
Let $J_1,J_2,\ldots,J_n$ be a basis of $\Jts$ such that 
$J_1(t_0),\ldots,J_k(t_0)$ are a basis for $\Jts[t_0]$
and $J_{i}(t_0)=0$ for $i\ge k+1$. 

The vectors $J'_{k+1}(t_0),\ldots,J'_n(t_0)$ are a basis of $\Jts[t_0]^\perp$.
To prove this, we first observe that they belong to $\Jts[t_0]^\perp$; namely,
by \eqref{eq:identity}, if $i\in\{k+1,\ldots,n\}$ and $j\in\{1,\ldots,k\}$, we have
\[g(J_i'(t_0),J_j(t_0))=g(J_i(t_0),J_j'(t_0))=g(0,J'_j(t_0))=0.\]
To prove the claim, we need to show that the vectors
$J'_{k+1}(t_0),\ldots,J'_n(t_0)$  are linearly independent,
because ${\rm dim}(\Jts[t_0]^\perp)=n-k$, by \eqref{eq:dimcodim}. 
To see this, observe that the fields
$J_{k+1},\ldots,J_n$ are linearly independent in $\Jts$, hence the pairs
\[(J_{k+1}(t_0),J_{k+1}'(t_0)),\ldots,(J_n(t_0),J_n'(t_0))\] are linearly
independent in $\R^{2n}$. The conclusion follows from the fact that
$J_{k+1}(t_0)=\ldots=J_n(t_0)=0$.

We now define a family of {\em continuous\/} vector fields $\tilde J_1,\ldots,\tilde J_n$
along $\gamma$, by setting:
\[\tilde J_j=J_j,\quad\text{for}\ j=1,\ldots,k;\]
and
\[\tilde J_i(t)=\left\{\begin{array}{lr}
{\displaystyle\frac{J_i(t)}{t-t_0}},&\text{if}\ t\ne t_0,\\ \\
J'_i(t_0),&\text{if}\ t=t_0,
\end{array}\right.\qquad\text{for}\
i=k+1,\ldots,n.\]
The vectors $\tilde J_1(t_0),\ldots,\tilde J_n(t_0)$ are now a basis for $\R^n$.

Namely, the first $k$ vectors $\tilde J_1(t_0),\ldots,\tilde J_k(t_0)$ are
a basis for $\Jts[t_0]$, and the remaining $n-k$ vectors
$\tilde J_{k+1}(t_0),\ldots, \tilde J_n(t_0)$ are a basis for $\Jts[t_0]^\perp$;
moreover, $g$ is non degenerate on $\Jts[t_0]$, which implies that $\R^n=\Jts[t_0]\oplus
\Jts[t_0]^\perp$.

By continuity, the vectors $\tilde J_1(t),\ldots,\tilde J_n(t)$ are a basis
for $\R^n$ for $t$ sufficiently close to $t_0$. But that implies that, for $t$
sufficiently close to $t_0$ and $t\ne t_0$ the vectors $J_1(t),\ldots, J_n(t)$
are a basis for $\R^n$, which implies that there are no $(P,S)$-focal instants
around $t_0$.

The case $t_0=a$ is treated similarly, observing that $\Jts[a]=P$ and considering that
$g$ is non degenerate on $P$.
\end{proof}
We have the following immediate Corollary:
\begin{riemlorcausal}\label{thm:riemlorcausal}
Let $(\M,\mathfrak  g,\gamma, {\mathcal P})$ be an admissible
quadruple for the geometric problem.
 Assume $\Mg$ is Riemannian or Lorentzian, and in the latter
case, that $\gamma$ is non spacelike. Then, there are only a finite number
of $\mathcal P$-conjugate points along $\gamma$. Moreover, the
focal index and the geometrical index of $\gamma$ coincide:
\begin{equation}\label{eq:eqindexes}
{\rm i}_{\rm geom}(\gamma)={\rm i}_{\rm foc}(\gamma).
\end{equation}
\end{riemlorcausal}
\begin{proof}
It is an easy consequence of Remark~\ref{thm:unifying}, Proposition~\ref{thm:matatudo} 
and Proposition~\ref{thm:isolated}.
\end{proof}
\begin{analytic}\label{thm:analytic}
The $(P,S)$-focal instants coincide precisely with the
zeroes of the function $r(t)={\rm det}(J_1(t),\ldots,J_n(t))$, where
$J_1,\ldots,J_n$ is a basis of $\Jts$.
If $(g,R,P,S)$ is an admissible quadruple with $R(t)$ real analytic
on $[a,b]$, then $r(t)$ is also analytic, and so its zeros are isolated.
Observe indeed that $r(t)$ cannot vanish identically on $[a,b]$ because,
by Proposition~\ref{thm:isolated}, $r(t)$ is non zero for $t$ sufficiently
close to  $a$, $t\ne a$.
It follows easily that, if $(\M,\mathfrak  g,\gamma, {\mathcal P})$
is an admissible quadruple for the geometric problem with
$\Mg$ analytic, then the set of 
$\mathcal P$-focal points along $\gamma$ is finite.
\end{analytic}

\end{section}
\begin{section}{Geometry of the Lagrangian Grassmannian}\label{sec:geometry}
We have seen in Section~\ref{sec:preliminaries}
that the set $\Jts$ can be identified with a {\em Lagrangian\/}
subspace of the symplectic space $(\R^{2n},\omega)$, i.e., a maximal
subspace of $\R^{2n}$ on which $\omega$ vanishes. In view to future applications,
in this Section we present the main properties and we discuss the
geometrical structure of the collection of all Lagrangian subspaces of 
a symplectic space.
\smallskip

Throughout this section we will assume that $V$ is a $2n$-dimensional real
vector space, equipped with a symplectic form $\omega$, i.e., a skew symmetric
non degenerate bilinear form on $V$. 
\mysubsection{Generalities on symplectic spaces}
\label{sub:generalities}
A {\em symplectic basis\/} of
$(V,\omega)$ is a vector space basis $e_1,\ldots,e_{2n}$ of $V$ such
that \[\omega[e_{n+j},e_j]=-\omega[e_j,e_{n+j}]=1\] for all $j=1,\ldots,n$,
and $\omega[e_i,e_j]=0$ otherwise; the existence of a symplectic basis
in $(V,\omega)$ is standard. We recall that a {\em complex structure\/} for
$V$ is  a linear operator ${\mathcal I}:V\mapsto V$ such that ${\mathcal I}^2=-{\rm Id}$.
A complex structure $\mathcal I$ on $V$ induces a complex vector space structure
on $V$, and $\mathcal I$ becomes the scalar multiplication by the imaginary unit $i$.
A complex structure $\mathcal I$ is {\em compatible\/} with the symplectic form
$\omega$ if the bilinear form $\omega[{\mathcal I}\cdot,\cdot]$ is symmetric
and positive definite on $V$. 

If $(V_i,\omega_i)$, $i=1,2$, are
symplectic spaces of the same dimension $2n$, 
a linear map $T:V_1\mapsto V_2$ is called a 
{\em symplectomorphism\/}
if $\omega_2(Tx,Ty)=\omega_1(x,y)$ for all $x,y\in V_1$. Observe that
a symplectomorphism $T$ is always an isomorphism; namely, the $n$-th exterior 
powers $\omega_i^{n}$ are {\em volume forms\/} in $V_i$, $i=1,2$, which
are preserved by $T$.

We identify $\R^{2n}$ with $\C^n$ by considering the first $n$ coordinates to
be the real part, and the remaining coordinates to be the imaginary part.
Therefore, we get a complex structure ${\mathcal I}_0$ given by 
${\mathcal I}_0(e_j)=e_{n+j}$, ${\mathcal I}_0(e_{n+j})=-e_j$, for
$j=1,\ldots,n$, where $\{e_i\}_{i=1}^{2n}$ is the canonical basis of $\R^{2n}$.
For $x,y\in\R^{2n}$, we denote by $x\cdot y$ the Euclidean inner product, and
by $\langle x,y\rangle$ the Hermitian product in $\C^n\simeq\R^{2n}$
whose real part is $x\cdot y$ and which is conjugate in the second variable.
The {\em canonical\/} symplectic form $\omega_0$ in $\R^{2n}$ is the imaginary
part of the Hermitian product. Observe that the canonical basis is a symplectic
basis for $\omega_0$ and ${\mathcal I}_0$ is compatible with $\omega_0$.
\smallskip

A subspace $W$ of $V$ will be called {\em isotropic\/} if $\omega$ vanishes 
identically on $W$ (by this we mean on $W\times W$); an $n$-dimensional isotropic 
subspace $W$ will
be called a {\em Lagrangian subspace\/} of $(V,\omega)$. It is easy to see that
the Lagrangian subspaces coincide with the {\em maximal\/}
isotropic subspaces of $(V,\omega)$. 

Given a Lagrangian direct sum decomposition $V=L_0\oplus L_1$, i.e., both
subspaces $L_0$ and $L_1$ are Lagrangian, we denote by $\Iota{L_0}{L_1}$
the isomorphism from $L_1$ to the dual space $L_0^*$ given by:
\begin{equation}\label{eq:defIota}
\Iota{L_0}{L_1}(v)=\omega[v,\cdot]\lower5truept\hbox{$\Big\vert$}_{L_0},\quad\forall\,v\in L_1.
\end{equation}
The injectivity of $\Iota{L_0}{L_1}$ follows immediately from the non degeneracy
of $\omega$. We observe that, by the anti-symmetry of $\omega$, the following
identity holds:
\begin{equation}\label{eq:Iotas}
\Iota{L_1}{L_0}=-(\Iota{L_0}{L_1})^*.
\end{equation}

\begin{remcomplLagr}\label{thm:remcomplLagr}
The existence of a complex structure  compatible with $(V,\omega)$ is
proven easily. Namely, a complex structure compatible with $(V,\omega)$ is 
obtained as the pull-back of
${\mathcal I}_0$ by the symplectomorphism $V\mapsto\R^{2n}$ defined by 
a symplectic basis of $(V,\omega)$. Using a compatible complex structure
${\mathcal I}$, we can now prove that every Lagrangian subspace
$L_0$ of $V$ admits a complementary Lagrangian subspace $L_1$. 
Namely, just define $L_1={\mathcal I}(L_0)$. Given any Lagrangian
direct sum decomposition $V=L_0\oplus L_1$, we construct a symplectic
basis $\{e_1,\ldots,e_n,f_1,\ldots,f_n\}$ of $V$ by taking any linear 
basis $\{e_1,\ldots,e_n\}$ of $L_0$ and the linear basis $\{f_1,\ldots,f_n\}$
of $L_1$ whose image by $\Iota{L_0}{L_1}$ is the dual basis of $\{e_1,\ldots,e_n\}$.
This implies that every linear isomorphism $\psi:L_0\mapsto\R^n\oplus\{0\}$
extends to a symplectomorphism $\psi$ from $(V,\omega)$ to $(\R^{2n},\omega_0)$
which carries $L_1$ to $\{0\}\oplus\R^n$.
\end{remcomplLagr}

The {\em symplectic group\/} ${\rm Sp}(V,\omega)$ is the Lie subgroup of ${\rm GL}(V)$
consisting of symplectomorphisms of $(V,\omega)$; its Lie algebra ${\rm sp}(V,\omega)$
consists of all linear maps $H:V\mapsto V$ such that:
\begin{equation}\label{eq:omegaantisym}
\omega(Hx,y)+\omega(x,Hy)=0,\quad\forall\;x,y\in V.
\end{equation} 
Equation \eqref{eq:omegaantisym} is equivalent to the symmetry of the bilinear form
$\omega(H\cdot,\cdot)$ on $V$.
 
The group ${\rm Sp}(\R^{2n},\omega_0)$
is also denoted by ${\rm Sp}(n,\R)$; the subgroup of $GL({2n},\R)$ consisting
of unitary transformations with respect to the canonical Hermitian product
is denoted by ${\rm U}(n)$. Since $\omega_0$ is the imaginary part of the
Hermitian product which is preserved by elements in ${\rm U}(n)$, we see that
${\rm U}(n)$ is a subgroup of ${\rm Sp}(n,\R)$.

By ${\rm O}(n)$ we mean the orthogonal group in $\R^n$, and
by ${\rm SO}(n)$ the subgroup of ${\rm O}(n)$ consisting
of matrices with determinant equal to $1$. Every linear map
$\psi:\R^n\mapsto\R^n$ has a unique $\C$-linear extension to a map $\psi^\C:
\C^n\mapsto\C^n$. If $\psi\in {\rm O}(n)$, then $\psi^\C\in {\rm U}(n)$, which
identifies ${\rm O}(n)$ with the subgroup of ${\rm U}(n)$ consisting of those
maps that preserve the subspace $\R^n\oplus\{0\}$ in $\R^{2n}$.

It is well known that ${\rm U}(n)$, ${\rm O}(n)$ and ${\rm SO}(n)$
are compact Lie groups,
and ${\rm Sp}(n,\R)$, or more in general ${\rm Sp}(V,\omega)$, is a non compact
Lie group.

\mysubsection{The Lagrangian Grassmannian}
\label{sub:lagrangian}

For $k=0,\ldots,2n$, we denote by $G_k(V)$ the Grassmannian of all the $k$-dimen\-sional
subspaces of $V$. We will be interested in the subset $\Lambda=\Lambda(V,\omega)\subset G_n(V)$
consisting of all the Lagrangian subspaces of $(V,\omega)$:
\[\Lambda=\Lambda(V,\omega)=\Big\{L: L\ \text{is a Lagrangian subspace of}\ (V,\omega)\Big\}.\]
For simplicity, we will omit the argument $(V,\omega)$ whenever there is
no risk of confusion, and we will write simply $\Lambda$.

We recall that $G_k(V)$ has the structure of a real analytic manifold of dimension
$k(2n-k)$; given a direct sum decomposition $V=W_0\oplus W_1$, where
${\rm dim}(W_0)=k$, a local chart of $G_k(V)$ is defined in an open neighborhood
of $W_0$ taking values in the vector space ${\mathcal L}(W_0,W_1)$ of all
linear maps $T:W_0\mapsto W_1$. Namely, to every $W\in G_k(V)$ which is {\em transversal\/}
to $W_1$, i.e., $W\cap W_1=\{0\}$, one associates the unique $T\in {\mathcal L}(W_0,W_1)$
whose graph in $W_0\oplus W_1=V$ is $W$.
We now give a description of the restrictions to $\Lambda$
of the local charts defined on $G_n(V)$ by this construction.

Given any real vector space $Z$, we denote by ${\rm B}(Z,\R)$ and
${\rm B}_{\rm sym}(Z,\R)$ respectively the space of bilinear forms
and symmetric bilinear forms on $Z$. There is an identification
of ${\rm B}(Z,\R)$ with ${\mathcal L}(Z,Z^*)$ obtained by associating
to each $B\in {\rm B}(Z,\R)$ the map $v\mapsto B(v,\cdot)$.
\begin{cartas}\label{thm:cartas}
Given a Lagrangian direct sum decomposition $V=L_0\oplus L_1$, for
all $W\in G_n(V)$ transverse to $L_1$, i.e., $W\cap L_1=\{0\}$, we
define $\phi_{L_0,L_1}(W)\in {\rm B}(L_0,\R)\simeq {\mathcal L}(L_0,L_0^*)$
by \[\phi_{L_0,L_1}(W)=\Iota{L_0}{L_1}\circ T,\]
where $T$ is the unique linear operator $T:L_0\mapsto L_1$ whose graph
in $V=L_0\oplus L_1$ is $W$.

The map $\phi_{L_0,L_1}$ is a diffeomorphism from the open set of $G_n(V)$
consisting of subspaces transverse to $L_1$ onto ${\rm B}(L_0,\R)$.
\end{cartas}
Observe that $\phi_{L_0,L_1}$ is simply one of the local charts on $G_n(V)$ described
above, up to the composition with the linear isomorphism $\Iota{L_0}{L_1}$.
We now show how the maps $\phi_{L_0,L_1}$ induce a submanifold structure
on $\Lambda$.
\begin{geometriaLambda}\label{thm:geometriaLambda}
The set $\Lambda$ is an analytic embedded submanifold of $G_n(V)$
of dimension $\frac12 n(n+1)$; 
each map $\phi_{L_0,L_1}$ restricts to a local chart on $\Lambda$ which
maps the open set of Lagrangian  subspaces transverse to $L_1$ onto
${\rm B}_{\rm sym}(L_0,\R)$. 

For all $L_0\in\Lambda$ the tangent space
$T_{L_0}\Lambda$ is  canonically isomorphic to ${\rm B}_{\rm sym}(L_0,\!\R)$;
more precisely, this isomorphism is given by the differential at $L_0$ of
any coordinate map $\phi_{L_0,L_1}$, and this isomorphism does not depend
on the choice of the complementary Lagrangian $L_1$. 

Moreover, the isomorphisms $T_{L_0}\Lambda\simeq {\rm B}_{\rm sym}(L_0,\R)$
are {\em natural\/} in the sense that, given a symplectomorphism $\psi$
of $(V,\omega)$, we have the following commutative diagram:
\begin{equation}\label{eq:CommDiag}
\begin{CD}
T_{L_0}\Lambda@>{\rm d}\hat\psi_{L_0}>>T_{\psi(L_0)}\Lambda\\
@VVV @VVV\\
{\rm B}_{\rm sym}(L_0,\R)@>\psi_*>> {\rm B}_{\rm sym}(\psi(L_0),\R),
\end{CD}
\end{equation}
where the vertical arrows are the canonical isomorphisms, $\hat\psi:\Lambda\mapsto\Lambda$
is the diffeomorphism given by $L\mapsto\psi(L)$, and $\psi_*$ is the push-forward
operator given by $B\mapsto B(\psi^{-1}.,\psi^{-1}.)$.
\end{geometriaLambda}
\begin{proof}
Let $L\in G_n(V)$ be transverse to $L_1$, and let $T:L_0\mapsto L_1$ be the linear
operator whose graph in $V=L_0\oplus L_1$ is $L$. Then, $L$ is Lagrangian
if and only if $\omega[v+T(v),w+T(w)]=0$ for all $v,w\in L_0$, i.e., if and
only if \[\omega[v,T(w)]+\omega[T(v),w]=0, \quad\forall\,v,w\in L_0.\] This is just the symmetry
of the bilinear form $\phi_{L_0,L_1}(L)=\Iota{L_0}{L_1}\circ T$.

We now prove that the differential ${\rm d}\phi_{L_0,L_1}(L_0)$ does
not depend on the choice of the complementary Lagrangian $L_1$; observe
that, by Remark~\ref{thm:remcomplLagr}, we can always find complementary Lagrangians
to $L_0$. To prove the claim, let $L_1$ and $L_2$ be two complementary Lagrangians
to $L_0$; the two charts $\phi_{L_0,L_1}$ and $\phi_{L_0,L_2}$ map $L_0$ to
the zero bilinear map. We have to prove that the differential of the 
transition map from $\phi_{L_0,L_1}$ to $\phi_{L_0,L_2}$ at $0$ is the
identity of ${\rm B}_{\rm sym}(L_0,\R)$. The transition map is given by:
\begin{equation}\label{eq:transition}
{\mathcal L}(L_0,L_0^*)\ni B\longmapsto B\circ({\rm Id}+\rho\circ\Iota{L_0}{L_1}^{-1}
\circ B)^{-1}\in{\mathcal L}(L_0,L_0^*),
\end{equation}
where $\rho$ is the restriction to $L_1$ of the projection $L_0\oplus L_2\mapsto L_0$
and ${\rm Id}$ is the identity on $L_0$. The differential of \eqref{eq:transition}
at $B=0$ is easily computed to be the identity.

It remains to prove the commutativity of \eqref{eq:CommDiag}. Let $\Omega$
and $\Omega'$ be the domains of the charts  $\phi_{L_0,L_1}$ and
$\phi_{\psi(L_0),\psi(L_1)}$ respectively. Then, it is easy to check
the commutativity of the diagram:
\begin{equation}\label{eq:CommDiag1}
\begin{CD}
\Lambda\supset\Omega\phantom{\Lambda\supset\,}@>\hat\psi>>\phantom{\,\subset\Lambda}
\Omega'\subset\Lambda\\
@V{\phi_{L_0,L_1}}VV @VV{\phi_{\psi(L_0),\psi(L_1)}}V\\
{\rm B}_{\rm sym}(L_0,\R)@>\psi_*>> {\rm B}_{\rm sym}(\psi(L_0),\R).
\end{CD}
\end{equation}
The conclusion follows by differentiating \eqref{eq:CommDiag1}.
\end{proof}
The action of ${\rm Sp}(V,\omega)$ on $\Lambda$ induces a map ${\rm sp}(V,\omega)\mapsto
T_{L_0}\Lambda$ for every $L_0$ in $\Lambda$. This map is described in the following:
\begin{LieAl}\label{thm:LieAl}
Let $L_0\in\Lambda$; define the map $\kappa_{L_0}:{\rm Sp}(V,\omega)\mapsto\Lambda$
by $\kappa_{L_0}(\psi)=\psi(L_0)$. The differential ${\rm d}\kappa_{L_0}({\rm Id})$
of $\kappa_{L_0}$ at the neutral element ${\rm Id}\in{\rm Sp}(V,\omega)$ maps each
$H\in{\rm sp}(V,\omega)$ to the symmetric bilinear form ${\rm d}\kappa_{L_0}({\rm Id})[H]
\in{\rm B}_{\rm sym}(L_0,\R)$ given by the restriction of $\omega(H\cdot,\cdot)$
to $L_0$.
\end{LieAl}
\begin{proof}
Let $L_1$ be any complementary Lagrangian to $L_0$, $V=L_0\oplus L_1$, 
and let $\phi_{L_0,L_1}$ be
the corresponding coordinate map around $L_0$. Recall that the differential ${\rm
d}\phi_{L_0,L_1}$ at $L_0$ is the isomorphism used to identify $T_{L_0}\Lambda$ with ${\rm
B}_{\rm sym}(L_0,\R)$ (see Proposition~\ref{thm:geometriaLambda}). Let $\pi_0:V\mapsto L_0$ and
$\pi_1:V\mapsto L_1$ be the projections onto the summands. 

In the chart $\phi_{L_0,L_1}$, the map $\kappa_{L_0}$ is given by:
\begin{equation}\label{eq:kappachart}
\psi\longmapsto\phi_{L_0,L_1}\circ\kappa_{L_0}(\psi)=
\Iota{L_0}{L_1}\circ\psi_{10}\circ\psi_{00}^{-1},
\end{equation}
where $\psi_{00}=\pi_0\circ(\psi\vert_{L_0})$ and 
$\psi_{10}=\pi_1\circ(\psi\vert_{L_0})$.
Formula \eqref{eq:kappachart} holds for $\psi$ in a neighborhood of
${\rm Id}\in{\rm Sp}(V,\omega)$, where $\psi_{00}$ is invertible. 

The differential of \eqref{eq:kappachart} is then easily computed as:
\[{\rm sp}(V,\omega)\ni H\longmapsto {\rm d}\kappa_{L_0}({\rm Id})[H]=
\Iota{L_0}{L_1}\circ H_{10},\]
where $H_{10}=\pi_1\circ (H\vert_{L_0})$.

The conclusion follows at once from the definition of $\Iota{L_0}{L_1}$.
\end{proof}
We recall that, if $f:M\mapsto N$ is a smooth map between differentiable
manifolds, two smooth vector fields $X$ and $Y$ on $M$ and $N$ respectively
are said to be {\em $f$-related\/} if ${\rm d}f(p)[X(p)]=Y(f(p))$ 
for all $p\in M$. If $X$ and $Y$ are $f$-related, then $f$ maps integral
curves of $X$ into integral curves of $Y$.

If a Lie group $G$ acts on the left on the manifold $M$, then to each $X$
in the Lie algebra of $G$ we associate a vector field $X^*$ in $M$
given by $X^*(p)={\rm d}\kappa_p(1)[X]$, where $\kappa_p:G\mapsto M$
is the map $g\mapsto g\cdot p$ and $1$ is the neutral element of $G$.
For all $p\in M$, the vector field $X^*$ is $\kappa_p$-related to the 
{\em right invariant\/} vector field on $G$ associated to $X$.

Considering $G={\rm Sp}(V,\omega)$ and $M=\Lambda$, we are not motivated
to give the following definition:
\begin{defH*}\label{thm:defH*}
Let $H\in{\rm sp}(V,\omega)$, the vector field $H^*$ in $\Lambda$ associate
to each $L\in\Lambda$ the vector $H^*(L)\in T_L\Lambda\simeq{\rm B}_{\rm sym}(L,\R)$
given by the restriction of $\omega(H\cdot,\cdot)$ to $L$.
\end{defH*}
The vector fields $H^*$ will be used to {\em project\/} differential
equations in ${\rm Sp}(V,\omega)$ to differential equations in $\Lambda$.
\smallskip

Using group actions, we now give a description of the geometrical 
structure of $\Lambda$ as a homogeneous space.
\begin{UnOn}\label{thm:UnOn}
$\Lambda$ is diffeomorphic to ${\rm U}(n)/{\rm O}(n)$; in particular,
$\Lambda$ is compact and connected.
\end{UnOn}
\begin{proof}
By choosing a symplectic basis for $(V,\omega)$, we reduce the problem
to the case $V=\R^{2n}$ and $\omega=\omega_0$. The
group ${\rm Sp}(n,\R)$ acts smoothly on $\Lambda$; 
we show that the restriction of this action
to ${\rm U}(n)$ is transitive on $\Lambda$. Let $L_0,L_1\in\Lambda$
be fixed; we consider bases ${\mathcal B}_0$ and ${\mathcal B}_1$
of $L_0$ and $L_1$ respectively, which are orthonormal relatively to the
Euclidean inner product of $\R^{2n}$. Since the imaginary part of the
Hermitian product is $\omega_0$, and $\omega_0$ vanishes on both $L_0$ and
$L_1$, then ${\mathcal B}_0$ and ${\mathcal B}_1$ are orthonormal
basis of $\C^n\simeq \R^{2n}$ with respect to the Hermitian product.
Hence, there exists an element of ${\rm U}(n)$ that carries
${\mathcal B}_0$ to ${\mathcal B}_1$, and ${\rm U}(n)$ acts
transitively on $\Lambda$.

Obviously, the isotropy group of $L_0=\R^n\oplus\{0\}$ is ${\rm O}(n)$,
which concludes the proof.
\end{proof}
We now give the following definition:
\begin{defLambdak}\label{thm:defLambdak}
Let $L_0\in\Lambda$ and $k=0,1,\ldots,n$ be fixed. We denote by $\Lambda_k(L_0)$ the
subset of $\Lambda$ consisting of Lagrangians $L$ with ${\rm dim}(L\cap L_0)=k$.
We also define the sets $\Lambda_{\le k}(L_0)$ and $\Lambda_{\ge k}(L_0)$
by:
\[\Lambda_{\le k}(L_0)=\bigcup_{i=0}^k\Lambda_i(L_0),\quad
\Lambda_{\ge k}(L_0)=\bigcup_{i=k}^n\Lambda_i(L_0).\]

\end{defLambdak}
\begin{abertofechado}\label{thm:abertofechado}
Clearly, $\Lambda_0(L_0)$ is precisely the set of all Lagrangians complementary
to $L_0$. It is an open set of $\Lambda$, since it is the domain of
any coordinate map $\phi_{L_1,L_0}$; moreover, it is diffeomorphic to a 
vector space by Proposition~\ref{thm:geometriaLambda}. 
For $k=0,\ldots,n$, we observe that $\Lambda_{\le k}(L_0)$ 
is open,  and so $\Lambda_{\ge k}(L_0)$ is closed in $\Lambda$. Namely, 
let $L\in \Lambda_{\le k}(L_0)$; we prove that $L$ admits a neighborhood
in $G_n(V)$ consisting only of subspaces $W$ such that ${\rm dim}(W\cap L_0)\le k$.
For, simply consider a subspace $W_1$ of $V$ which is complementary to both
$L_0$ and $L$; then, given a linear operator  $T:L_0\mapsto W_1$, its graph
in $L_0\oplus W_1=V$ intercepts $L_0$ in a subspace of dimension less than or equal
to $k$ if and only if ${\rm dim}({\rm Ker}(T))\le k$.
The conclusion follows easily by observing that the set of linear operators
$T\in{\mathcal L}(L_0,W_1)$ such that ${\rm dim}({\rm Ker}(T))\le k$ is open.
\end{abertofechado}
Given $L_0\in\Lambda$, we denote by ${\rm Sp}(V,\omega,L_0)$ the closed subgroup
of ${\rm Sp}(V,\omega)$ consisting of elements $\psi$ such that $\psi(L_0)=L_0$;
by ${\rm Sp}_+(V,\omega,L_0)$ we mean the subgroup of ${\rm Sp}(V,\omega,L_0)$
consisting of those $\psi$ whose restriction to $L_0$ is orientation preserving.
The Lie algebra ${\rm sp}(V,\omega,L_0)$ of both ${\rm Sp}(V,\omega,L_0)$
and ${\rm Sp}_+(V,\omega,L_0)$ is the subalgebra of ${\rm sp}(V,\omega)$
consisting of those $H$ such that $H(L_0)\subset L_0$.

Clearly, ${\rm Sp}(V,\omega,L_0)$ and ${\rm Sp}_+(V,\omega,L_0)$ act on 
all the spaces
$\Lambda_*(L_0)$ introduced in Definition~\ref{thm:defLambdak}. These actions 
are transitive on each
$\Lambda_k(L_0)$, as we prove in the following:
\begin{transitividade}\label{thm:transitividade}
For all $k=0,\ldots,n$, the group ${\rm Sp}_+(V,\omega,L_0)$ acts transitively
on $\Lambda_k(L_0)$.
\end{transitividade}
\begin{proof}
By choosing a symplectic basis of $(V,\omega)$, we can reduce to the case
$V=\R^{2n}$,  $\omega=\omega_0$ and $L_0=\R^n\oplus\{0\}$ (see
Remark~\ref{thm:remcomplLagr}); let $\{e_1,\ldots,e_{2n}\}$ be the canonical
basis of $\R^{2n}$. 
Let $L$ be any Lagrangian such that
${\rm dim}(L\cap L_0)=k$; we show that there is an element $\psi\in{\rm Sp}_+(V,\omega,L_0)$
such that $\psi(L)\cap L_0=\R^k\oplus\{0\}$. Let $\psi\in{\rm SO}(n)$ be a linear isometry
of $\R^n$ such that $\psi(L\cap L_0)=\R^k\oplus\{0\}$; now consider the complex linear extension
of $\psi$ to $\C^n\simeq\R^{2n}$. Such a map has the required property.

Let $L_1$ be the subspace generated by $\{e_1,\ldots,e_k,e_{n+k+1},\ldots,e_{2n}\}$.
Then, $L_1$ is Lagrangian, and $L_1\cap L_0=\R^k\oplus\{0\}$.
It remains to prove that, given a Lagrangian $L$ with $L\cap L_0=\R^k\oplus\{0\}$,
there exists an element $\psi\in{\rm Sp}_+(V,\omega,L_0)$ such that 
$\psi(L)=L_1$. 

To prove this claim, we define the following spaces. 
Let $V_1$ be  the space generated by $\{e_1,\ldots,e_k,e_{n+1},\ldots,e_{n+k}\}$;
$V_2$ be generated by $\{e_{k+1},\ldots,e_n,e_{n+k+1},\ldots,e_{2n}\}$ 
and $S$ be generated by $\{e_1,\ldots,e_n,e_{n+k+1},\ldots,e_{2n}\}$.
Observe that $S$ is the orthogonal complement
of $\R^k\oplus\{0\}$ with respect to $\omega_0$; also, $\R^{2n}= V_1\oplus V_2$,
and $\omega_0$ restricts to the canonical symplectic forms of $V_1\simeq\R^{2k}$
and of $V_2\simeq\R^{2(n-k)}$, that will be still denoted by $\omega_0$. Let
$\pi:S\mapsto V_2$ be the restriction to $S$ of the projection $V_1\oplus V_2\mapsto V_2$.
It is easy to check that $\omega_0(\pi(x),\pi(y))=\omega_0(x,y)$ for all $x,y\in S$.
Since $L$ is Lagrangian, we have $L\subset S$; moreover, it is easily seen that
$\pi(L)$ is Lagrangian in $V_2$. Since $L\cap L_0=\R^k\oplus\{0\}$,
we have that $\pi(L)$ is complementary to $\R^{n-k}\oplus\{0\}$ in $V_2
\simeq\R^{n-k}\oplus\R^{n-k}$.
By Remark~\ref{thm:remcomplLagr}, there exists a symplectomorphism $\varphi$
of $(V_2,\omega_0)$ that is the identity on $\R^{n-k}\oplus\{0\}$ and
carries $\pi(L)$ into $\{0\}\oplus\R^{n-k}=\pi(L_1)$. 

Finally, the required element $\psi\in{\rm Sp}_+(V,\omega,L_0)$ is given by:
\[\psi\lower5truept\hbox{$\big\vert$}_{V_1}={\rm Id},\quad
\psi\lower5truept\hbox{$\big\vert$}_{V_2}=\varphi.\]
Indeed, $\psi(L)=L_1$, because $\psi(L)$ and $L_1$ are both subspaces
of $S$ containing ${\rm Ker}(\pi)$ that have the same image under $\pi$.
This concludes the proof.
%
%
%
\end{proof}
\begin{complcomun}\label{thm:complcomun}
Given any two Lagrangians $L_0$ and $L$ in $\Lambda$, there exists $L_1\in\Lambda$
which is complementary to both $L_0$ and $L$. In particular, the domain of the
coordinate map $\phi_{L_0,L_1}$ contains both $L_0$ and $L$.
\end{complcomun}
\begin{proof}
By choosing a symplectic basis of $(V,\omega)$, we can reduce to the case
$V=\R^{2n}$,  $\omega=\omega_0$ and $L_0=\R^n\oplus\{0\}$ (see
Remark~\ref{thm:remcomplLagr}). 

Let $\{e_1,\ldots,e_{2n}\}$ be the canonical
basis of $\R^{2n}$ and  $L_2$ be the subspace generated 
by $\{e_1,\ldots,e_k,e_{n+k+1},\ldots,e_{2n}\}$, where $k={\rm dim}(L_0\cap L)$.
Since $L_2$ and $L$ are both in $\Lambda_k(L_0)$, Proposition~\ref{thm:transitividade}
gives a symplectomorphism $\psi$ of $(\R^{2n},\omega_0)$ such that $\psi(L_0)=L_0$
and $\psi(L_2)=L$. Observe that the diagonal \[\Delta=\{(v,v):v\in\R^n\}\] is
a Lagrangian subspace of $\R^{2n}$ which is complementary to both $L_0$ and $L_2$;
the desired Lagrangian $L_1$ is, for instance, $\psi(\Delta)$.
\end{proof}
Although we will not need it, we observe that the existence of complementary
Lagrangians can be proven in a much more general situation. Namely, using
Baire's Theorem, one proves that, given a sequence $\{L_k\}_{k\in\N}$ of
Lagrangians in $\Lambda$, the set $\bigcap\limits_{k\in\N}\Lambda_0(L_k)$
of their common complementary Lagrangians is dense in $\Lambda$. Each $\Lambda_0(L_k)$
is open dense because its complement in $\Lambda$ is a finite union of
embedded submanifolds of lower dimension, as we will see next.
\begin{geometriaLambdak}\label{thm:geometriaLambdak}
For all $k=0,\ldots,n$ and all $L_0\in\Lambda$, $\Lambda_k(L_0)$
is a connected embedded analytic submanifold of $\Lambda$ having codimension 
equal to $\frac12k(k+1)$.
For $L\in \Lambda_k(L_0)$, the tangent space $T_L\Lambda_k(L_0)\subset
{\rm B}_{\rm sym}(L,\R)$ is equal to the space of symmetric bilinear forms
on $L$ that vanish in $L\cap L_0$.

The submanifold $\Lambda_1(L_0)$, that has codimension $1$ in $\Lambda$
has a {\em transverse orientation\/} in $\Lambda$, namely, for $L\in \Lambda_1(L_0)$,
a vector $B\in {\rm B}_{\rm sym}(L)\simeq T_L\Lambda$ is {\em positive\/}
if $B$ is positive definite on the one-dimensional space $L\cap L_0$.
Moreover, the transverse orientation of $\Lambda_1(L_0)$ in $\Lambda$
is {\em natural\/} in the sense that, given $\psi\in{\rm Sp}(V,\omega,L_0)$,
the diffeomorphism $L\mapsto\psi(L)$ of $\Lambda$ is orientation preserving.
\end{geometriaLambdak}
\begin{proof}
To prove that $\Lambda_k(L_0)$ is an embedded submanifold of $\Lambda$, 
observe first that, by Proposition~\ref{thm:transitividade}, 
$\Lambda_k(L_0)$ is an orbit of the action of 
${\rm Sp}_+(V,\omega,L_0)$. It follows that $\Lambda_k(L_0)$ is an {\em immersed\/}
submanifold, i.e., it does not necessarily have the relative topology. 
By \cite[Theorem~2.9.7]{Vj}, an orbit is embedded if and only if it is {\em locally closed},
i.e., it is the intersection of an open and a closed set. Now, recall
Remark~\ref{thm:abertofechado} and simply
observe that $\Lambda_k(L_0)=\Lambda_{\ge k}(L_0)\cap \Lambda_{\le k}(L_0)$.

We now compute the codimension of $\Lambda_k(L_0)$ in $\Lambda$.

Let $L_1$ be any Lagrangian complementary to $L_0$;
the Lie group ${\rm Sp}(V,\omega,L_0)$ is diffeomorphic to ${\rm GL}(L_0)\times
{\rm B}_{\rm sym}(L_1,\R)$. Namely, we have a diffeomorphism:
\[F:{\rm GL}(L_0)\times{\rm B}_{\rm sym}(L_1,\R)\longmapsto{\rm Sp}(V,\omega,L_0)\]
that associates to each pair $(\alpha,\beta)$ the symplectomorphism $\psi=F(\alpha,\beta)$
of $(V,\omega)$ whose restriction to $L_0$ is $\alpha$ and whose restriction to
$L_1$ is equal to:
\[\psi\big\vert_{L_1}=\alpha\circ\Iota{L_1}{L_0}^{-1}\circ\beta+
\Iota{L_0}{L_1}^{-1}\circ(\alpha^*)^{-1}\circ\Iota{L_0}{L_1},\]
where $\alpha^*\in{\rm GL}(L_0^*,\R)$ denotes the transpose map of
$\alpha$, and $\beta$ is seen as a linear map $\beta:L_1\mapsto L_1^*$.

It follows that the dimension of ${\rm Sp}(V,\omega,L_0)$ is equal to $n^2+\frac12 n(n+1)$.
The group ${\rm Sp}_+(V,\omega,L_0)$ is the image under $F$ of the product
${\rm GL}_+(L_0)\times{\rm B}_{\rm sym}(L_1,\R)$, where ${\rm GL}_+(L_0)$ is the group of
orientation preserving isomorphisms of $L_0$. It follows that ${\rm Sp}_+(V,\omega,L_0)$ and
hence $\Lambda_k(L_0)$ is connected.

Now, we choose an element $L\in\Lambda_k(L_0)$ and we calculate the dimension
of its isotropy group in ${\rm Sp}(V,\omega,L_0)$. To this aim, let $S\subset L_0$ 
be any $k$-dimensional subspace and let $S'\subset L_1$ be the image under
$\Iota{L_0}{L_1}^{-1}$ of the annihilator of $S$ in $L_0^*$. Then, $L=S\oplus S'$
is a Lagrangian in $V$ and $L\in\Lambda_k(L_0)$.

The isotropy group of $L$ is the image under $F$ of the set of pairs $(\alpha,\beta)$
such that $\alpha(S)\subset S$ and $\beta$ vanishes on $S'$. It follows that the
dimension of this isotropy group is $n^2+\frac12 k(k+1)$. Hence, using
Proposition~\ref{thm:geometriaLambda}, the codimension
of $\Lambda_k(L_0)$ in $\Lambda$ is computed as $\frac12 k(k+1)$.

We now compute the tangent space $T_L\Lambda_k(L_0)$ at any point $L\in \Lambda_k(L_0)$.
Such a space is given by the image of ${\rm sp}(V,\omega,L_0)$ under
the differential ${\rm d}\kappa_L({\rm Id})$, defined in Proposition~\ref{thm:LieAl}: 
\[T_L\Lambda_k(L_0)=\Big\{\omega(H\cdot,\cdot)\big\vert_L:H\in
{\rm sp}(V,\omega,L_0)\Big\}.\]
The elements of $T_L\Lambda_k(L_0)$ vanish on $L_0\cap L$. A simple dimension
counting shows that $T_L\Lambda_k(L_0)$ consists precisely of those elements.
This completes the proof of the first part of the statement.

We now consider the submanifold $\Lambda_1(L_0)$; by the formula computed above, its
codimension in $\Lambda$ is equal to $1$. The transverse orientation is
well defined in the statement of the Proposition, and the naturality follows
easily from  the commutative diagram~\eqref{eq:CommDiag} in
Proposition~\ref{thm:geometriaLambda}.
\end{proof}
\begin{remtangent}\label{thm:remtangent}
In Proposition~\ref{thm:geometriaLambdak} we have given a description
of a tangent space $T_L\Lambda_k(L_0)$ as a subspace of ${\rm B}_{\rm sym}(L,\R)$,
where ${\rm B}_{\rm sym}(L,\R)$ is identified with $T_L\Lambda$ by means of a
coordinate map $\phi_{L,L_1}$ (see Proposition~\ref{thm:geometriaLambda}).

In many situations we will have to deal with curves $L (t)$ of Lagrangians,
and to study the tangent space $T_{L(t)}\Lambda$ it will be more convenient
to work with a {\em fixed\/} coordinate map $\phi_{L_0,L_1}$ rather than 
using variable charts $\phi_{L(t),L_1}$. 

For this reason, we now describe the transition map from a coordinate
map $\phi_{L_0,L_1}$ to $\phi_{L,L_1}$, where $L_1$ is a complementary Lagrangian
to both $L_0$ and $L$ (see Corollary~\ref{thm:complcomun}). 

Let $\eta:L_0\mapsto L$ be the isomorphism obtained by the restriction to
$L_0$ of the projection $L\oplus L_1\mapsto L$. The transition map from
$\phi_{L_0,L_1}$ to $\phi_{L,L_1}$ is now easily computed as:
\[{\rm B}_{\rm sym}(L_0,\R)\ni B\longmapsto\phi_{L,L_1}(L_0)+\eta_*( B),\]
where $\eta_*$ is the push-forward operator given by $\eta_*(B)=B(\eta^{-1}.,\eta^{-1}
.)$. Thus, the transition map is $\eta_*$ plus a translation by a fixed element,
and so its differential at any point is given by $\eta_*$.

Observe that $\eta$ is the identity in $L_0\cap L$, therefore we get
\begin{equation}\label{eq:identification}
{\rm d}\phi_{L_0,L_1}(L)\big[T_L\Lambda_k(L_0)\big]=\Big\{B\in{\rm B}_{\rm sym}(L_0,\R):
B\ \text{vanishes on}\ L\cap L_0\Big\}.
\end{equation}
The reader should compare formula \eqref{eq:identification} with
the description of $T_L\Lambda_k(L_0)$ given in the statement of
Proposition~\ref{thm:geometriaLambdak}. Observe also that, for $L\in\Lambda_1(L_0)$,
since the push-forward operator does not affect the positivity of a bilinear form, 
a given $B\in{\rm B}_{\rm sym}(L_0,\R)$ is such that ${\rm d}\phi_{L_0,L_1}(L)^{-1}[B]$
is a positive vector in the transverse orientation of $\Lambda_1(L_0)$ if and
only if $B$ is positive definite on $L\cap L_0$.
\end{remtangent}
\end{section}
\begin{section}{Intersection Theory: the Maslov Index}\label{sec:intersection}
The purpose of this Section is to associate an integer number to each 
pair $(\ell,L_0)$, where $\ell$ is a continuous curve in the Lagrangian 
Grassmannian $\Lambda$ studied in the previous
section, and $L_0\in\Lambda$. Such a number, that will be defined 
to be the {\em Maslov Index\/} of $\ell$ with respect to $L_0$, in 
the generic case will count (algebraically) the number of intersections
of $\ell$ with $\Lambda_{\ge1}(L_0)$.
\smallskip

We will assume throughout the Section that $(V,\omega)$ is a fixed
symplectic space of dimension $2n$, and that $\Lambda$ is the associated
Lagrangian Grassmannian.
\mysubsection{The fundamental group of the Lagrangian Grassmannian}
\label{sub:fundamental}

We begin with an easy result  on the fundamental group of homogeneous spaces:
\begin{quotient}\label{thm:quotient}
Let $G$ be a connected Lie group and $K$ be a closed subgroup
of $G$; we denote by $p:G\mapsto G/K$ the quotient map. 
Let $q:\tilde G\mapsto G$ be the universal covering group of $G$,
$\tilde K=q^{-1}(K)$ and $\tilde K_0$ be the connected component of the
neutral element $1\in\tilde K$. Then, the fundamental group
$\pi_1(G/K)$ is isomorphic to the quotient $\tilde K/\tilde K_0$.
The isomorphism 
\[\zeta:\tilde K/\tilde K_0\mapsto \pi_1(G/K,p(1))\]
is defined as follows. If $g\tilde K_0$ is any element of $\tilde K/\tilde K_0$,
let $c:[0,1]\mapsto \tilde G$ be any continuous curve such that
$c(0)=1\in\tilde G$ and $c(1)=g^{-1}$. Then, $\zeta(g\tilde K_0)$
is the homotopy class of the loop $p\circ q\circ c:[0,1]\mapsto G/K$
based in $p(1)$.
\end{quotient}
\begin{proof}
We start by determining the universal covering of the quotient $G/K$.
Towards this goal, we consider the transitive action of $\tilde G$ on $G/K$
given by $g\cdot(xK)=(q(g)x)K$, for all $g\in\tilde G$ and $x\in G$.
The isotropy group of $p(1)=1K$ is $\tilde K=q^{-1}(K)$; we have therefore
a diffeomorphism $\tilde G/\tilde K\mapsto G/K$ given by $g\tilde K\mapsto q(g)K$,
for all $g\in\tilde G$.

Since $\tilde K/\tilde K_0$ is discrete, then the map $\tilde G/\tilde K_0\mapsto
\tilde G/\tilde K$ given by $g\tilde K_0\mapsto g\tilde K$ is a covering map.

Considering the composition of the two maps above, we obtain
a covering map $\overline q:\tilde G/\tilde K_0\mapsto G/K$ given by
$\overline q(g\tilde K_0)=q(g)K$, $g\in\tilde G$.

Since $\tilde G$ is simply connected and $\tilde K_0$ is connected, the quotient
$\tilde G/\tilde K_0$ is simply connected, and so $\overline q$ is the
universal covering map of $G/K$.

We now determine the group of covering automorphisms of $\overline q$,
which is isomorphic to $\pi_1(G/K)$. 
We recall that an automorphism $\varphi$ of $\overline q$ is a homeomorphism
of $\tilde G/\tilde K_0$ such that $\overline q\circ\varphi=\overline q$.
For all $g\in\tilde K$, the map $x\tilde K_0\mapsto xg^{-1}\tilde K_0$ is 
an automorphism of $\overline q$ which is trivial if $g\in\tilde K_0$. 
Thus, we  have an action of $\tilde K/\tilde K_0$ in $\tilde G/\tilde K_0$
by automorphisms of $\overline q$, which is transitive and simple 
(i.e., without fixed points) on the fibers of $\overline q$.

It follows that $\tilde K/\tilde K_0$ is the group of covering automorphisms
of $\overline q$, concluding the proof of the first part of the statement.

To construct an explicit isomorphism between $\tilde K/\tilde K_0$ and
$\pi_1(G/K,p(1))$, one uses the standard procedure of taking homotopy classes
of loops obtained as the images under $\overline q$ of curves in $\tilde G/\tilde K_0$
that connect the point $1\tilde K_0$ and a generic point in the fiber
$\overline q^{-1}(1K)$.
\end{proof}
We use Lemma~\ref{thm:quotient} to compute the fundamental group of the
Lagrangian Grassmannian $\Lambda$: 
\begin{pi1}\label{thm:pi1}
Let $\{e_1,\ldots,e_{2n}\}$ be a symplectic
basis of $(V,\omega)$
and let $L_0\in\Lambda$ be the Lagrangian subspace generated by $\{e_1,\ldots,e_n\}$. 
Then, the fundamental group of $\Lambda$ with base point
$L_0$, $\pi_1(\Lambda, L_0)$, is isomorphic to $\Z$. 
A generator of
$\pi_1(\Lambda,L_0)$ is given by the homotopy class of the loop $\ell:[0,1]\mapsto\Lambda$, 
where $\ell (t)$ is the Lagrangian generated by the vectors 
$e_1\cdot\cos\pi t-e_{n+1}\cdot\sin\pi t, e_2,\ldots,e_n$.
\end{pi1}
\begin{proof}
We can clearly assume that $V=\R^{2n}$, $\omega=\omega_0$, with
$\{e_1,\ldots,e_{2n}\}$ the canonical basis of $\R^{2n}$, hence,
$L_0=\R^n\oplus\{0\}$.
We apply Lemma~\ref{thm:quotient} to $G={\rm U}(n)$, $K={\rm O}(n)$; 
by Proposition~\ref{thm:UnOn}, we can identify $\Lambda$ with $G/K$, and
the quotient map $p:{\rm U}(n)\mapsto\Lambda$ is given by ${\rm U}(n)\ni\psi\mapsto
\psi(L_0)\in\Lambda$.

Let ${\rm SU}(n)$ be the Lie group of unitary $n\times n$ complex matrices
having determinant equal to $1$; the universal covering
group of $G$ is $\tilde G={\rm SU}(n)\times\R$, with covering map
$q(A,t)=e^{it}\cdot A$.

The group $\tilde K=q^{-1}(K)$ is easily computed as:
\begin{equation}\label{eq:union}
\tilde K=\bigcup_{k\in \Z}\Big[e^{-i\frac{k\pi}n}\,{\rm
O}_k(n)\Big]\times
\Big\{\frac{k\pi}n\Big\}, 
\end{equation}
where by ${\rm O}_k(n)$ we mean ${\rm SO}(n)$ if $k$ is even and its complement
${\rm O}(n)\setminus{\rm SO}(n) $ if $k$ is odd. \smallskip

The connected component $\tilde K_0$ equals ${\rm SO}(n)\times\{0\}$, and
$\tilde K/\tilde K_0$ is isomorphic to $\Z$. Such an isomorphism is given
by mapping each term of the union in formula \eqref{eq:union}
into the integer $k$.

As a generator for $\tilde K/ \tilde K_0$, we choose the term in \eqref{eq:union}
corresponding to $k=1$; such element is of the form $g\tilde K_0$, where $g$
is chosen to be the pair $g=(e^{-i\,A},\frac\pi n)$, with $A$ the diagonal matrix:
\[A=\left(\begin{array}{cccc}\pi(\frac{1-n}n)& & & \\
& \frac\pi n & & \\
& & \ddots & \\
& & & \frac\pi n\end{array}\right).\]
Observe that $A$ is a traceless Hermitian matrix, so that $i\,A$ belongs to the Lie
algebra ${\rm su}(n)$ of ${\rm SU}(n)$.

\noindent
In order to determine a generator for $\pi_1(\Lambda,L_0)$, we choose the curve
$c:[0,1]\mapsto\tilde G$ given by $c(t)=(e^{i\,t\,A},-t\,\frac\pi n)$, connecting
the neutral element of $\tilde G$ with $g^{-1}$. 

The curve $\ell= p\circ q\circ c$ in $\Lambda$ is now easily computed as given
in the statement of the Corollary.
\end{proof}
\mysubsection{The intersection theory and the construction of the Maslov index}
\label{sub:construction}
In order to develop our intersection theory, we are interested in
the singular homology groups of $\Lambda$. The first homology group
of a path connected topological space is isomorphic to the abelianization
of its fundamental group, and therefore, it follows from Corollary~\ref{thm:pi1}
that $H_1(\Lambda)$ is isomorphic to $\Z$. 

Each loop $\ell:[a,b]\mapsto\Lambda$ defines a homology class in
$H_1(\Lambda)$. Given any $L_0\in\Lambda$,
we obtain a homomorphism $\pi_1(\Lambda,L_0)\mapsto H_1(\Lambda)$, 
which associates to the homotopy class of a loop based in $L_0$ its 
homology class. This is called the {\em Hurewicz homomorphism\/}
and it is well known that it is surjective, and its kernel is the
commutator subgroup of $\pi_1(\Lambda,L_0)$ (see~\cite[Proposition~4.21]{Vi}).

The homology class of the curve $\ell:[0,1]\mapsto\Lambda$ defined in Corollary~\ref{thm:pi1}
is therefore a generator of $H_1(\Lambda)\simeq\Z$, and, up to the choice
of a sign, to each loop in $\Lambda$ we have a well defined integer associated
to it. Such a number is to be interpreted as a sort of {\em winding number\/}
of the loop around $\Lambda$.

Using the relative homology groups, we now show how to extend the above
construction to curves in $\Lambda$ that are not necessarily closed.
Let $L_0$ be a fixed Lagrangian in $\Lambda$; we will consider the relative
homology group $H_1(\Lambda,\Lambda_0(L_0))$. We recall that $\Lambda_0(L_0)$
is the complement in $\Lambda$ of the set $\Lambda_{\ge1}(L_0)$; every continuous
curve $\ell:[a,b]\mapsto\Lambda$ with endpoints outside $\Lambda_{\ge1}(L_0)$ defines
a homology class in $H_1(\Lambda,\Lambda_0(L_0))$.

We observe that $\Lambda_0(L_0)$ is {\em contractible}, since any coordinate
map $\phi_{L_1,L_0}$ maps $\Lambda_0(L_0)$ diffeomorphically onto a vector space
(see Remark~\ref{thm:abertofechado}). Hence, by the long exact homology sequence
of the pair $(\Lambda,\Lambda_0(L_0))$, we have an isomorphism 
$H_1(\Lambda)\mapsto H_1(\Lambda,
\Lambda_0)$ induced by the inclusion 
$(\Lambda,\emptyset)\mapsto(\Lambda,\Lambda_0(L_0))$. This implies
that $H_1(\Lambda,\Lambda_0(L_0))\simeq\Z$, and a generator of this group
is the homology class of the curve $\ell$ given in the statement of 
Corollary~\ref{thm:pi1}.

We make some remarks that follow from elementary properties
of the homology theory. Let $\ell:[a,b]\mapsto\Lambda$
be a continuous curve with endpoints outside $\Lambda_1(L_0)$, and
 $\sigma:[c,d]\mapsto[a,b]$  
a continuous  map with $\sigma(c)=a$ and $\sigma(d)=b$. 
Then, the curves $\ell $ and $\ell\circ\sigma$ are homologous  
in $H_1(\Lambda,\Lambda_0(L_0))$. If $\sigma(d)=a$ and $\sigma(c)=b$,
then $\ell\circ\sigma$ is homologous to the singular $1$-chain $-\ell$
in $H_1(\Lambda,\Lambda_0(L_0))$. If $u\in\,]a,b[$ is such that
$\ell(u)\not\in\Lambda_{\ge1}(L_0)$, then $\ell$ is homologous
to the singular $1$-chain $\ell\vert_{[a,u]}+\ell\vert_{[u,b]}$.
Clearly, if the image of $\ell$ does not intersect $\Lambda_{\ge1}(L_0)$,
then $\ell$ is homologous to zero in $H_1(\Lambda,\Lambda_0(L_0))$.
Finally, if $\ell_1,\ell_2:[a,b]\mapsto\Lambda$ are continuous curves
with endpoints outside $\Lambda_{\ge1}(L_0)$ and that are homotopic 
through curves with endpoints outside $\Lambda_{\ge1}(L_0)$, then,
they are homologous in $H_1(\Lambda,\Lambda_0(L_0))$.
\begin{induced}\label{thm:induced}
Each element in ${\rm Sp}(V,\omega)$ induces the identity homomorphism
in the homology of $\Lambda$, and each element of ${\rm Sp}(V,\omega,L_0)$
induces the identity in the relative homology of the pair $(\Lambda,\Lambda_0(L_0))$,
in the following sense. Recall that a continuous map $f$ between (pairs of) topological
spaces induces homomorphisms between their (relative) homology groups, that
will denoted by $(f)_*$.

If $\psi\in{\rm Sp}(V,\omega)$ and $\hat\psi:\Lambda\mapsto\Lambda$ is the diffeomorphism
given by $L\mapsto\psi(L)$, then $(\hat\psi)_*:H_1(\Lambda)\mapsto H_1(\Lambda)$
is the identity map. For, it is well known that ${\rm Sp}(V,\omega)\simeq{\rm Sp}(n,\R)$ 
is connected, hence every $\psi$ can be continuously connected to the neutral
element of ${\rm Sp}(V,\omega)$, which gives a homotopy between
$\hat\psi$ and the identity of $\Lambda$.

Moreover, if $\psi\in{\rm Sp}(V,\omega,L_0)$, then $\hat\psi$ carries $\Lambda_0(L_0)$
onto itself, and $(\hat\psi)_*$ is again the identity on $H_1(\Lambda,\Lambda_0(L_0))$.
To see this, it suffices to observe that the following diagram commutes, by standard
functoriality properties in homology:
\[
\begin{CD}
H_1(\Lambda,\Lambda_0(L_0))@>(\hat\psi)_*>>H_1(\Lambda,\Lambda_0(L_0))\\
@A{(i)_*}AA @AA{(i)_*}A\\
H_1(\Lambda)@>>{\rm Id}=(\hat\psi)_*>H_1(\Lambda),
\end{CD}
\]
where $i:(\Lambda,\emptyset)\mapsto(\Lambda,\Lambda_0)$ is the inclusion.
\end{induced}

We now give the following sufficient condition for two curves to be
homologous in $H_1(\Lambda,\Lambda_0(L_0))$:
\begin{preliminary}\label{thm:preliminary}
Let $\ell_1,\ell_2:[a,b]\mapsto\Lambda$ be continuous curves
with endpoints outside $\Lambda_{\ge1}(L_0)$. Suppose that
there exists a Lagrangian $L_1$ complementary to $L_0$
such that the images of both $\ell_1$ and $\ell_2$ are entirely contained
in the domain $\Lambda_0(L_1)$ of $\phi_{L_0,L_1}$. Let 
$\beta_i=\phi_{L_0,L_1}\circ\ell_i$, $i=1,2$; then, if 
$n_+(\beta_1(t))=n_+(\beta_2(t))$ for $t=a,b$, it follows that
$\ell_1$ and $\ell_2$ are homologous in $H_1(\Lambda,\Lambda_0(L_0))$.
\end{preliminary}
\begin{proof}
We define the space:
\begin{equation}\label{eq:defBksym}
\Bge1{L_0}=\Big\{B\in{\rm B}_{\rm sym}(L_0,\R):\dgn(B)\ge1\Big\},
\end{equation}
it is easy to see that $\phi_{L_0,L_1}\big(\Lambda_{\ge1}(L_0)\cap\Lambda_0(L_1)\big)=
\Bge1{L_0}$.  We also denote by ${\rm B}_{\rm sym}^0(L_0,\R)$
the complement of $\Bge1{L_0}$ in ${\rm B}_{\rm sym}(L_0,\R)$;
${\rm B}_{\rm sym}^0(L_0,\R)$ is given by the union of $n+1$ open
connected components ${\rm B}_{\rm sym}^{0,i}(L_0,\R)$, given by
\begin{equation}\label{eq:defBkisym}
{\rm B}_{\rm sym}^{0,i}(L_0,\R)=\Big\{B\in{\rm B}_{\rm sym}(L_0,\R):n_+(B)=i
\Big\},\quad i=0,\ldots,n.
\end{equation} 
Observe that each ${\rm B}_{\rm sym}^{0,i}(L_0,\R)$ is indeed path connected, because,
by Sylvester's Inertia Theorem, it admits a transitive action of the connected
group ${\rm GL}_+(n,\R)$.

If we set $i=n_+(\beta_1(a))=n_+(\beta_2(a))$ and $j=n_+(\beta_1(b))=n_+(\beta_2(b))$,
then we can find a continuous 
curve $\beta_3$  (and $\beta_4$) in ${\rm B}_{\rm sym}^{0,i}(L_0,\R)$
(in ${\rm B}_{\rm sym}^{0,j}(L_0,\R)$) from $\beta_2(a)$ to $\beta_1(a)$
(from $\beta_1(b)$ to $\beta_2(b)$). 

Define $\ell_i=\phi_{L_0,L_1}^{-1}\circ\beta_i$, $i=3,4$. Then, since $\ell_3$
and $\ell_4$ do not intersect $\Lambda_{\ge1}(L_0)$, the {\em concatenation\/}
$\ell=\ell_3\diamond\ell_1\diamond\ell_4$ is homologous to $\ell_1 $
in $H_1(\Lambda,\Lambda_0(L_0))$.

Let $\beta=\phi_{L_0,L_1}\circ\ell$; then, 
$\beta=\beta_3\diamond\beta_1\diamond\beta_4$. Since $\beta$ and $\beta_2$
have the same endpoints, they admit a fixed endpoint homotopy in the vector space
${\rm B}_{\rm sym}(L_0,\R)$. The composition of such homotopy with
$\phi_{L_0,L_1}^{-1}$ gives a homotopy between $\ell$ and $\ell_2$ through
curves with endpoints outside $\Lambda_{\ge1}(L_0)$.

Hence, $\ell$ and $\ell_2$ are homologous in $H_1(\Lambda,\Lambda_0(L_0))$, and
we are done.
\end{proof}
In the next Lemma, we show how to compute the variation of the
type numbers for a differentiable curve of symmetric 
bilinear forms:
\begin{elementary}\label{thm:elementary}
Let $Z$ be a finite dimensional real vector space and 
$A:[0,r]\mapsto {\rm B}_{\rm sym}(Z,\R)$
be a map of class $C^1$. Suppose that the restriction $\tilde A$ of the derivative
$A'(0)$ to the kernel ${\rm Ker}(A(0))$ is non degenerate. Then, for
$t>0$ sufficiently small, $A(t)$ is non degenerate, and we have:
\begin{equation}\label{eq:sumn+-}
n_+(A(t))=n_+(A(0))+n_+(\tilde A),\quad n_-(A(t))= n_-(A(0))+n_-(\tilde A).
\end{equation}
\end{elementary}
\begin{proof}
Let $N={\rm Ker}(A(0))$.
We start with the case where $A(0)$ is positive semidefinite, i.e.,
$n_-(A(0))=0$, and $\tilde A$ is positive definite, i.e.,
$n_-(\tilde A)=\dgn(\tilde A)=0$. Let $S$ be a subspace of $Z$ which is
complementary to $N$ and such that $A(0)$ is positive definite on $S$.
We need to show that $A(t)$ is positive definite
on $Z=N\oplus S$ for $t$ small enough. First, since $A(0)$ is positive
definite on $S$, there is $\varepsilon>0$ such that $A(t)$ is positive
definite on $S$ for $t\in[0,\varepsilon]$ (the set of positive definite 
symmetric bilinear forms is open). Let $\Vert\cdot\Vert$ be
an arbitrary norm on $Z$ and define:
\begin{equation}\label{eq:first}
c_0=\inf_{\stackrel{t\in[0,\varepsilon]}{x\in S,\ \Vert
x\Vert=1}}A(t)[x,x]>0,
\quad c_1=\inf_{\stackrel{y\in N}{\Vert y\Vert=1}}A'(0)[y,y]>0.
\end{equation}
It is easy to see that, for all $t>0$
small enough, we have 
\begin{equation}\label{eq:third}
\phantom{\quad\forall\, y\in N,\,\Vert y\Vert=1,}
A(t)[y,y]\ge\frac12 c_1t,\quad\forall\, y\in N,\,\Vert y\Vert=1,
\end{equation} so that $A(t)$ is positive 
definite on both $N$ and $S$ for $t>0$ small enough. We want to show that,
if $t>0$ is small enough, then for all $x\in S\setminus\{0\}$ and $y\in N\setminus\{0\}$,
$A(t)$ is positive definite on the two dimensional subspace of $Z$ generated by
$x$ and $y$. By the positivity on $S$ and $N$, it suffices to prove that,
for $t>0$ small enough, the 
following inequality holds:
\begin{equation}\label{eq:schwartz}
A(t)[x,y]^2\le A(t)[x,x]\cdot A(t)[y,y],
\end{equation}
for all $x\in S$, $y\in N$. Obviously, we can assume
$\Vert x\Vert=\Vert y\Vert=1$. As $A(0)$ vanishes on $N\times S$, 
there exists $c_2>0$ such that, for all $t>0$ small enough, we have:
\begin{equation}\label{eq:second}
\vert A(t)[x,y]\vert\le c_2\cdot t,
\end{equation}
for all $x\in S$, $y\in N$ with $\Vert x\Vert=\Vert y\Vert=1$.
By \eqref{eq:first}, \eqref{eq:third} and \eqref{eq:second},
for all $t>0$ small enough we get:
\[A(t)[x,y]^2\le c_2^2t^2\le\frac12 c_0\,c_1\,t\le A(t)[x,x]\cdot A(t)[y,y],\]
for all $x\in S$, $y\in N$ with $\Vert x\Vert=\Vert y\Vert=1$.
This yields \eqref{eq:schwartz} and concludes the first part of the proof.
\smallskip

For the general case, we consider decompositions $Z=S_+\oplus S_-\oplus N$
and $N=N_+\oplus N_-$, where $A(0)$ is positive definite on $S_+$,
negative definite on $S_-$, and $A'(0)$ is positive definite in $N_+$
and negative definite in $N_-$. We then apply the result proven 
in the first part of the proof to the restriction of $A(t)$ to $S_+\oplus N_+$ 
once, and again to the restriction of $-A(t)$ to $S_-\oplus
N_-$. The conclusion follows by observing that $A(t)$ is positive definite
on $S_+\oplus N_+$ and negative definite on $S_-\oplus N_-$, which
implies that $n_+(A(t))={\rm dim}(S_+\oplus N_+)$ and $n_-(A(t))=
{\rm dim}(S_-\oplus N_-)$ for $t>0$ small enough. Clearly, this also implies
that $A(t)$ is non degenerate.
\end{proof}

%
%

We now go back to the study of the homology of the pair $(\Lambda,\Lambda_0(L_0))$
and of the intersection theory.
 
The choice of an isomorphism $H_1(\Lambda,\Lambda_0(L_0))\simeq\Z$
is equivalent to the choice of one of the two generators of 
$H_1(\Lambda,\Lambda_0(L_0))$. Using the canonical transverse
orientation of $\Lambda_1(L_0)$ (see Proposition~\ref{thm:geometriaLambdak}), 
we now show how such a choice will be made.

In order to give a precise statement of our next Proposition, we give the following
definitions. Let $\ell:[a,b]\mapsto\Lambda$ be a smooth curve, with $\ell (t_0)\in\Lambda_{\ge1}
(L_0)$ for some $t_0\in[a,b]$. We say that $\ell$ {\em intercepts $\Lambda_{\ge1}(L_0)$
transversally\/} at the instant $t_0$ if $\ell(t_0)\in\Lambda_1(L_0)$ and
the tangent vector $\ell'(t_0)$ does not belong to $T_{\ell(t_0)}\Lambda_1(L_0)$.
If $\ell $ intercepts $\Lambda_{\ge1}(L_0)$ transversally at $t_0$, we say
that this intersection is {\em positive\/} or {\em negative\/}
if $\ell'(t_0)$ is respectively a positive or a negative vector in the sense of 
the transversal orientation of $\Lambda_1(L_0)$ (see Proposition~\ref{thm:geometriaLambdak}).

\begin{homologous}\label{thm:homologous}
Let $\ell_1$ and $\ell_2$ be smooth curves in $\Lambda$ with both endpoints
outside $\Lambda_{\ge1}(L_0)$. Suppose that both curves intercept
$\Lambda_{\ge1}(L_0)$ only once, and  that such intersections
are both transverse and positive (or both negative). Then $\ell_1$ and $\ell_2$ define the same
homology class in $H_1(\Lambda,\Lambda_0(L_0))$. Moreover, this homology class is a generator
of $H_1(\Lambda,\Lambda_0(L_0))\simeq\Z$.
\end{homologous}
\begin{proof}
We consider the case of positive intersections; the other case
is then easily  obtained by passing to the backwards orientation.
By reparameterizing, we can assume that both curves intercept $\Lambda_1(L_0)$
at the same instant $t_0\in\,]a,b[$. By Proposition~\ref{thm:transitividade},
there exists $\psi\in{\rm Sp}_+(V,\omega,L_0)$ such that $\psi(\ell_2(t_0))=
\ell_1(t_0)$. Let $\hat\psi$ denote the  diffeomorphism of $\Lambda$ given
by $L\mapsto\psi(L)$; we deduce from Proposition~\ref{thm:geometriaLambdak}
that the curve $\hat\psi\circ\ell_2$ has a unique intersection
with $\Lambda_{\ge1}(L_0)$, which is transverse and positive. Moreover,
by Remark~\ref{thm:induced}, $\hat\psi\circ\ell_2$ and $\ell_2$ are homologous
in $H_1(\Lambda,\Lambda_0(L_0))$.

Without loss of generality, we can therefore assume that $\ell_1(t_0)=\ell_2(t_0)$.
Let $L_1\in\Lambda$ be a Lagrangian which is complementary to both
$L_0$ and $\ell_1(t_0)=\ell_2(t_0)$ (see Corollary~\ref{thm:complcomun}).
Since $\ell_i$, $i=1,2$, has a unique intersection with $\Lambda_{\ge1}(L_0)$
at $t_0$, then the restriction of $\ell_i$ to any closed subinterval  containing $t_0$
in its interior is homologous to $\ell_i$ in $H_1(\Lambda,\Lambda_0(L_0))$.
Thus, we can assume that the images of $\ell_1$ and $\ell_2$ are contained in
the domain $\Lambda_0(L_1)$ of the chart $\phi_{L_0,L_1}$.

Let  $\beta_i=\phi_{L_0,L_1}\circ\ell_i$, $i=1,2$; one checks easily that
${\rm Ker}(\beta_i(t_0))=\ell_i(t_0)\cap L_0$ and, since the intersection
of $\ell_i$ with $\Lambda_{\ge1}(L_0)$ is unique, $\beta_i(t)$ is
non degenerate for all $t\ne t_0$. It follows that $n_+(\beta_i(t))$ is
constant for $t\in[a,t_0\,[$ and for $t\in \,]t_0,b]$. 

By the positivity of the intersection, the restriction of $\ell_i'(t_0)$
to the one dimensional subspace $\ell_i(t_0)\cap L_0$ is positive
definite (see Proposition~\ref{thm:geometriaLambdak}). Moreover, by
Remark~\ref{thm:remtangent}, it follows that $\beta'_i(t_0)$
is positive definite on $\ell_i(t_0)\cap L_0$.

The fact that $\ell_1$ is homologous to $\ell_2$ in $H_1(\Lambda,\Lambda_0(L_0))$
will follow from Lemma \ref{thm:preliminary} once we prove that
$n_+(\beta_1(a))=n_+(\beta_2(a))$ and that $n_+(\beta_1(b))=n_+(\beta_2(b))$.

Applying twice Lemma~\ref{thm:elementary} around $t_0$, we obtain the following
equalities for each $i=1,2$:
\[n_+(\beta_i(b))=n_+(\beta_i(t_0))+1,\quad n_+(\beta_i(a))=n_+(\beta_i(t_0));\]
for the second equality we have applied Lemma~\ref{thm:elementary} to 
the curve $\beta_i$ reparameterized {\em backwards}. The conclusion follows from
the fact that $\beta_1(t_0)=\beta_2(t_0)$.
\smallskip

Using the first part of the Proposition, to conclude the proof
we need to exhibit a smooth curve $\ell$ whose homology class
is a generator of $H_1(\Lambda,\Lambda_0(L_0))$, and that intersects
$\Lambda_{\ge1}(L_0)$ exactly once, with such intersection 
transverse.

To this aim, let $\{e_1,\ldots,e_{2n}\}$ be a symplectic basis of $(V,\omega)$
such that $L_0$ is the Lagrangian generated by $\{e_{n+1},\ldots,e_{2n}\}$
(see Remark~\ref{thm:remcomplLagr}). Consider the curve $\ell:[0,1]\mapsto\Lambda$ 
introduced in
the statement of Corollary~\ref{thm:pi1}. It intersects $\Lambda_{\ge1}(L_0)$ only 
at the instant
$t_0=\frac12$ and $\ell(\frac12)\in\Lambda_1(L_0)$, because $\ell(\frac12)\cap L_0=
\R\cdot e_{n+1}$.

To check the transversality, we make computations using the chart
$\phi_{L_0,L_1}$ (rather than $\phi_{\ell(\frac12),L_1}$, see
Remark~\ref{thm:remtangent}), where we choose $L_1$ to be the
Lagrangian generated by the vectors $e_i+e_{n+i}$, $i=1,\ldots,n$, which is complementary
to both $L_0$ and $\ell(\frac12)$. We set $\beta=\phi_{L_0,L_1}\circ\ell$, and
we obtain for each $t\in[0,1]$ a symmetric bilinear form in $L_0$ which, in the basis
$\{e_{n+1},\ldots,e_{2n}\}$ is given by the diagonal matrix:
\[\beta(t)\sim\left(
\begin{array}{cccc}
f(t)& & &  \\
    &1& &  \\
    & &\ddots& \\
    & & & 1\\
\end{array}\right),\]
with
\[f(t)=\frac{\cos\pi t}{\cos\pi t+\sin\pi t}.\]
Since $\ell'(\frac12)(e_{n+1},e_{n+1})=f'(\frac12)=-\pi\ne0$, then $\ell$ intersects
$\Lambda_{\ge1}(L_0)$ transversally (with negative intersection), and the proof
is complete.
\end{proof}
Let $\mu_{L_0}:H_1(\Lambda,\Lambda_0(L_0))\mapsto\Z$ be the unique isomorphism
such that $\mu_{L_0}(\mathfrak h)=1$ where $\mathfrak h$ is the homology class
of any smooth curve $\ell$ in $\Lambda$, with endpoints outside $\Lambda_{\ge1}(L_0)$
and intersecting only once $\Lambda_{\ge1}(L_0)$, such intersection being transverse
and positive. The fact that $\mu_{L_0}$ is well defined and that it is indeed an
isomorphisms follows directly from Proposition~\ref{thm:homologous}.

We can now define the Maslov index of a curve in $\Lambda$.

\begin{defmaslov}\label{thm:defmaslov}
Let $\ell$ be any continuous curve in $\Lambda$ with endpoints outside
$\Lambda_{\ge1}(L_0)$. The {\em Maslov index\/} of $\ell$ (relatively 
to $L_0$) is the
value of $\mu_{L_0}$ in the homology class of $\ell$. The Maslov index of
$\ell$ will be denoted by $\mu_{L_0}(\ell)$.

If $\ell:[a,b]\mapsto\Lambda$ is any continuous curve such that 
$\{t\in\,]a,b[:\ell(t)\in\Lambda_{\ge1}(L_0)\}$ is contained in some
closed interval $[c,d]\subset\,]a,b[$, then the {\em Maslov index\/}
 $\mu_{L_0}(\ell)$ of $\ell$ is defined to be the Maslov index of the 
restriction of $\ell$ to any such $[c,d]$.
\end{defmaslov}
Since the homology class of a concatenation of curves is equal
to the sum of their homology classes, and since $\mu_{L_0}$ is
a group homomorphism, it follows that the Maslov index of curves
is additive by concatenation.
Moreover, Proposition~\ref{thm:homologous}
gives us the following geometrical interpretation of the
Maslov index of a curve. If $\ell$ is a smooth curve in $\Lambda$, with
endpoints outside $\Lambda_{\ge1}(L_0)$ and having only transverse
intersections with $\Lambda_{\ge1}(L_0)$,  
then the Maslov index of $\ell$ is the number of positive intersections
minus the number of negative intersections of $\ell$ with $\Lambda_{\ge1}(L_0)$.

If either one of the endpoints of $\ell$ do belong to $\Lambda_{\ge1}(L_0)$,
our definition of Maslov index simply says that these intersections are
not counted. 
\mysubsection{Computation of the Maslov index}
\label{sub:computation}

We now show how to compute the Maslov index of a curve having image 
entirely contained in the domain of a fixed chart.
\begin{method}\label{thm:method}
Let $\ell:[a,b]\mapsto\Lambda$ be any continuous curve with endpoints outside
$\Lambda_{\ge1}(L_0)$. If there exists a Lagrangian subspace $L_1$ complementary to
$L_0$ and such that the image of $\ell$ is entirely contained in the
domain $\Lambda_0(L_1)$ of the chart $\phi_{L_0,L_1}$, then:
\begin{equation}\label{eq:eqmaslov}
\mu_{L_0}(\ell)=n_+(\beta(b))-n_+(\beta(a)),
\end{equation}
where $\beta=\phi_{L_0,L_1}\circ\ell$.
\end{method}
\begin{proof}
We start observing that, by Lemma~\ref{thm:preliminary},
the Maslov index $\mu_{L_0}(\ell)$ depends only on the
numbers $n_+(\beta(b))$ and $n_+(\beta(a))$.

To prove the statement, it suffices to exhibit for each $i,j=0,\ldots,n$
a curve $\beta_{i,j}:[a,b]\mapsto {\rm B}_{\rm sym}(L_0,\R)$, such that
$n_+(\beta_{i,j}(a))=i$, $n_+(\beta_{i,j}(b))=j$, and such that the curve
$\ell_{i,j}=\phi_{L_0,L_1}^{-1}\circ\beta_{i,j}$ has Maslov index equal
to $j-i$. Clearly, since we can consider curves reparameterized
backwards, it suffices to consider the case $i\le j$. For $i=j$, a 
constant curve with positive type number equal to $i$ would do the job.
It is indeed sufficient to exhibit curves $\beta_{i,i+1}$ as above
for all $i=0,\ldots,n-1$. To prove this claim, observe in first place
that, if such a curve $\beta_{i,i+1}$ is found and 
$\tilde\beta_{i,i+1}$ is any other curve having the same positive type
numbers at the endpoints, then the corresponding curves in $\Lambda$
have the same Maslov index. Now, if the curves $\tilde\beta_{i,i+1}$ are chosen 
in such a way that the endpoint of $\tilde\beta_{i,i+1}$
coincides with the initial point of $\tilde\beta_{i+1,i+2}$, then
the concatenation $\beta_{i,j}=\tilde\beta_{i,i+1}\diamond\tilde\beta_{i+1,i+2}
\diamond\cdots\diamond\tilde\beta_{j-1,j}$ has the desired properties.

To complete the proof, we now show how to construct the curves $\beta_{i,i+1}$
as above. Choose any basis of $L_0$ and define a curve
$\beta_{i,i+1}:[-1,1]\mapsto{\rm B}_{\rm sym}(L_0,\R)$ 
such that $\beta_{i,i+1}(t)$ is given in the chosen basis by
the diagonal $n\times n$ matrix having diagonal vector $(\underbrace{1,\ldots,1}_{
\scriptstyle i\ \rm times},t,-1,\ldots,-1)$. Let $\ell_{i,i+1}=\phi_{L_0,L_1}^{-1}\circ
\beta_{i,i+1}$; we need to show that $\mu_{L_0}(\ell_{i,i+1})=1$.
It is easy to see that every $\ell_{i,i+1}$ intersects $\Lambda_{\ge1}(L_0)$ only
once at $t_0=0$, and that $\ell_{i,i+1}(0)\cap L_0={\rm Ker}(\beta_{i,i+1}(0))$.
Since $\beta'_{i,i+1}(0)$ is positive definite on the one dimensional
space ${\rm Ker}(\beta_{i,i+1})$, the intersection of $\ell_{i,i+1}$
with $\Lambda_{\ge1}(L_0)$ is transverse and positive (see
Proposition~\ref{thm:geometriaLambdak} and Remark~\ref{thm:remtangent}). 
By definition, the Maslov index of $\ell_{i,i+1}$ is equal to $1$, and
we are done.
\end{proof}
It is now easy to prove the following estimate for the Maslov index:
\begin{estimate}\label{thm:estimate}
Let $\ell:[a,b]\mapsto\Lambda$ be any continuous curve with endpoints outside
$\Lambda_{\ge1}(L_0)$. Then,
\begin{equation}\label{eq:estimate}
\vert\mu_{L_0}(\ell)\vert\le \sum_{t\in\,]a,b[}\mathrm{dim}\big(\ell(t)\cap L_0\big).
\end{equation}
\end{estimate}
\begin{proof}
If there are infinitely many $t\in\,]a,b[$ such that $\ell(t)\in\Lambda_{\ge1}(L_0)$,
then the right hand side of \eqref{eq:estimate} is infinite, and the statement
of the Corollary is trivial. 
Otherwise, let $t_0\in\,]a,b[$ be such that $\ell(t_0)\in\Lambda_{\ge1}(L_0)$ and
let $L_1\in\Lambda$ be a Lagrangian complementary to both $L_0$ and $\ell(t_0)$
(see Corollary~\ref{thm:complcomun}). Set $\beta=\phi_{L_0,L_1}\circ\ell$;
then, $\beta$ is a curve in ${\rm B}_{\rm sym}(L_0,\R)$ defined in a neighborhood of 
$t_0$. It is easily seen that ${\rm Ker}(\beta(t_0))=\ell(t_0)\cap L_0$.
By elementary arguments, we have that, for $t$ sufficiently close to
$t_0$, the following inequality holds:
\[n_+(\beta(t_0))\le n_+(\beta(t))\le n_+(\beta(t_0))+\dgn(\beta(t_0)).\]
Hence, for $\varepsilon>0$ small enough, we get:
\[\vert n_+(\beta(t_0+\varepsilon))-n_+(\beta(t_0-\varepsilon))\vert\le \dgn(\beta(t_0)). \]
The conclusion follows easily from Proposition~\ref{thm:method}.
\end{proof}
Under a non degeneracy assumption, the Maslov index can be computed
as a sum of signatures:
\begin{signature}\label{thm:signature}
Let $\ell:[a,b]\mapsto\Lambda$ be a curve of class $C^1$ having endpoints
outside $\Lambda_{\ge1}(L_0)$. If for all $t\in\,]a,b[$ such 
that $\ell(t)\in\Lambda_{\ge1}(L_0)$  we have that $\ell'(t)$ is non degenerate
on $\ell(t)\cap L_0$, then the number of intersections of $\ell$ with
$\Lambda_{\ge1}(L_0)$ is finite, and:
\begin{equation}\label{eq:eqmaslov1}
\mu_{L_0}(\ell)=\sum_{t\in\,]a,b[}\sgn\big(\ell'(t)\big\vert_{\ell(t)\cap L_0}\big).
\end{equation}
\end{signature}
\begin{proof}
Let $t_0\in\,]a,b[$ be such that $\ell(t_0)\in\Lambda_{\ge1}(L_0)$ and
let $L_1\in\Lambda$ be a Lagrangian complementary to both $L_0$ and $\ell(t_0)$
(see Corollary~\ref{thm:complcomun}). Set $\beta=\phi_{L_0,L_1}\circ\ell$;
then, $\beta$ is a curve in ${\rm B}_{\rm sym}(L_0,\R)$ defined in a neighborhood of 
$t_0$. It is easily seen that ${\rm Ker}(\beta(t_0))=\ell(t_0)\cap L_0$, and
it follows from Remark~\ref{thm:remtangent} that $\beta'(t_0)$ and $\ell'(t_0)$ coincide in
$\ell(t_0)\cap L_0$. 

Applying Lemma~\ref{thm:elementary} around $t_0$, once to $\beta$ and again to a backwards 
reparameterization of $\beta$, we conclude that if $\varepsilon>0$ is small
enough, then $\beta(t)$ is non degenerate for $t\in[t_0-\varepsilon,t_0+\varepsilon]
\setminus\{t_0\}$, and that:
\begin{eqnarray}
\label{eq:1st} 
&& n_+(\beta(t_0+\varepsilon))=n_+\big(\beta(t_0))+n_+(\beta'(t_0)\big\vert_{
{\rm Ker}(\beta(t_0))}\big),
\\
\label{eq:2nd}
&& n_+(\beta(t_0-\varepsilon))=n_+\big(\beta(t_0))+n_-(\beta'(t_0)\big\vert_{%
{\rm Ker}(\beta(t_0))}\big).
\end{eqnarray}
Subtracting \eqref{eq:2nd} from \eqref{eq:1st}, we get
\begin{equation}\label{eq:3rd}
n_+(\beta(t_0+\varepsilon))-n_+(\beta(t_0-\varepsilon))=\sgn\big(\beta'(t_0)\big\vert_{%
{\rm Ker}(\beta(t_0))}\big).
\end{equation}
We have proven that the intersection of $\ell$ with $\Lambda_{\ge1}(L_0)$
at $t_0$ is isolated, and, using Proposition~\ref{thm:method}, it follows
from \eqref{eq:3rd} that:
\[\mu_{L_0}\big(\ell\big\vert_{[t_0-\varepsilon,t_0+\varepsilon]}\big)=
\sgn\big(\ell'(t_0)\big\vert_{\ell(t_0)\cap L_0}\big).\]
The conclusion follows from the additivity of the Maslov index with respect to
concatenation. 
\end{proof}

Obvious modifications can be made to the statements of Proposition~\ref{thm:method}
and Corollaries~\ref{thm:estimate} and \ref{thm:signature} to adapt them to the case of curves
$\ell$ with endpoints in $\Lambda_{\ge1}(L_0)$ for which the Maslov index $\mu_{L_0}(\ell)$ is
defined.

\end{section}
\begin{section}{Applications of the Maslov Index:\\ Stability of the Geometric Index}
\label{sec:applications}
In this Section we apply the abstract theory developed in Sections~\ref{sec:geometry}
and \ref{sec:intersection} to the study of the indexes of the quadruples introduced
in Section~\ref{sec:preliminaries}.
\mysubsection{The Maslov index of a differential problem}
\label{sub:maslovdiff}

Let $(g,R,P,S)$ be an admissible quadruple for the differential problem in $\R^n$;
we recall that, associated to $(g,R,P,S)$, we have constructed
a symplectic form $\omega$ in $\R^{2n}$ and, for each $t\in[a,b]$ a
symplectomorphism $\Psi(t)$ of $(\R^{2n},\omega)$ (Definition~\ref{thm:defomega} and
equation~\eqref{eq:defPsit}). 
The map $t\mapsto\Psi(t)$ is a curve
of class $C^1$ in ${\rm Sp}(\R^{2n},\omega)$.

\noindent\ \ 
Rewriting equation~\eqref{eq:DE} as a first order linear system, we get the following
Cauchy problem satisfied by $\Psi$:
\begin{equation}\label{eq:DEbis}
\Psi'(t)=H(t)\,\Psi(t),\quad\Psi(0)={\rm Id},
\end{equation}
where $H(t):\R^{2n}\mapsto\R^{2n}$ is the linear map defined for all $t\in[a,b]$ by:
\begin{equation}\label{eq:defH}
\phantom{\quad\forall\,x,y\in\R^n.}H(t)[(x,y)]=(y,R(t)[x]),\quad\forall\,x,y\in\R^n.
\end{equation}
Since $R(t)$ is $g$-symmetric, $t\mapsto H(t)$ defines a continuous curve
in ${\rm sp}(\R^{2n},\omega)$; equation~\eqref{eq:DEbis} says that  $\Psi'(t)$ is equal to the
evaluation at $\Psi(t)$ of the right invariant vector field determined by $H(t)$ in ${\rm Sp}(
\R^{2n},\omega)$.\smallskip

We now define the following Lagrangian subspaces of $(\R^{2n},\omega)$:
\begin{equation}\label{eq:ell0}
\ell_0=\big\{(x,y)\in\R^{2n}:x\in P,\ y+S[x]\in P^\perp\big\},
\end{equation}
and
\begin{equation}\label{eq:L0}
L_0=\{0\}\oplus\R^n.
\end{equation}
Observe that $\ell_0$ is the subspace of $\R^{2n}$ determined by the
initial conditions \eqref{eq:IC};
the fact that it is a Lagrangian subspace is proven in formula
\eqref{eq:ell0Lagr}. We also define the $C^1$-curve $\ell:[a,b]\mapsto\Lambda$ by:
\begin{equation}\label{eq:defellt}
\ell(t)=\Psi(t)[\ell_0].
\end{equation}
The crucial observation here is that the curve
$\ell$ intercepts $\Lambda_{\ge1}(L_0)$ at $t_0\in\,]a,b]$ if and only if $t_0$ is
a $(P,S)$-focal instant. Observe also that $\ell(a)\in\Lambda_{\ge1}(L_0)$,
unless $P=\R^n$.
More in general, if $t_0\in[a,b]$, then we have:
\begin{equation}\label{eq:interse}
\ell(t_0)\cap L_0=\Big\{(0,\Jt'(t_0)):\Jt\in\Jts,\
\Jt(t_0)=0\Big\}=\{0\}\oplus\Jts[t_0]^\perp.
\end{equation}
The last equality in \eqref{eq:interse} follows easily from \eqref{eq:identity},
arguing as in the proof of Proposition~\ref{thm:isolated}. 
Recalling \eqref{eq:dimcodim}, we have that, for $t_0\in\,]a,b]$, 
$\ell(t_0)\in\Lambda_k(L_0)$ if and only if $t_0$
is a $(P,S)$-focal instant of multiplicity $k$.
\smallskip

It follows from Proposition~\ref{thm:LieAl} and formula~\eqref{eq:DEbis}
that $\ell$ satisfies the following Cauchy problem:
\begin{equation}\label{eq:DEell}
\ell'(t)=H(t)^*\,\big[\ell(t)\big],\quad\ell(a)=\ell_0;
\end{equation}
where $H^*$ is the vector field introduced in Definition~\ref{thm:defH*}.
In the notation of Proposition~~\ref{thm:LieAl}, the curve $\ell$
is equal to $\kappa_{\ell_0}\circ\Psi$.
\smallskip

By Proposition~\ref{thm:isolated}, $t_0=a$ is an isolated intersection
of $\ell$ with $\Lambda_{\ge1}(L_0)$, and therefore we can give the following
definition:
\begin{MaslovDiff}\label{thm:MaslovDiff}
Let $(g,R,P,S)$ be an admissible quadruple for the differential problem
such that the final instant $t_0=b$ is not $(P,S)$-focal.
Its {\em Maslov index\/} $\mu(g,R,P,S)$ is the Maslov
index $\mu_{L_0}(\ell)$, where $L_0$
is the Lagrangian defined in \eqref{eq:L0} and $\ell$ is 
the curve defined in \eqref{eq:defellt}.
\end{MaslovDiff}
The condition that $t_0=b$ is not $(P,S)$-focal means that the curve $\ell$
does not intersect $\Lambda_{\ge1}(L_0)$ at its final endpoint; possibly, one
could extend the definition of Maslov index for quadruples where $t_0=b$
is an isolated $(P,S)$-focal instant. \smallskip

We have the following relation between the Maslov index and the focal index
of a quadruple $(g,R,P,S)$:
\begin{MaslovFocal}\label{thm:MaslovFocal}
Let $(g,R,P,S)$ be an admissible quadruple for the differential
problem in $\R^n$ such that:
\begin{enumerate}
\item\label{itm:1} $t_0=b$ is not a $(P,S)$-focal instant;
\item\label{itm:2} for every $(P,S)$-focal instant $t_0\in\,]a,b[$, the restriction
of $g$ to $\Jts[t_0]$ is non degenerate.
\end{enumerate}
Then, the focal index ${\rm i}_{\rm foc}$ of $(g,R,P,S)$ is well defined, and
it equals the Maslov index:
\begin{equation}\label{eq:maslovfocal}
{\rm i}_{\rm foc}=\mu(g,R,P,S).
\end{equation}
\end{MaslovFocal}
\begin{proof}
Using Proposition~\ref{thm:isolated}, hypothesis~\ref{itm:2} implies that the number
of $(P,S)$-focal instants is finite, and so ${\rm i}_{\rm foc}$ is well defined
(see Definition~\ref{thm:defsigdiff}). Now, using
equations~\eqref{eq:defH}, \eqref{eq:interse}, \eqref{eq:DEell}
and Definitions~\ref{thm:defomega} and \ref{thm:defH*}, we compute as follows:
\begin{equation}\label{eq:conta}
\begin{split}
\ell'(t_0)[(0,x),(0,y)]&=H(t_0)^*\,\big[\ell(t_0)\big][(0,x),(0,y)]=\\&=\omega\big[
H(t_0)[(0,x)],(0,y)\big]=\\ &= \omega\big[(x,R(t_0)[0]),(0,y)\big]=g(x,y),
\end{split}
\end{equation}
for all $t_0\in[a,b]$ and 
for all $x,y\in\R^n$ such that the pairs $(0,x)$ and $(0,y)$ belong to $\ell(t_0)\cap L_0=
\{0\}\oplus\Jts[t_0]^\perp$, i.e., for all $x,y\in\Jts[t_0]^\perp$.

By hypothesis~\ref{itm:2} and equation~\eqref{eq:conta}, $\ell'(t_0)$ is non degenerate
on $\ell(t_0)\cap L_0$; moreover, we have:
\[\sgn\big(\ell'(t_0)\big\vert_{\ell(t_0)\cap L_0}\big)=
\sgn\big(g\big\vert_{\Jts[t_0]^\perp}\big).\]
The conclusion follows from Corollary~\ref{thm:signature}.
\end{proof}
Observe that Theorem~\ref{thm:MaslovFocal} gives also an alternative proof 
of Proposition~\ref{thm:isolated}, as Corollary~\ref{thm:signature}
guarantees that the number of intersections of $\ell$ with $\Lambda_{\ge1}(L_0)$
is finite.\smallskip

We now apply the above result to Riemannian or causal Lorentzian geodesics,
obtaining the following:

\begin{MaslovCausal}\label{thm:MaslovCausal}
Let $(\M,\mathfrak  g, {\mathcal P},\gamma)$ be an admissible quadruple
for the geometric problem such that $\gamma(b)$ is not
a ${\mathcal P}$-focal point. Assume $\Mg$ is Riemannian or Lorentzian,
and in the latter case, that $\gamma$ is non spacelike. 
Let $(g,R,P,S)$ be any associated quadruple to $(\M,\mathfrak  g, {\mathcal P},\gamma)$.
Then, the geometric index of $\gamma$ equals the Maslov index of $(g,R,P,S)$:
\begin{equation}\label{eq:GeoMaslov}
{\rm i}_{\rm geom}(\gamma)=\mu(g,R,P,S).
\end{equation}
\end{MaslovCausal}
\begin{proof}
By Proposition~\ref{thm:matatudo}, since $\gamma(b)$ is not $\mathcal P$-focal,
then $t_0=b$ is not a $(P,S)$-focal instant. The conclusion follows at once
from Proposition~\ref{thm:matatudo}, Remark~\ref{thm:unifying}, 
Corollary~\ref{thm:riemlorcausal} and  Theorem~\ref{thm:MaslovFocal}.
\end{proof}
\mysubsection{Stability of the indexes}
\label{sub:stability}
We now want to study the stability of the Maslov and the focal index
for the differential problem and for the geometrical problem. 
We begin by introducing a
notion of convergence for quadruples $(g,R,P,S)$; in particular, we 
will describe the topological structure of the set of pairs $(P,S)$
as a suitable fiber bundle.

For $k=0,\ldots,n$, let ${\rm GB}_k(n,\R)$ be the set of pairs $(P,S)$, where
$P\subset\R^n$ is a $k$-dimensional subspace and $S\in{\rm B}_{\rm sym}(P,\R)$. 
We define in ${\rm GB}_k(n,\R)$ the structure of a vector bundle over the Grassmannian
$G_k(\R^n)$, whose fiber over $P\in G_k(\R^n)$ is the vector space ${\rm B}_{\rm sym}(P,\R)$.
To define local  trivializations of ${\rm GB}_k(n,\R)$ we argue as follows.
Let $\R^n=W_0\oplus W_1$ be a direct sum decomposition, where $W_0$ is a $k$-dimensional
subspace. As in Section~\ref{sec:geometry}, we define a chart in $G_k(\R^n)$ by associating
to each $P\in G_k(\R^n)$ transverse to $W_1$ the only linear map $T:W_0\mapsto W_1$
whose graph in $W_0\oplus W_1=\R^n$ is $P$. Then, a local trivialization
of ${\rm GB}_k(n,\R)$ is defined by mapping each $S\in{\rm B}_{\rm sym}(P,\R)$
to the bilinear map $\hat T_*(S)=S(\hat T^{-1}\cdot,\hat T^{-1}\cdot)\in{\rm B}_{\rm
sym}(W_0,\R)$,
where the isomorphism $\hat T:W_0\mapsto P$ is given by $v\mapsto v+T(v)$.

As to the geometrical problem, we now define the following space.
Let $M$ be any smooth manifold of dimension $m$ and let $k=0,\ldots,m$ be fixed.
We denote by ${\rm GB}_k(M)$ the set of triples $(p,P,S)$, where $p$ is a point of
$M$, $P$ is a $k$-dimensional subspace of $T_pM$ and $S$ is a symmetric
bilinear form on $P$. The space ${\rm GB}_k(M)$ has an obvious structure
of a fiber bundle over $M$ with projection $(p,P,S)\mapsto p$; namely, any
local trivialization of the tangent bundle $TM$ around $p_0\in M$ induces
a bijection from the fiber of ${\rm GB}_k(M)$ over $p$ (in a neighborhood of
$p_0$) and the manifold ${\rm GB}_k(m,\R)$. These bijections give a local trivialization
of ${\rm GB}_k(M)$ around $p_0$. Observe that, in the case $k=0$, the typical fiber
${\rm GB}_0(m,\R)$ reduces to a point and the fiber bundle ${\rm GB}_0(M)$ is diffeomorphic
to $M$.

The notion of convergence in the bundles ${\rm GB}_k(n,\R)$ and ${\rm GB}_k(M)$
can be described in elementary terms, using convergence of linear basis and
matrices. Namely, a sequence $(P_j,S_j)$ in ${\rm GB}_k(n,\R)$ converges to
$(P,S)$ if and only if for each $j\in\N$ there exists a basis $\{e^j_1,\ldots,e^j_k\}$
of $P_j$ such that $e^j_i\to e_i$ as $j\to\infty$ for all $i$, with $\{e_1,\ldots,e_k\}$
a basis of $P$, and such that $S_j(e^j_\alpha,e^j_\beta) \to S(e_\alpha,e_\beta)$
as $j\to\infty$ for all $\alpha,\beta=1,\ldots,k$.

The convergence $(p_j,P_j,S_j)$ to $(p,P,S)$ in ${\rm BG}_k(M)$ is equivalent
to the convergence of $p_j$ to $p$ in $M$ and to the convergence of $(P_j,S_j)$
to $(P,S)$ in ${\rm GB}_k(m,\R)$, when one considers a local trivialization
of the tangent bundle $TM$ around $p$.

Alternatively, the manifold structure of ${\rm GB}_k(n,\R)$ and ${\rm GB}_k(M)$
can be described in terms of principle and associated bundles.%
\footnote{%
The open subset ${\mathcal L}_{\rm inj}(\R^k,\R^n)$ in the vector space ${\mathcal L}(\R^k,\R^n)$
consisting of injective linear maps is the total space of a ${\rm GL}(k,\R)$-principal bundle
over ${\rm G}_k(\R^n)$. Namely, the projection is given by $T\mapsto T(\R^k)$ and the action of
${\rm GL}(k,\R)$ on ${\mathcal L}_{\rm inj}(\R^k,\R^n)$ is given by composition on the right.
Moreover, we have an action of ${\rm GL}(k,\R)$ on the left on ${\rm B}_{\rm sym}(\R^k,\R)$ given
by $\varphi\cdot B=B(\varphi^{-1}.,\varphi^{-1}.)$. It's easily seen that the associated bundle
obtained from this principal bundle and this action is (isomorphic to) the vector bundle ${\rm
GB}_k(n,\R)$. The fiber bundle ${\rm GB}_k(M)$ can also be seen as an associated bundle to the
${\rm GL}(m,\R)$-principal bundle of referentials in $M$ and to the action of ${\rm GL}(m,R)$ on
the manifold ${\rm GB}_k(m,\R)$ on the left defined in the obvious way.}

\smallskip

We can now prove the
following results about  the stability of the Maslov index in the
differential problem and of the focal index in the geometrical problem:
\begin{LimitMaslov}\label{thm:LimitMaslov}
For each $j\in\N\cup\{\infty\}$, let $(g_j,R_j,P_j,S_j)$ be an 
admissible quadruple for the differential problem in $\R^n$. 

Assume that $(g_j,R_j,P_j,S_j)$ {\em
tends\/} to $(g_\infty,R_\infty,P_\infty,S_\infty)$ as
$j\to\infty$, in the following sense:
\begin{enumerate}
\item ${\rm dim}(P_j)=k$ for all $j=1,\ldots,\infty$;
\smallskip
\item $(P_j,S_j)\mapsto (P_\infty,S_\infty)$ in ${\rm GB}_k(n,\R)$ as $j\to\infty$;
\smallskip
\item $g_j\mapsto g_\infty$ in ${\rm B}_{\rm sym}(\R^n,\R)$ as $j\to\infty$;
\smallskip
\item $R_j\mapsto R_\infty$ uniformly on $[a,b]$ as $j\to\infty$.
\end{enumerate}
If $t_0=b$ is not $(P_\infty,S_\infty)$-focal for $(g_\infty,R_\infty,P_\infty,S_\infty)$, 
then, for $j\in\N$ sufficiently large, $t_0=b$ is not $(P_j,S_j)$-focal for  
$(g_j,R_j,P_j,S_j)$, and:
\[\mu(g_j,R_j,P_j,S_j)=\mu(g_\infty,R_\infty,P_\infty,S_\infty).\]
\end{LimitMaslov}
\begin{proof}
For each $j=1,\ldots,\infty$, define the {\em objects\/} $\Psi_j$, $H_j$, $(\ell_0)_j$ and 
$\ell_j$  relative to the quadruple $(g_j,R_j,P_j,S_j)$ as in formulas \eqref{eq:defPsit},
\eqref{eq:defH}, \eqref{eq:ell0} and \eqref{eq:defellt} respectively.
A simple calculation using the charts described for ${\rm GB}_k(n,\R)$ and
for $G_n(\R^{2n})$ shows that $(\ell_0)_j\to(\ell_0)_\infty$ in $G_n(\R^{2n})$, and
therefore in $\Lambda$.

Obviously, $H_j$ tends to $H_\infty$ uniformly on $[a,b]$; by standard results about 
the continuous dependence on the data for ordinary differential equations, from
\eqref{eq:DEbis} we get that $\Psi_j$ tends to $\Psi_\infty$ uniformly (actually, in the
$C^1$-topology) as $j\to\infty$.

By the continuity of the action of ${\rm Sp}(\R^{2n},\omega)$ in $\Lambda$,
it follows that $\ell_j$ tends to $\ell_\infty$ in the compact-open topology.
Since $\ell_\infty(b)\in\Lambda_0(L_0)$ and $\Lambda_0(L_0)$ is open in $\Lambda$, 
we have that $\ell_j(b)\in\Lambda_0(L_0)$, i.e.,
$t_0=b$ is not $(P_j,S_j)$-focal for $(g_j,R_j,P_j,S_j)$ for $j\in\N$ sufficiently large.

It is not hard to prove (see Remark~\ref{thm:uniformity} below)
that there exists an $\varepsilon>0$ such that there are no $(P_j,S_j)$-focal instants
on the interval $]a,a+\varepsilon]$ relatively to the quadruple $(g_j,R_j,P_j,S_j)$,
for all $j=1,\ldots,\infty$.
Hence, the curve $\ell_j$ does not intercept
$\Lambda_{\ge1}(L_0)$ in the interval $\,]a,a+\varepsilon\,]$, 
for all $j=1,\ldots,\infty$.
The Maslov index $\mu(g_j,R_j,P_j,S_j)$ is by definition equal to the Maslov
index $\mu_{L_0}$ of the restriction of $\ell_j$ to $[a+\varepsilon,b]$.

Since $\Lambda$ and $\Lambda_0(L_0)$ are locally path connected, and
since $\Lambda$ is locally simply connected, the convergence of $\ell_j$ to
$\ell_\infty$ (over the interval $[a+\varepsilon,b]$) in the compact-open topology 
implies that, for $j\in\N$ sufficiently large, $\ell_j$ is homotopic to $\ell_\infty$
through curves with endpoints outside $\Lambda_{\ge1}(L_0)$.

Therefore, $\mu_{L_0}(\ell_j)=\mu_{L_0}(\ell_\infty)$ for $j\in\N$ large enough, and we are done.
\end{proof}
\begin{LimitCausal}\label{thm:LimitCausal}
Let $\Mg$ be a Riemannian or Lorentzian manifold; for each $j\in\N\cup\{\infty\}$
let ${\mathcal P}_j$ be a $k$-dimensional smooth submanifold of $\M$ and 
let $\gamma_j:[a,b]
\mapsto\M$ be a non constant geodesic in $\M$, with $\gamma_j(a)\in{\mathcal P}_j$
and $\gamma'_j(a)\in T_{\gamma_j(a)}{\mathcal P}_j^\perp$. 

If $\Mg$ is Lorentzian, we also assume that  
$\gamma_j$ is non spacelike and that $\gamma_j'(a)\not\in T_{\gamma_j(a)}{\mathcal P}_j$
for all $j=1,\ldots,\infty$.

Let ${\mathcal S}_j$ denote the second fundamental form of ${\mathcal P}_j$
at $\gamma_j(a)$ in the normal direction $\gamma_j'(a)$.
Suppose that 
\begin{itemize}
\item $\lim\limits_{j\to\infty}\gamma'_j(a)=\gamma'_\infty(a)$  in $T\M$, 
\item $\lim\limits_{j\to\infty}(\gamma_j(a),T_{\gamma_j(a)}{\mathcal P}_j,{\mathcal S}_j)=
(\gamma_\infty(a),T_{\gamma_\infty(a)}{\mathcal P}_\infty,{\mathcal S}_\infty)$ in
 ${\rm GB}_k(\M)$.
\end{itemize}
Then, if $\gamma_\infty(b)$ is not ${\mathcal P}_\infty$-focal, it follows that,
for $j\in\N$ sufficiently large, $\gamma_j(b)$ is not ${\mathcal P}_j$-focal,
and the geometrical index of $\gamma_j$ relative to ${\mathcal P}_j$ is equal
to the geometrical index of $\gamma_\infty$ relative to ${\mathcal P}_\infty$:
\[{\rm i}_{\rm geom}(\gamma_j)={\rm i}_{\rm geom}(\gamma_\infty),\quad\forall\,j\gg0.\]
\end{LimitCausal}
\begin{proof}
We choose a local trivialization of the tangent bundle $T\M$ around
$\gamma_\infty(a)$ by linearly independent smooth vector fields
$X_1,\ldots,X_m$. Since $\gamma_j(a)\to\gamma_\infty(a)$ as $j\to\infty$,
we can assume without loss of generality that $\gamma_j(a)$ is in the domain
of the $X_i$'s, for all $j=1,\ldots,\infty$.

Now, we trivialize the tangent bundle along each $\gamma_j$,
$j=1,\ldots,\infty$, by considering the parallel transport of
the vectors $X_i(\gamma_j(a))$ along $\gamma_j$. Associated to these trivializations,
we  produce quadruples $(g_j,R_j,P_j,S_j)$ admissible for the differential problem 
in $\R^m$, $j=1,\ldots,\infty$. We emphasize that we are considering  
trivializations of the {\em entire\/} tangent bundle along the geodesics $\gamma_j$;
recall Remark~\ref{thm:unifying} for a discussion about this issue.
We also observe that the condition that $\gamma_j'(a)\not\in T_{\gamma_j(a)}{\mathcal P}_j$
implies in particular that $\mathfrak g$ is non degenerate on $T_{\gamma_j(a)}{\mathcal P}_j$.

Clearly, under our hypothesis,  $(g_j,R_j,P_j,S_j)$ tends
to $(g_\infty,R_\infty,P_\infty,S_\infty)$ as $j\to\infty$ in the sense
of Proposition~\ref{thm:LimitMaslov}.

The conclusion follows now easily from Corollary~\ref{thm:MaslovCausal}
and Proposition~\ref{thm:LimitMaslov}.
\end{proof}

\begin{uniformity}\label{thm:uniformity}
Let $(g_\lambda,R_\lambda,P_\lambda,S_\lambda)$ be an admissible
quadruple for the differential problem in $\R^n$, that depends {\em continuously\/} 
on a parameter $\lambda$ varying in a compact topological space. This means
that ${\rm dim}(P_\lambda)=k\in\N$ for all $\lambda$ and that the maps $\lambda\mapsto
g_\lambda\in {\rm B}_{\rm sym}(\R^n,\R)$, 
$(t,\lambda)\mapsto R_\lambda(t)\in{\mathcal L}(\R^n,\R^n)$ and
$\lambda\mapsto (P_\lambda,S_\lambda)\in {\rm GB}_k(\R^n) $ are continuous.
A minor modification in the argument of the proof of Proposition~\ref{thm:isolated}
shows that we can find $\varepsilon>0$ such that there are no $(P_\lambda,S_\lambda)$-focal
instants on the interval $]a,a+\varepsilon]$ for all $\lambda$. Namely, the vector
fields $J_i$ and $\tilde J_i$ appearing in the proof of Proposition~\ref{thm:isolated}
may be chosen to depend continuously on $(t,\lambda)$, and the conclusion follows easily.
\end{uniformity}
\end{section}
\begin{section}{The Spectral Index.\\ Some Remarks on a Possible Extension
of the Morse Index Theorem}\label{sec:spectral}
In this section we will define the {\em spectral index\/} of an admissible
quadruple for the differential problem $(g,R,P,S)$. Such number is related
to the spectral properties of the unbounded operator associated to
the differential equation \eqref{eq:DE} with boundary conditions \eqref{eq:IC}
and $\Jt(b)=0$. Under suitable hypotheses, we will prove that this index
equals the Maslov index of the quadruple. If $(g,R,P,S)$ arises from a Riemannian
or non spacelike Lorentzian geodesic, the equality of the spectral index and the
Maslov index of $(g,R,P,S)$ gives an equivalent form of the classical Morse
index theorem.
\mysubsection{Eigenvalues of the differential problem and the spectral index}
\label{sub:spectral}

Let's fix an admissible  quadruple $(g,R,P,S)$ for the differential problem in $\R^n$;
we will consider the space ${\mathcal H}=L^2([a,b],\R^n)$ of
$\R^n$-valued square-integrable vector fields on $[a,b]$; rather than 
choosing a specific inner product on
$\mathcal H$, we will only regard it as a {\em Hilbertable\/} space, since all our
statements on $\mathcal H$ will only depend on its topological structure.

On $\mathcal H$ we define the following bounded symmetric bilinear form $\hat g$:
\begin{equation}\label{eq:hatg}
\hat g(u,v)=\int_a^bg(u(t),v(t))\;{\rm d}t;
\end{equation}
from the nondegeneracy of $g$ and the fundamental theorem of Calculus of Variations,
it follows easily that $\hat g$ is non degenerate on $\mathcal H$.

Let $\hat R:{\mathcal H}\mapsto{\mathcal H}$ be the bounded linear operator
given by:
\begin{equation}\label{eq:defhatR}
\hat R[v](t)=R(t)[v(t)],\quad t\in[a,b],
\end{equation}
and let $\mathcal A$ be the densely defined unbounded operator given by
\begin{equation}\label{eq:defAcal}
{\mathcal A}=-\ddt+\hat R,
\end{equation}
defined in the domain $D\subset{\mathcal H}$:
\begin{equation}\label{eq:defdomain}
 D=\Big\{u\in C^2([a,b],\R^n):
u(a)\in P,\ u'(a)+S[u(a)]\in P^\perp,\ u(b)=0\Big\}.
\end{equation}
It is easily seen that the operator $\mathcal A$ is $\hat g$-symmetric, 
in the sense that \begin{equation}\label{eq:gsym}
\hat g({\mathcal A}u,v)=\hat g(u,{\mathcal
A}v),\end{equation} for all $u,v\in D$. However, it is in general impossible to choose
a Hilbert space product on $\mathcal H$ that makes $\mathcal A$ symmetric.
Hence, the spectrum of $\mathcal A$ will not in general be real, and for
this reason we need to introduce a complexification of~$\mathcal H$.
Indeed, we need to investigate the holomorphy properties of our
differential problem in order to establish the discreteness of the set
of eigenvalues of $\mathcal A$.

Let ${\mathcal H}^\C$ be the complex Hilbertable space $L^2([a,b],\C^n)$;
we regard ${\mathcal H}$ as a subspace of ${\mathcal H}^\C$. The space
${\mathcal H}^\C$ is a {\em complexification\/} of ${\mathcal H}$, in the
sense that ${\mathcal H}^\C={\mathcal H}\oplus i{\mathcal H}$.
To each subspace ${\mathcal W}\subset {\mathcal H}$ we associate its complexification
${\mathcal W}^\C={\mathcal W}\oplus i{\mathcal W}$, which is the complex
subspace of ${\mathcal H}^\C$ generated by $\mathcal W$.

Moreover, every linear operator on $\mathcal H$ (bounded or unbounded) has a unique
complex linear extension to ${\mathcal H}^\C$. For simplicity, we will maintain
the same notations for linear operators on $\mathcal H$ and their complex
linear extensions to ${\mathcal H}^\C$. In particular, we will consider the
complex linear extension of $\hat R$ to ${\mathcal H}^\C$ and of ${\mathcal A}$
to $D^\C$.

Let $\lambda\in\C$; we consider the eigenvalue problem for $\mathcal A$
in $D^\C$:
\begin{equation}\label{eq:EigProbl}
u''=\hat R[u]-\lambda u.
\end{equation}
Observe that \eqref{eq:EigProbl} is the differential equation \eqref{eq:DE}
corresponding to the quadruple $(g,R_\lambda, P^\C,S)$, where 
\begin{equation}\label{eq:defRC}
R_\lambda(t)=R(t,\lambda)=R(t)-\lambda\cdot{\rm Id}.
\end{equation}
Here, we are considering an obvious extension of the notion
of admissible quadruples for the differential problem to complex spaces.
For such an extension, one identifies the space ${\mathcal L}(\R^k,\R^k)$
with the subspace of ${\mathcal L}_\C(\C^k,\C^k)$ consisting of all the
complex linear operators on $\C^k$ that preserve the real subspace $\R^k$.
For instance, $R(t)$ and $R_\lambda(t)$ are seen as  complex linear
operators on $\C^n$.

In analogy with \eqref{eq:defPsit} and \eqref{eq:defH}, we define complex
linear operators $\Psi_\lambda(t)$ and $H_\lambda(t)$ on $\C^{2n}$, given by:
\begin{equation}\label{eq:defpsiC}
\Psi_\lambda(t)[(u(a),u'(a))]=(u(t),u'(t)),\quad
H_\lambda(t)[(x,y)]=(y,R_\lambda(t)[x]),
\end{equation}
where $u:[a,b]\mapsto\C^n$ is a solution of \eqref{eq:EigProbl} and $x,y\in\C^n$.
Formulas \eqref{eq:defRC} and \eqref{eq:defpsiC} define maps:
\[R:[a,b]\times\C\mapsto{\mathcal L}_\C(\C^n,\C^n),\quad \Psi,H:[a,b]\times\C
\mapsto{\mathcal L}_\C(\C^{2n},\C^{2n});\]
these maps are continuous, and they are holomorphic on the second variable.
For the holomorphy of $\Psi$, we are using well known regularity results
for the solutions of differential equations. Indeed, $\Psi$ and
$\frac{\partial\Psi}{\partial\lambda}$ satisfy
the following Cauchy problems:
\begin{eqnarray}
\label{eq:EDPsi}
&& \frac{\rm d}{{\rm d}t}\Psi(t,\lambda)=H(t,\lambda)\Psi(t,\lambda),\quad\Psi(a,\lambda)={\rm
Id};
\\
\label{eq:EDPsilambda}
&& \frac{\rm d}{{\rm d} t}\frac{\partial\Psi}{\partial\lambda}(t,\lambda)=
\frac{\partial H}{\partial\lambda}(t,\lambda)\Psi(t,\lambda)+H(t,\lambda)\frac{
\partial\Psi}{\partial\lambda}(t,\lambda),\quad \frac{\partial\Psi}{\partial\lambda}
(a,\lambda)=0.
\end{eqnarray}
From \eqref{eq:EDPsi} and \eqref{eq:EDPsilambda}, we see that $\Psi(t,\lambda)$
is differentiable in $t$, the derivative $\frac{\rm d}{{\rm d}t}\Psi(t,\lambda)$
is jointly continuous in the two variables and holomorphic in $\lambda$.

The eigenvalues of $\mathcal A$ in $D^\C$, i.e., the complex numbers
$\lambda$ for which equation \eqref{eq:EigProbl} admits non trivial
solutions in $D^\C$, can be described as the zeroes of a suitable
entire function. For instance, they are the zeroes of the function
\[r(\lambda)={\rm det}(\pi\circ\Psi(b,\lambda)[e_1],\ldots, \pi\circ\Psi(b,\lambda)[e_n]),\]
where $e_1,\ldots,e_n$ is a basis of the vector space $\ell_0\subset\R^{2n}$
defined in \eqref{eq:ell0} and $\pi:\C^{2n}\mapsto\C^n$ is the projection onto the first
$n$ coordinates. Hence, the set of eigenvalues of $\mathcal A$ in $D^\C$
is either $\C$ or a {\em discrete\/} subset of $\C$. We will establish next
that the real eigenvalues of $\mathcal A$ in $D^\C$ (or equivalently in $D$) 
are bounded from below, and so the set of eigenvalues of $\mathcal A$ in $D^\C$
is discrete in  $\C$.

We need the following technical Lemma:
\begin{technical}\label{thm:technical}
Let $Z$ be a finite dimensional real (or complex)
vector space equip\-ped with a positive definite 
inner (or Hermitian) product 
$\langle\cdot,\cdot\rangle$, and let the corresponding norm
be denoted by $\Vert\cdot\Vert$.
Let $u:[a,b]\mapsto Z$ be a $C^1$-function such that $u(b)=0$. 
Then, the following inequality holds:
\begin{equation}\label{eq:techn}
\int_a^b\Vert u(t)\Vert^2\;{\rm d}t\cdot
\int_a^b\Vert u'(t)\Vert^2\;{\rm d}t\ge \frac 14\Vert
u(a)\Vert^4.
\end{equation}
\end{technical}
\begin{proof}
%
Using the Cauchy--Schwarz inequality, we compute easily:
\[\begin{split}
\Vert u(a)\Vert^2&=-\int_a^b\frac{{\rm d}}{{\rm d}t}\langle u(t),u(t)\rangle
\;{\rm d}t=-2\int_a^b\langle u'(t),u(t)\rangle\;{\rm d}t\le\\
&\le2\left(\int_a^b\Vert u(t)\Vert^2\;{\rm d}t\right)^\frac12
\left(\int_a^b\Vert u'(t)\Vert^2\;{\rm d}t\right)^\frac12.
\end{split}
\]
The inequality \eqref{eq:techn} follows easily.
\end{proof}
We can now prove the following:
\begin{eigenvalues}\label{thm:eigenvalues}
The real part of the eigenvalues of $\mathcal A$ in $D^\C$ is bounded from below.
\end{eigenvalues}
\begin{proof}
Since $g$ is non degenerate on $P$, then $P$ and $P^\perp$ are complementary subspaces
in $\R^n$; let $\langle\cdot,\cdot\rangle$ be any positive definite inner product
on $\R^n$ which makes $P$ and $P^\perp$ orthogonal, and denote also by
$\langle\cdot,\cdot\rangle$ its extension to a Hermitian product in
$\C^n$. We denote by
$\langle\cdot,\cdot\rangle_2$ the corresponding Hermitian product in $L^2([a,b],\C^n)$.

Using integration by parts, for all $u\in D^\C$ we have the following:
\begin{equation}\label{eq:p}
\langle-\ddt u,u\rangle_2=\int_a^b\langle u'(t),u'(t)\rangle\;{\rm d}t+
\langle u(a),u'(a)\rangle.
\end{equation}
Moreover, for all $u\in D^\C$, we have:
\begin{equation}\label{eq:s}
\vert \langle u(a),u'(a)\rangle\vert=\vert\langle S[u(a)],u(a)\rangle\vert\le\Vert
S\Vert\cdot
\Vert u(a)\Vert ^2.
\end{equation}
Let $u\in D^\C$ be such that $\Vert u\Vert_2^2=\int_a^b\Vert u(t)\Vert^2\;{\rm d}t=1$;
we apply to such a function Lemma~\ref{thm:technical}, obtaining:
\begin{equation}\label{eq:t}
\int_a^b\Vert u'(t)\Vert^2\;{\rm d}t\ge \frac{\Vert
u(a)\Vert^4}{4},\quad
\forall\,u\in D^\C,\ \Vert u\Vert_2=1.
\end{equation}
Using \eqref{eq:s} and \eqref{eq:t}, we obtain that the right side of \eqref{eq:p}
is bounded from below for $u\in D^\C$ with $\Vert u\Vert_2=1$, i.e.,
there exists $k_0\in\R$ such that:
\begin{equation}\label{eq:q}
\langle-\ddt u,u\rangle_2\ge k_0,\quad\forall u\in D^\C,\ \Vert u\Vert_2=1.
\end{equation}
Let now $\lambda$ be any eigenvalue of $\mathcal A$ in $D$ and $u\in D^\C$
be a corresponding eigenvector with $\Vert u\Vert_2=1$. From \eqref{eq:q}
we compute easily:
\[{\rm Re}(\lambda)={\rm Re}(\langle \lambda u,u\rangle_2)=\langle-\ddt u,u\rangle_2+
{\rm Re}(\langle \hat R[u],u\rangle_2)\ge k_0-\Vert\hat R\Vert_2>-\infty,\]
where $\Vert\hat R\Vert_2$ is the operator norm of $\hat R$
in $L^2([a,b],\C^n)$. This concludes the proof.
\end{proof}
As we have observed previously, the set of eigenvalues of $\mathcal A$ in
$D^\C$ is discrete, and so Proposition~\ref{thm:eigenvalues} gives us
the following corollary:
\begin{finiteeigen}\label{thm:finiteeigen}
The operator $\mathcal A$ has only a finite number of real negative
eigenvalues in $D$.\qed
\end{finiteeigen}
\begin{remuniform}\label{thm:remuniform}
In what follows, we will have to consider the operators ${\mathcal A}_t$ and $\hat R_t$
on $L^2([a,t],\R^n)$, for a fixed $t\in\,]a,b]$,  defined in analogy with \eqref{eq:defhatR}
and \eqref{eq:defAcal} considering the restriction of $R$ to $[a,t]$.
The domain of ${\mathcal A}_t$ is meant to be the subspace $D_t$ defined as in
\eqref{eq:defdomain} by replacing the endpoint $b$ with $t$.

Clearly, Proposition~\ref{thm:eigenvalues} and Corollary~\ref{thm:finiteeigen} 
remain valid for ${\mathcal A}_t$; as a matter of fact, one can choose a lower bound for
the real eigenvalues of ${\mathcal A}_t$ which is independent of $t\in\,]a,b]$.
This can be easily seen by considering that the constant $k_0$ in the inequality
\eqref{eq:q} does not depend on $t$ and that $\Vert\hat R_t\Vert$ is bounded from above
by the supremum norm of $R:[a,b]\mapsto{\mathcal L}(\R^n,\R^n)$. 
\end{remuniform}

From now on, we will disregard the complexified spaces introduced, and
we will only deal with the real eigenvalues of $\mathcal A$ in $D$. 
So, we will look at the maps $H$ and $\Psi$ only in their real domains
and counterdomains:
\[H:[a,b]\times\R\mapsto {\rm sp}(\R^{2n},\omega),\quad
\Psi:[a,b]\times\R\mapsto {\rm Sp}(\R^{2n},\omega),\]
where $\omega$ is the symplectic form of Definition~\ref{thm:defomega}.

Keeping in mind formulas \eqref{eq:ell0}, \eqref{eq:L0} and \eqref{eq:defellt},
we define $\ell:[a,b]\times\R\mapsto\Lambda$ by:
\[\ell(t,\lambda)=\Psi(t,\lambda)[\ell_0].\]
If $\lambda$ is a real eigenvalue of ${\mathcal A}$ in $D$, we denote
by ${\mathcal H}_\lambda={\rm Ker}({\mathcal
A}-\lambda\cdot{\rm Id})\subset D$ the corresponding eigenspace.
We observe that ${\mathcal H}_\lambda$ is the set of $(P,S)$-solutions
relative to the quadruple $(g,R_\lambda,P,S)$ vanishing at $t_0=b$. It follows
that $\lambda\in\R$ is an eigenvalue of $\mathcal A$ if and only if 
$t_0=b$ is a $(P,S)$-focal instant for such a quadruple. Moreover, the dimension
of ${\mathcal H}_\lambda$ coincides with the multiplicity of $t_0=b$ as a $(P,S)$-focal
instant for the quadruple $(g,R_\lambda,P,S)$, and therefore it is finite:
\begin{equation}
\label{eq:dimHlambda}
{\rm dim}({\mathcal H}_\lambda)\le n<+\infty.
\end{equation}

We now look at the Maslov index $\mu_{L_0}$ of the curve $\lambda\mapsto
\ell(b,\lambda)$; we observe that $\ell(b,\lambda)\in\Lambda_{\ge1}(L_0)$ if and
only if $\lambda$ is an eigenvalue of $\mathcal A$ in $D$. Moreover, in analogy
with \eqref{eq:interse}, we have:
\begin{equation}\label{eq:sobre}
\ell(b,\lambda)\cap L_0=\Big\{(0,u'(b)):u\in\Jts_\lambda,\
u(b)=0\Big\}=\{0\}\oplus\Jts_\lambda[b]^\perp,
\end{equation}
where $\Jts_\lambda$ is the space of all $(P,S)$-solutions relative
to the quadruple $(g,R_\lambda,P,S)$.

We can now give the following:
\begin{defindspectral}\label{thm:defindspectral}
The {\em spectral index\/} ${\rm i}_{\rm spec}$ of the quadruple $(g,R,P,S)$ 
is the sum of the signatures of the restrictions of $\hat g$ to the eigenspaces
relative to the negative eigenvalues of $\mathcal A$:
\[{\rm i}_{\rm spec}=\sum_{\lambda<0}\sgn(\hat g \,\vert_{{\mathcal H}_\lambda}).\]
\end{defindspectral}
Observe that, by Corollary~\ref{thm:finiteeigen} and formula \eqref{eq:dimHlambda},
${\rm i}_{\rm spec}$ is a finite integer number.

\mysubsection{A generalized Morse Index Theorem}
\label{sub:morse}

We want to prove that, under suitable hypotheses, ${\rm i}_{\rm spec}$ is
equal to the Maslov index $\mu_{L_0}$ of the curve $\lambda\mapsto\ell(b,\lambda)$.
We start with the following:
\begin{restrictions}\label{thm:restrictions}
Let $\lambda\in\R$ be an eigenvalue of $\mathcal A$ in $D$.
The map ${\mathcal H}_\lambda\ni u\mapsto (0,u'(b))\in\ell(b,\lambda)\cap L_0$
is a linear isomorphism which carries the restriction of $\hat g$ to the restriction
of the symmetric bilinear form $\frac{\partial\ell}{\partial\lambda}(b,\lambda)$.
\end{restrictions}
\begin{proof}
The map $u\mapsto(0,u'(b))$ is clearly injective on ${\mathcal H}_\lambda$, 
and it is onto by \eqref{eq:sobre}.

We compute the derivative $\frac{\partial\ell}{\partial\lambda}$;
recalling Proposition~\ref{thm:LieAl} and Definition~\ref{thm:defH*},
we have:
\[\frac{\partial\ell}{\partial\lambda}(t,\lambda)=\left[\frac{\partial\Psi}{\partial\lambda}
\,\Psi^{-1}\right]^*\,[\ell(t,\lambda)]=\omega(\frac{\partial\Psi}{\partial\lambda}
\,\Psi^{-1}\,.\,,.\,)\big\vert_{\ell(t,\lambda)};\]
the pull-back of $\frac{\partial\ell}{\partial\lambda}(t,\lambda)$ by $\Psi(t,\lambda)$
is a symmetric bilinear form on $\ell_0$ given by:
\begin{equation}\label{eq:pullback}
\frac{\partial\ell}{\partial\lambda}(t,\lambda)[\Psi\;\cdot,\Psi\;\cdot]=
\omega(\frac{\partial\Psi}{\partial\lambda}\;\cdot,\Psi\;\cdot)\big\vert_{\ell_0}.
\end{equation}
We want to calculate the derivative of the pull-back \eqref{eq:pullback}
with respect to $t$. First, we differentiate $\Psi^{-1}\frac{\partial\Psi}{\partial\lambda}$:
\begin{equation}\label{eq:derPsiPsil}
\begin{split}
\frac{\rm d}{{\rm d}t}\,\left[\Psi^{-1}\frac{\partial\Psi}{\partial\lambda}\right]&=
-\Psi^{-1}\frac{\partial\Psi}{\partial t}\Psi^{-1}\frac{\partial\Psi}{\partial\lambda}+
\Psi^{-1}\frac{\partial H}{\partial\lambda}\Psi+\Psi^{-1} H\frac{\partial\Psi}{\partial\lambda}
=\\
&=-\Psi^{-1}H\frac{\partial\Psi}{\partial\lambda}+\Psi^{-1}
\frac{\partial H}{\partial\lambda}\Psi+\Psi^{-1}H\frac{\partial\Psi}{\partial\lambda}=\\
&=\Psi^{-1}
\frac{\partial H}{\partial\lambda}\Psi,
\end{split}
\end{equation}
where in the first equality we have used \eqref{eq:EDPsilambda} and in the second one we 
have used 
\eqref{eq:EDPsi}.
Hence, the derivative of the pull-back \eqref{eq:pullback} is given by:
\begin{equation}\label{eq:derpullback}
\begin{split}
\frac{\rm d}{{\rm d}t}\,\omega(\frac{\partial\Psi}{\partial\lambda}\,.,\Psi\,.)\big\vert_{\ell_0}
&=
\frac{\rm d}{{\rm d}t}\,\omega(\Psi^{-1}\frac{\partial\Psi}{\partial\lambda}\,.,.\,)
\big\vert_{\ell_0}=\omega(\Psi^{-1}\frac{\partial H}{\partial\lambda}\Psi\,.,.\,)
\big\vert_{\ell_0}=\\
&=\omega(\frac{\partial H}{\partial\lambda}\Psi\,.,\Psi\,.\,)
\big\vert_{\ell_0},
\end{split}
\end{equation}
where in the first and in the third equality we have used the fact that $\Psi$ is a
symplectomorphism and in the second one we have used \eqref{eq:derPsiPsil}.

Observe that, by \eqref{eq:defRC} and \eqref{eq:defpsiC}, 
$\frac{\partial H}{\partial \lambda}$ is
simply given by:
\begin{equation}\label{eq:Hl}
\frac{\partial H}{\partial
\lambda}[(x,y)]=(0,-x),\quad\forall\,x,y\in\R^n.
\end{equation} Integration of \eqref{eq:derpullback} on the
interval $[a,b]$  using
\eqref{eq:pullback} gives:
\begin{equation}\label{eq:integration}
\int_a^b\omega(\frac{\partial H}{\partial\lambda}\Psi\,.,\Psi\,.\,)\big\vert_{\ell_0}\;{\rm d}t=
\frac{\partial\ell}{\partial\lambda}(b,\lambda)[\Psi(b,\lambda)\,\cdot,\Psi(b,\lambda)\,\cdot],
\end{equation}
because $\frac{\partial\Psi}{\partial\lambda}(a,\lambda)=0$ (see equation \eqref{eq:EDPsilambda}).

Finally, let $u,v$ be elements of ${\mathcal H}_\lambda$; evaluating
\eqref{eq:integration} in the pairs $(u(a),u'(a))$ and $(v(a),v'(a))$ in $\ell_0$, we
get:
\[
\begin{split}
\int_a^bg(u(t),v(t))\;{\rm d}t&=\frac{\partial\ell}{\partial\lambda}(b,\lambda)
[(u(b),u'(b)),(v(b),v'(b))]=\\
&=\frac{\partial\ell}{\partial\lambda}(b,\lambda)
[(0,u'(b)),(0,v'(b))],
\end{split}
\]
which was obtained by using Definition~\ref{thm:defomega} and formulas \eqref{eq:defpsiC}
and \eqref{eq:Hl}.

This concludes the proof.
\end{proof}
\begin{remerror}\label{thm:remerror}
The computation presented in the proof of Lemma~\ref{thm:restrictions}
also appears in the proof of \cite[Proposition~6.2]{Hel1}. However,
in drawing the final conclusion of the Proposition, the author 
identifies the intersection $\ell(t,\lambda)\cap L_0$ with the
{\em generalized\/} eigenspace $\bigcup_{k\ge1}\mathrm{Ker}(\mathcal A-\lambda)^k$. 
This is clearly {\em not\/} the case if $\lambda$ is not a simple eigenvalue
of $\mathcal A$; it is not clear whether the case of non simple
eigenvalues can be treated by other arguments.
\end{remerror}
\begin{indlambda}\label{thm:indlambda}
Suppose that the restriction of $\hat g$ to ${\mathcal H}_\lambda$
is non degenerate for every {\em negative\/} eigenvalue
$\lambda$ of ${\mathcal A}$ in $D$.
Then, if $\lambda_0<0$ is smaller than the minimum eigenvalue
of $\mathcal A$ in $D$,   the
Maslov index $\mu_{L_0}$ of the curve $[\lambda_0,0]\ni\lambda\mapsto\ell(b,\lambda)$
is equal to the spectral index ${\rm i}_{\rm spec}$.
\end{indlambda}
\begin{proof}
It follows immediately from Corollary~\ref{thm:signature} and Lemma~\ref{thm:restrictions}.
\end{proof}
In the following theorem we relate the spectral index with the Maslov index
of a quadruple $(g,R,P,S)$. Recalling Theorem~\ref{thm:MaslovFocal}, we then
obtain an equality also of the focal and the spectral indexes.
\begin{Morse}\label{thm:Morse}
Let $(g,R,P,S)$ be an admissible quadruple for the differential
problem in $\R^n$; assume that $t_0=b$ is {\em not\/} a $(P,S)$-focal
instant. Suppose that the restriction of $\hat g$ to ${\mathcal H}_\lambda$
is non degenerate for every negative  eigenvalue
$\lambda$ of ${\mathcal A}$ in $D$. Then, the spectral and the Maslov indexes
of $(g,R,P,S)$ coincide:
\[{\rm i}_{\rm spec}=\mu(g,R,P,S).\]
\end{Morse}
\begin{proof}
Let $\lambda_0<0$ be chosen so that  $\ell(t,\lambda)\not\in\Lambda_{\ge1}(L_0)$
for all $t\in\,]a,b]$ and all $\lambda\le\lambda_0$. By Remark~\ref{thm:remuniform}, 
to such purpose it suffices to take $\vert\lambda_0\vert$ large enough.
By Remark~\ref{thm:uniformity}, we can find $\varepsilon>0$ small enough, so that
$\ell(t,\lambda)\not\in\Lambda_{\ge1}(L_0)$ for all $t\in\,]a,a+\varepsilon]$ and for
all $\lambda\in[\lambda_0,0]$. We will consider the restriction of $\ell$ to the
rectangle $[a+\varepsilon,b]\times[\lambda_0,0]$.
Now, the Maslov index $\mu_{L_0}$ of the curve $t\mapsto\ell(t,0)$, $t\in[a+\varepsilon,b]$
is by definition the Maslov index of the quadruple $(g,R,P,S)$; the Maslov
index $\mu_{L_0}$ of the curve $\lambda\mapsto\ell(b,\lambda)$, $\lambda\in[\lambda_0,0]$
is equal to the spectral index, by Corollary~\ref{thm:indlambda}.
Finally, the image by $\ell$ of the remaining two sides of the rectangle
$[a+\varepsilon,b]\times[\lambda_0,0]$ is disjoint from $\Lambda_{\ge1}(L_0)$
by our choice of $\varepsilon$ and $\lambda_0$. The conclusion follows from
the homotopy invariance of the Maslov index.
\end{proof}
\begin{remMorse}\label{thm:remMorse}
Theorem~\ref{thm:Morse} can be seen as a generalization of the classical
Morse Index Theorem in Riemannian or Lorentzian geometry, in the following
sense. Let's assume that $(g,R,P,S)$ is a given admissible quadruple for
the differential problem in $\R^n$, with the property that $g$ is 
{\em positive definite\/}
in $\R^n$. This is the case when $(g,R,P,S)$ arises
from an admissible quadruple $(\M,\mathfrak g,\gamma,{\mathcal P})$, with
$\Mg$ Riemannian or Lorentzian, and in the latter case, with $\gamma$ 
non spacelike (see Proposition~\ref{thm:correspondence} and Remark~\ref{thm:unifying}). 
In the Lorentzian case, the bilinear form
$g$ is positive definite when one considers trivializations of the normal bundle along $\gamma$,
if
$\gamma$ is timelike, and of the quotient bundle $\mathcal N$ for a lightlike
geodesic $\gamma$ (see Section~\ref{sec:preliminaries}). 

For a positive definite $g$, the corresponding bilinear form $\hat g$ is a
Hilbert space inner product on $\mathcal H$; by \eqref{eq:gsym}, $\mathcal A$
is symmetric, and it admits a closed self-adjoint extension ${\mathcal A}_o$
to some suitable Sobolev space $D_o$. 
More precisely, $D_o$ is easily seen to be the space of all $C^1$-functions
$u:[a,b]\mapsto\R^n$ with absolutely continuous derivative and square integrable second
derivative, satisfying the boundary conditions:
\begin{equation}\label{eq:boundary}
u(a)\in P,\quad u'(a)+S[u(a)]\in P^\perp,\quad u(b)=0.
\end{equation}
An explicit integral formula for the resolvent
of ${\mathcal A}_o$,  using the method
of variations of constants, shows that the spectrum $\sigma({\mathcal A}_o)$ coincides
with the set of eigenvalues of $\mathcal A$ in $D$.

By the spectral theorem for (unbounded) self-adjoint operators, we get a direct sum decomposition
\[{\mathcal H}=\bigoplus_{\lambda\in\sigma({\mathcal A}_o)}{\mathcal H}_\lambda,\] 
where ${\mathcal H}_\lambda$ is, as before, the eigenspace of ${\mathcal A}$ corresponding to 
$\lambda$.

We introduce the {\em index form\/} $I$ by:
\begin{equation}\label{eq:dfindform}
I(u,v)=\int_a^b\left[g(u',v')+g(R[u],v)\right]\;{\rm d}t-S[u(a),v(a)],
\end{equation} 
which is a symmetric bilinear form  in $D_o$; we observe that a simple integration by parts
shows that:
\[I(u,v)=\hat g({\mathcal A}_o[u],v)=\int_a^bg(-u''+R[u],v)\;{\rm d}t,
\quad\forall\,u,v\in D_o.\]
From the spectral decomposition it follows easily that the index of the
bilinear form $I$ on $D_o$ (in the sense of Definition~\ref{thm:defsym})
is the sum of the dimensions of ${\mathcal H}_\lambda$ for negative $\lambda$.
As $\hat g$ is positive definite, this number coincides with the spectral
index of $(g,R,P,S)$.

The bilinear form $I$ of formula \eqref{eq:dfindform} can be naturally extended
to a {\em continuous\/} bilinear form on the Hilbert space ${\mathcal H}^1$
consisting of all absolutely continuous functions with square integrable
first derivative and satisfying the boundary conditions~\eqref{eq:boundary}. 
The Hilbert space inner product on ${\mathcal H}^1$ that makes $I$ continuous
is, for instance, $(u,v)\mapsto \hat g(u',v')$. The space $D_o$ is a dense
linear subspace of ${\mathcal H}^1$; since $I$ is continuous in ${\mathcal H}^1$,
a simple density argument shows that the index of $I$ on ${\mathcal H}^1$ is
the same as the index of $I$ on any dense subspace of ${\mathcal H}^1$.
For instance, in the classical proof of the Morse Index Theorem (see~\cite{dC,M}),
one considers the space of piecewise smooth functions.

Now, if $t_0=b$ is not a $(P,S)$-focal instant, Theorem~\ref{thm:Morse}
gives us an equality between the spectral index and the Maslov index of $(g,R,P,S)$;
observe that the non degeneracy assumption of $\hat g$ on each ${\mathcal H}_\lambda$ 
is automatically satisfied.

If $(g,R,P,S)$ is associated to the quadruple $(\M,\mathfrak g,\gamma,{\mathcal P})$,
where $\Mg$ is either Riemannian or Lorentzian, and in the latter case, $\gamma$
is non spacelike, then, under the assumption that $\gamma(b)$ is not a ${\mathcal P}$-focal
point, Corollary~\ref{thm:MaslovCausal} gives the equality between the geometrical
index of $\gamma$ and the Maslov index of $(g,R,P,S)$, i.e., the index of
the bilinear form $I$ on the space ${\mathcal H}^1$ (or any of its dense subspaces).

It is not difficult to prove that the index of $I$ on ${\mathcal H}^1$ is equal
to the index of the {\em second variation\/} of the {\em energy functional\/}
on the set of curves connecting the submanifold ${\mathcal P}$ with the point
$\gamma(b)$. The equality of this index with the geometrical index of $\gamma$ is
precisely the statement of the classical Morse Index Theorem.
\end{remMorse}
\end{section}
\begin{section}{Curves of Lagrangians originating
from differential problems}\label{sec:final}
In this section we discuss some necessary conditions for a given curve of
Lagrangians $\ell:[a,b]\mapsto\Lambda$ to arise from an admissible quadruple
for the differential problem. In some cases, in order to produce examples
or counterexamples we will also give sufficient conditions, and in particular
we will exhibit a counterexample to Theorem~\ref{thm:MaslovFocal} when the
hypothesis (\ref{itm:2}) is dropped. Namely, we give an example of an admissible
quadruple $(g,R,P,S)$ in $\R^2$, with $R$ real analytic, $n_+(g)=1$, $P=\{0\}$,
having a unique $(P,S)$-focal instant, and whose Maslov index is equal to $-1$,
while the focal index is zero. The {\em geometric realization\/} of this example
(see Proposition~\ref{thm:correspondence}) is given by a spacelike geodesic
in a three-dimensional real analytic Lorentzian manifold.
\mysubsection{Differential problems determine curves of Lagrangians that are
tangent to distributions of affine spaces}\label{sub:distributions}
Let $g$ be a fixed non degenerate symmetric bilinear form in $\R^n$, and
let $\omega$ be the symplectic form in $\R^{2n}$ given in Definition~\ref{thm:defomega}.
Let $\Psi:[a,b]\mapsto{\rm Sp}(\R^{2n},\omega)$ be a a $C^1$-curve 
such that $\Psi(a)={\rm Id}$; 
a necessary and sufficient condition for such a curve to arise from a quadruple
$(g,R,P,S)$ in the sense of \eqref{eq:defPsit} can be given as follows.
Let $\hat{\mathcal D}_0$ be the right invariant distribution of vector spaces
in ${\rm Sp}(\R^{2n},\omega)$ whose value $\hat{\mathcal D}_0({\rm Id})$
at the identity element is the subspace of ${\rm sp}(\R^{2n},\omega)$ consisting
of the linear operators $H:\R^{2n}\mapsto\R^{2n}$ of the form:
\[H(x,y)=(0,R[x]),\]
for some $g$-symmetric linear operator $R:\R^n\mapsto\R^n$. Observe that $\hat{\mathcal D}_0$
is a distribution of rank $\frac12n(n+1)$. We also define a right invariant
distribution $\hat{\mathcal D}$ of {\em affine subspaces\/} in ${\rm Sp}(\R^{2n},\omega)$
whose value $\hat{\mathcal D}({\rm Id})$ is the affine translation of the vector
space $\hat{\mathcal D}_0({\rm Id})$  by the vector $\bar H\in{\rm sp}(\R^{2n},\omega)$
given by \begin{equation}\label{eq:defHbar}\bar H(x,y)=(y,0).\end{equation}
Keeping in mind formulas \eqref{eq:DEbis} and \eqref{eq:defH}, it is easily
seen that a necessary and sufficient condition for $\Psi$ to arise from a
quadruple $(g,R,P,S)$ is that $\Psi'(t)\in \hat{\mathcal D}(\Psi(t))$ for all $t\in[a,b]$;
we will refer to this situation by saying that $\Psi$ is {\em horizontal\/} with respect to
the distribution $\hat{\mathcal D}$.
\smallskip

We are now going to {\em project\/} the distributions $\hat{\mathcal D}_0$  and
$\hat{\mathcal D}$ to  distributions $\mathcal D_0$ and $\mathcal D$ in $\Lambda$.
Recalling Proposition~\ref{thm:LieAl} and Definition~\ref{thm:defH*}, we consider 
the map:
\begin{equation}\label{eq:bundlemap}
\Lambda\times{\rm sp}(\R^{2n},\omega)\ni(L,H)\longmapsto H^*(L)={\rm d}\kappa_L({\rm Id})[H]
\in T\Lambda;
\end{equation}
it is obvious that this map is a vector bundle morphism, and, due to the
transitivity of the action of ${\rm Sp}(\R^{2n},\omega)$ on $\Lambda$, it
is surjective on each fiber. Moreover, for each $L\in\Lambda$, the kernel
of \eqref{eq:bundlemap} (restricted to the fiber over $L$) is the Lie algebra
${\rm sp}(\R^{2n},\omega,L)$ of the isotropy group of $L$. Namely, it consists
of those $H\in {\rm sp}(\R^{2n},\omega)$ such that $H(L)\subset L$.

We define $\mathcal D$ and $\mathcal D_0$ to be the images of $\Lambda\times\hat{\mathcal D}(
{\rm Id})$ and of $\Lambda\times\hat{\mathcal D}_0({\rm Id})$ respectively
under the map \eqref{eq:bundlemap}. Using the right invariance property
of $\hat{\mathcal D}_0$ and $\hat{\mathcal D}$, it is easily seen that, if
$\ell_0\in\Lambda$ and $\psi\in{\rm Sp}(\R^{2n},\omega)$ are given,
then setting $L=\psi(\ell_0)$, one has:
\[\mathcal D(L)={\rm d}\kappa_{\ell_0}(\psi)[\hat{\mathcal D}(\psi)],\quad
\text{and}\quad \mathcal D_0(L)={\rm d}\kappa_{\ell_0}(\psi)[\hat{\mathcal D_0}(\psi)].\]
Observe that, for each $L\in\Lambda$, $\mathcal D(L)$ is an affine subspace
of $T_L\Lambda$ whose parallel vector subspace is $\mathcal D_0(L)$; 
we emphasize that the dimension of $\mathcal D_0(L)$ is non constant
for $L\in\Lambda$, and so we have distributions of {\em non constant rank}.

Let now $P$ be a $g$-nondegenerate subspace of $\R^n$ and $S$ be a symmetric bilinear form
in $P$; we denote by $\ell_0$ the Lagrangian in $\Lambda$ determined
by $(P,S)$ as in \eqref{eq:ell0}.
Let $\ell:[a,b]\mapsto\Lambda$ be a $C^1$-curve with $\ell(a)=\ell_0$.
Clearly, a necessary condition for $\ell$ to arise from a quadruple $(g,R,P,S)$
in the sense of \eqref{eq:defellt} is that $\ell$ be horizontal with respect to
$\mathcal D$, i.e., $\ell'(t)\in\mathcal D(\ell(t))$ for all $t\in[a,b]$.
We will show that this condition is in general not sufficient.
\smallskip

We now compute explicitly $\mathcal D$ and $\mathcal D_0$.
Let $L\in\Lambda$ be fixed.
Using Definition~\ref{thm:defH*}, we see that $\mathcal D_0(L)$ consist of the restrictions
to $L$ of the bilinear forms $\omega(H\cdot,\cdot)$, where $H$ runs through
$\hat{\mathcal D}_0({\rm Id})$. We compute easily:
\begin{equation}\label{eq:contabis}
\big\{\omega(H\cdot,\cdot):H\in\hat{\mathcal D}_0({\rm Id})\big\}=
\big\{B\in{\rm B}_{\rm sym}(\R^{2n},\R):L_0\subset{\rm Ker}(B)\big\},
\end{equation}
where $L_0=\{0\}\oplus\R^n$. Moreover, the image of $(L,\bar H)$ under
\eqref{eq:bundlemap} is the restriction to $L$ of the bilinear form $0\oplus g$
in $\R^{2n}$. Hence, we have the following description of $\mathcal D_0$
and $\mathcal D$:
\begin{equation}\label{eq:D0D}
\mathcal D_0(L)=\big\{B\in{\rm B}_{\rm sym}(L,\R):L\cap L_0\subset
{\rm Ker}(B)\big\},\quad
\mathcal D=(0\oplus g)\vert_L+\mathcal D_0(L).
\end{equation}
Now it is easy to compute the dimensions of $\mathcal D_0(L)$, for varying $L\in\Lambda$.
Namely, if $L\in\Lambda_0(L_0)$, then $L\cap L_0=\{0\}$, and therefore
$\mathcal D_0(L)=\mathcal D(L)=T_L\Lambda$. More in general, if $L\in\Lambda_k(L_0)$,
then \[{\rm dim}(\mathcal D_0(L))={\rm dim}({\rm B}_{\rm sym}(\R^{n-k},\R))=\frac12
(n-k)(n-k+1),\quad L\in\Lambda_k(L_0);\]
moreover, comparing \eqref{eq:D0D} with Proposition~\ref{thm:geometriaLambdak}, 
we have that, for $L\in\Lambda_k(L_0)$, 
$\mathcal D_0(L)\subset T_L(\Lambda_k(L_0))$.

We consider the surjective linear map
\begin{equation}\label{eq:restriction}
\hat{\mathcal D}_0({\rm Id})\longmapsto\mathcal D_0(L)
\end{equation}
given by the restriction of \eqref{eq:bundlemap}. 

For $L\in\Lambda_0(L_0)$,
then ${\rm dim}(\hat{\mathcal D}_0({\rm Id}))={\rm dim}(\mathcal D_0(L))$, hence
\eqref{eq:restriction} is an isomorphism. More in general, for $L\in\Lambda_k(L_0)$,
the dimension of the kernel of \eqref{eq:restriction} is equal to the codimension
of $\mathcal D_0(L)$ in $T_L\Lambda$. We also consider the surjective affine
map:
\begin{equation}\label{eq:affinemap}
\hat{\mathcal D}({\rm Id})\longmapsto\mathcal D(L),
\end{equation}
defined similarly.

Suppose now that $\ell:[a,b]\mapsto\Lambda$ is a $C^1$-curve with $\ell(a)=\ell_0$
which is horizontal with respect to $\mathcal D$.
By the surjectivity of \eqref{eq:affinemap}, for all $t\in[a,b]$ there exists
a (possibly non unique) $H(t)\in \hat{\mathcal D}({\rm Id})$ mapped onto $\ell'(t)$.
Observe that every such element $H(t)$ defines uniquely a $g$-symmetric
linear map $R(t)$ on $\R^n$ via the formula \eqref{eq:defH}. The only obstruction
for $\ell$ to arise from a quadruple $(g,R,P,S)$ consists precisely in the fact
that one may not be able to make a {\em continuous\/} choice of the maps $H(t)$.
Such obstruction may only occur at the {\em jumps\/} of the function
${\rm dim}(\mathcal D(\ell(t)))$. 

If, for $t$ in a subinterval of $[a,b]$, $\ell(t)\in\Lambda_0(L_0)$, then
there is a unique choice of $H(t)$ on such interval, which is clearly continuous
(such $H(t)$ has the same regularity as $\ell'$). As a matter of facts, one can
prove easily that a continuous choice of $H(t)$ can be made on every interval
for which ${\rm dim}(\mathcal D_0(\ell(t)))$ is constant, even though the choice
of $H(t)$ may not be unique.
\mysubsection{A study of curves of Lagrangians in local coordinates}\label{sub:localcoordinates}
In order to determine sufficient conditions for a curve $\ell$  to arise from
a quadruple $(g,R,P,S)$, we need to study derivatives of $\ell$ of higher
order at the points of intersection with $\Lambda_{\ge1}(L_0)$.
To this aim, we now consider a local chart $\phi_{L_0,L_1}$ (see Definition~\ref{thm:cartas})
where $L_1$ is any Lagrangian complementary to $L_0$;  we will consider
a restriction of $\ell$ whose image lies in the domain $\Lambda_0(L_1)$
of $\phi_{L_0,L_1}$. 

Let $\beta$ be the composition $\phi_{L_0,L_1}\circ\ell$; we write differential
equation in \eqref{eq:DEell} (recall formula~\eqref{eq:defH}) in terms of $\beta$.

By Remark~\ref{thm:remtangent}, the isomorphism ${\rm B}_{\rm sym}(\ell(t),\R)
\simeq {\rm B}_{\rm sym}(L_0,\R)$ given by the differential ${\rm d}\phi_{L_0,L_1}(\ell(t))$
is the pull-back $\eta^*$ by the isomorphism $\eta:L_0\mapsto \ell(t)$ given by 
the restriction of the projection $\ell(t)\oplus L_1\mapsto \ell(t)$. To simplify
the notations, whenever possible we will omit the variable $t$ in the computations that follow.

Recalling Definition~\ref{thm:defH*}, the expression  in coordinates of the right side of the
differential equation  in \eqref{eq:DEell}  is therefore the bilinear form on $L_0$ given by
$\omega(H\eta\,\cdot,\eta\,\cdot)$. Writing $\ell(t)$ as the graph of a linear
map $T:L_0\mapsto L_1$, we have:
\begin{equation}\label{eq:betaT}
\beta=\Iota{L_0}{L_1}\circ T.
\end{equation}
It is now easily seen that $\eta(v)=v+Tv$ for all $v\in L_0$.
Let   $\pi_i:L_0\oplus L_1\mapsto L_i$, $i=0,1$, be the projections; we 
write $H_{ij}=(\pi_i\circ H)\vert_{L_j}$, for $i,j=0,1$.
We now compute $\omega(H\eta\, v,\eta w)$ for all $v,w\in L_0$ as follows:
\[\begin{split}
\omega(H\eta\, v,\eta\, w)=&\;\left[\Iota{L_0}{L_1}\circ H_{10}(v)\right](w)
+\left[\Iota{L_1}{L_0}\circ H_{00}(v)\right](Tw)+\\
&+\left[\Iota{L_0}{L_1}\circ H_{11}(Tv)\right](w)
+\left[\Iota{L_1}{L_0}\circ H_{01}(Tv)\right](Tw).
\end{split}
\]
Using \eqref{eq:Iotas} and \eqref{eq:betaT}, by the above formula we get:
\[\begin{split}
\omega(H\eta\,\cdot,\eta\,\cdot)=&\;\Iota{L_0}{L_1}\circ H_{10}-\beta\circ H_{00}
+\Iota{L_0}{L_1}\circ H_{11}\circ\Iota{L_0}{L_1}^{-1}\circ\beta+\\&-\beta\circ
H_{01}\circ\Iota{L_0}{L_1}^{-1}\circ\beta.
\end{split}
\]
Since $H\in{\rm sp}(\R^{2n},\omega)$, it follows easily $ \Iota{L_0}{L_1}\circ
H_{11}\circ\Iota{L_0}{L_1}^{-1}=-H_{00}^*$; hence the differential equation
for $\beta$ is given by:
\begin{equation}\label{eq:DEbeta}
\beta'=\Iota{L_0}{L_1}\circ H_{10}-\beta\circ H_{00}-H_{00}^*\circ\beta-
\beta\circ H_{01}\circ\Iota{L_0}{L_1}^{-1}\circ\beta.
\end{equation}
By Definition~\ref{thm:defomega} and \eqref{eq:defH}, it is easily checked that:
\begin{equation}\label{eq:primeirotermo} 
\Iota{L_0}{L_1}\circ H_{10}=0\oplus g.
\end{equation}
Moreover, writing $L_1$ as the graph of a $g$-symmetric linear map $Z:\R^n\mapsto\R^n$:
\begin{equation}\label{eq:defL1}
L_1=\big\{(x,Z(x)):x\in\R^n\big\},
\end{equation}
we get:
\begin{equation}\label{eq:segundotermo}
H_{00}=-0\oplus Z,\quad H_{01}\circ(\Iota{L_0}{L_1})^{-1}=0\oplus\big[(R-Z^2)\circ g^{-1}\big],
\end{equation}
where in the above formulas $g$ is considered as the map $\R^n\mapsto(\R^n)^*$
given by $x\mapsto g(x,\cdot)$.

Using \eqref{eq:primeirotermo} and \eqref{eq:segundotermo} where $L_0=\{0\}\oplus\R^n$ 
is identified with $\R^n$, we rewrite \eqref{eq:DEbeta} as:
\begin{equation}\label{eq:secundaria}
\beta'=g+\beta\circ Z+Z^*\circ\beta-\beta\circ(R-Z^2)\circ g^{-1}\circ\beta.
\end{equation}
Equation \eqref{eq:secundaria} is the translation in coordinates of the differential
equation in \eqref{eq:DEell}; a $C^1$-curve $\beta$ is such that the corresponding curve
$\ell$ is horizontal with respect to $\mathcal D$ if and only if for each $t$ there
exists a $g$-symmetric $R(t):\R^n\mapsto\R^n$ satisfying \eqref{eq:secundaria}.

We now concentrate our attention to the problem of determining conditions on $\beta$
that guarantee the existence of a continuous choice of maps $R(t)$ as above satisfying
\eqref{eq:secundaria}. A first necessary condition to the existence of $R(t)$
is obtained by evaluating \eqref{eq:secundaria} at a pair $(v,w)\in{\rm Ker}(\beta)\times
\R^n$:
\begin{equation}\label{eq:CN1}
\beta'(v,w)=g(v,w)+\beta(Zv,w),\quad\forall v\in{\rm Ker}(\beta),\ w\in\R^n.
\end{equation}
Condition \eqref{eq:CN1} is simply a coordinate version of the horizontality of $\ell$
(compare with \eqref{eq:D0D}). Let's assume now that $\beta$ is a curve
of class $C^2$; we determine a necessary condition for the existence of
a curve $R(t)$ of class $C^1$ satisfying \eqref{eq:secundaria}. We differentiate
\eqref{eq:secundaria} and we evaluate at a pair of vectors $v,w\in{\rm Ker}(\beta)$,
obtaining:
\[\beta''(v,w)=\beta'(Zv,w)+\beta'(v,Zw),\quad \forall\,v,w\in{\rm Ker}(\beta).\]
Using \eqref{eq:CN1}, the above formula becomes:
\begin{equation}\label{eq:CN2}
\beta''(v,w)=2\, g(Zv,w)+2\, \beta(Zv,Zw),\quad\forall\,v,w\in{\rm Ker}(\beta).
\end{equation}
At the instants $t$ where $\beta$ is invertible, the unique ($g$-symmetric) 
map $R(t)$ satisfying
\eqref{eq:secundaria} is computed as:
\begin{equation}\label{eq:principal}
R=\beta^{-1}\circ\left(g-\beta'+\beta\circ Z+Z^*\circ\beta\right)\circ\beta^{-1}\circ g+Z^2.
\end{equation}
Observe that ${\rm Ker}(\beta(t))=\ell(t)\cap L_0$; hence, the
condition that $\beta$ be invertible means that $\ell(t)\in\Lambda_0(L_0)$.
As we have observed earlier, in this case there is no obstruction to the
existence of the map $R(t)$.

For simplicity, we will now restrict to the case that $\beta$ is smooth, that
it is not invertible for only a finite number of instants $t$ and that
${\rm det}(\beta)$ has only zeroes of finite order. For instance,
this is the case if $\beta$ is real analytic and if ${\rm det}(\beta)$
is not identically zero.

Let $t_0$ be a fixed instant at which $\beta$ is {\em not\/} invertible.
For $t\sim t_0$, $t\ne t_0$, we write \eqref{eq:principal} in matrix form
(using a suitable basis of $\R^n$); the entries of $R(t)$ will then be given
by quotients of smooth functions of $t$. A necessary and sufficient
condition for the existence of a smooth extension of $R$ at the instant $t_0$ is
that in these quotients the order of zero of the functions appearing at the 
numerator be greater than or equal to the order of zero of the functions at the 
denominator. In this situation, it is obvious that necessary and sufficient
conditions for the existence of a smooth extension of $R$ can be given in terms
of certain nonlinear systems of equations involving higher order derivatives
of the coefficients of $\beta$ at $t=t_0$. It is interesting to observe that,
if $\beta$ is real analytic, then so is $R$.
\mysubsection{The case where $g$ is nondegenerate on ${\rm Ker}(\beta(t_0))$}
\label{sub:nondegenerate}
We temporarily make the extra assumption that $g$, or equivalently $\beta'(t_0)$,
be non degenerate on ${\rm Ker}(\beta(t_0))$. We can prove then that conditions
\eqref{eq:CN1} and \eqref{eq:CN2} are sufficient. Towards this goal,
let $e_1,\ldots,e_n$ be a basis of $\R^n$ such that $e_1,\ldots,e_k$ is
a basis of ${\rm Ker}(\beta(t_0))$; the restriction of $\beta(t_0)$ to 
the space spanned by the $e_{k+1},\ldots,e_n$ is clearly non degenerate.
We will now think of all our bilinear forms as matrices relative to this basis.

For $t\ne t_0$, let $\tilde\beta(t)$ be the matrix obtained by dividing 
the first $k$ columns
of $\beta(t)$ by $(t-t_0)$; we define $\tilde\beta(t_0)$ by replacing the
first $k$ columns of $\beta(t_0)$ by the first $k$ columns of $\beta'(t_0)$. 
It is easy to see that $\tilde\beta$ is smooth.

If we define $D_k(a)=\left(\begin{array}{cc} a\cdot I_k&0\\ 0&I_{n-k}\end{array}\right)$,
where $I_j$ denotes the $j\times j$ identity matrix and $a\in\R$, then we can write:
\[\beta(t)=\tilde\beta(t) D_k(t-t_0),\quad\forall\, t\ne t_0,\]
hence:
\begin{equation}\label{eq:Dkesq}
\beta(t)^{-1}=D_k\left(\frac1{t-t_0}\right)\,\tilde\beta(t)^{-1},\quad\forall\,
t\ne t_0.
\end{equation} 
Since $g$ is non degenerate on ${\rm Ker}(\beta(t_0))$, it follows from
\eqref{eq:CN1} that $\tilde\beta(t_0)$ is invertible, and so $\tilde\beta(t)^{-1}$ 
is smooth. By \eqref{eq:Dkesq}, the last $n-k$ lines of $\beta(t)^{-1}$ are smooth,
and the first $k$ lines have {\em a singularity of order at the most one\/} at $t=t_0$,
i.e., they are the quotient of smooth functions by $t-t_0$.

By the symmetry of $\beta(t)^{-1}$, it actually follows that the last
$n-k$ columns of $\beta(t)^{-1}$ are smooth, from which it follows that
the singularities of $\beta(t)^{-1}$ are concentrated in the upper left
$k\times k$ block, and all the singularities are of order at the most one. 

We denote by $Q$ the following symmetric bilinear form:
\[Q=g-\beta'+\beta\circ Z+Z^*\circ\beta;\]
formula \eqref{eq:principal} can be rewritten in terms of $Q$ as:
\begin{equation}\label{eq:novaR}
R=\beta^{-1}\circ Q \circ\beta^{-1}\circ g+Z^2.
\end{equation}
We observe that $Q$ is smooth, its coefficients of the first $k$ lines and
of the first $k$ columns have zeroes of order at least one at $t=t_0$, and
that the coefficients of the upper left $k\times k$ block have zeroes of order at least
two at $t=t_0$.

From \eqref{eq:novaR} it now follows easily that $R$ is smooth, which proves the claim.
\smallskip

For the above argument, the crucial hypothesis of nondegeneracy of the restriction of
$g$ to the kernel of $\beta(t_0)$ cannot be avoided; if this condition is not
satisfied, in order to get to the conclusion one needs to analyze the behavior of
derivatives of higher order of $\beta$ at $t=t_0$. 

Under the assumption that $g$ be positive definite, i.e., when $g$ is related to
a Riemannian or a causal Lorentzian geodesic problem, by the above argument
\eqref{eq:CN1} and \eqref{eq:CN2} characterize completely the
curves $\beta$ arising from a quadruple $(g,R,P,S)$.

\mysubsection{A counterexample to the equality $\mu(g,R,P,S)={\rm i}_{\rm foc}$}
\label{sub:counterexample}

As announced at the beginning of the section, we now pass to the construction
of a counterexample to Theorem~\ref{thm:MaslovFocal} when the
hypothesis (\ref{itm:2}) is dropped. 

We consider the following setup. Let $n=2$, the objects 
$g$, $Z$ and $P$ that we consider are:
\[g=\begin{pmatrix}0&1\\ 1&0\end{pmatrix},\quad
Z=\begin{pmatrix}0&0\\ 1&0 \end{pmatrix}\quad \text{and}\quad P=\{0\};\]
all the matrices involved are relative to the canonical basis of $\R^2$.
Observe that $Z$ is $g$-symmetric, i.e., the matrix $gZ$ is symmetric.

For our purposes, we will construct a curve $\ell$ in $\Lambda$ whose image
is entirely contained in the domain of the chart $\phi_{{L_0},{L_1}}$, where
$L_0=0\oplus\R^2$ and $L_1$ is defined by \eqref{eq:defL1}.
It now suffices to describe the curve $\beta$ in ${\rm B}_{\rm sym}(\R^2,\R)\simeq\R^3$;
we write:
\[\beta(t)=\begin{pmatrix}x(t)&z(t)\\z(t)&y(t) \end{pmatrix},\]
where $x$, $y$ and $z$ are real analytic scalar functions on an interval
$[a,b]$, with $a<0<b$, such that the 
following properties are satisfied:
\begin{itemize}
\item[(1)] $\beta(a)=0$, which means that the initial condition $\ell(a)=\ell_0=L_0$
is satisfied;
\item[(2)] $\beta'(a)=g$ and $\beta''(a)=2gZ=\begin{pmatrix}2&0\\ 0&0 \end{pmatrix}$, i.e.,
conditions \eqref{eq:CN1} and \eqref{eq:CN2} are satisfied at $t=a$;
\item[(3)] ${\rm det}(\beta(t))$ has zeroes precisely at $t=a$ and at $t=0$;
\item[(4)] ${\rm Ker}(\beta(0))$ is generated by the first vector of the
canonical basis of $\R^2$, i.e., $x(0)=z(0)=0$ and $y(0)\ne0$;
\item[(5)] $x'(0)=0$, $z'(0)=1+y(0)$ and $x''(0)=2+2\,y(0)$, i.e., conditions
\eqref{eq:CN1} and \eqref{eq:CN2} are satisfied at $t=0$;
\item[(6)] ${\rm det}(\beta(t))$ has a zero of order precisely $3$ at $t=0$, and
its third derivative is positive at $t=0$. This is 
equivalent to $y(0)=-1$ and  $x'''(0)<0$;
\item[(7)] the function $R$ given by \eqref{eq:novaR}, or equivalently
the function $\beta^{-1}\circ Q \circ\beta^{-1}$ is non singular
at $t=0$.
\end{itemize}
Conditions (5) e (6) imply also $z'(0)=x''(0)=0$; conditions (1) and (2) imply that
${\rm det}(\beta(t))$ has a zero of order two at $t=a$.
\smallskip

We will now show that it is possible to determine polynomial functions
$x$, $y$ and $z$ satisfying all the above conditions. We proceed by steps as follows.

Our interval $[a,b]$ will be of the form $[-1,b]$, with $b>0$ sufficiently small.
Once made a choice of functions $x$, $y$ and $z$ so that (1)---(7) are satisfied,
the endpoint $b$ will be chosen in such a way that ${\rm det}(\beta(t))$
is strictly positive in $]0,b]$ (observe that this is possible by condition (6)). 

We denote by $\hat\beta$ the following matrix:
\[\hat\beta(t)=\begin{pmatrix}y(t)&-z(t)\\-z(t)&x(t) \end{pmatrix},\]
so that $\beta^{-1}=({\rm det}(\beta))^{-1}\cdot\hat\beta$ whenever $\beta$
is invertible. By condition (6) above, ${\rm det}(\beta)$ has a zero
of order $3$ at $t=0$, so, in order to satisfy (7), a necessary and sufficient
condition is that the entries of the matrix $\hat\beta\circ Q\circ\hat\beta$ have 
zeroes of order at least $6$ at $t=0$. If we write:
\[\hat\beta\circ Q\circ\hat\beta\circ g=\begin{pmatrix}p_1(t)&p_2(t)\\
p_3(t)&p_1(t) \end{pmatrix},\]
we have:
\[
\begin{split}
& p_1=x'yz-yz^2+xy-xyz'+xy^2+z^2-z^2z'+xy'z,\\
& p_2=-x'y^2-2yz+2yzz'-y'z^2,\\
& p_3=-x'z^2+2z^3-2xz+2xzz'-2xyz-x^2y''.
\end{split}
\]
The conditions (4) through (7) above are satisfied, for instance, with the choice:
\begin{equation}\label{eq:xyz}
x(t)=-2t^3-\frac{54}5 t^5,\quad y(t)=-1-6t+18t^2-54t^3,\quad
z(t)=-3t^2,\quad t\sim0.
\end{equation}
It is easy now to see that it is possible to choose  smooth functions
$x$, $y$ and $z$ that coincide with the polynomials given in \eqref{eq:xyz}
around $t=0$ and such that conditions (1)---(3) are also satisfied.
By what we have observed so far, such choice provides a smooth counterexample
to Theorem~~\ref{thm:MaslovFocal} when the nondegeneracy assumption is dropped.


For the final step of our real analytic counterexample we now argue abstractly
using a density argument of polynomials, as follows.

Let $x$, $y$ and $z$ be given smooth functions so that conditions (1)---(7)
are satisfied; we start observing that if $\tilde x$, $\tilde y$ and $\tilde z$
are smooth functions having the first six derivatives at $t=0$ and the
first two derivatives at $t=-1$ equal to the corresponding derivatives
 of $x$, $y$ and $z$, then,
replacing $x$, $y$ and $z$ by $\tilde x$, $\tilde y$ and $\tilde z$, only
condition (3) may fail to hold. If such a replacement is done in such a way that
$\tilde x$, $\tilde y$ and $\tilde z$ are sufficiently close to $x$, $y$ and
$z$ in the $C^4$-topology, then also condition (3) will remain true.
To prove this, we apply the next lemma to the function $f=xy-z^2$ on the
interval $[-1,0]$ with $k=3$:
\begin{Ck}\label{thm:Ck}
Let $k\in\N$ and let $f:[a,b]\mapsto\R$ be a function of class $C^{k+1}$.
Assume that $f$
has zeroes precisely at the endpoints $a$,$b$ and that these zeroes
are of order at the most $k$. Then, there exists a neighborhood
$\mathcal U$ of $f$ in the $C^{k+1}$-topology such that,
every $g\in\mathcal U$ having the same order of zeroes as $f$ at $a$ and $b$
has no zeroes in $]a,b[$.
\end{Ck}
\begin{proof}
Let $i$ and $j$ be the order of zeroes of $f$ at $a$ and $b$ respectively,
$i,j\in\{1,\ldots,k\}$. Define the following constant:
\[M=\max\{\Vert f^{(i+1)}\Vert_\infty,\ \Vert f^{(j+1)}\Vert_\infty\},\]
and let $\delta_1,\delta_2>0$ be such that $a<a+\delta_1<b-\delta_2<b$
and
\[\delta_1<\frac{(i+1)\vert f^{(i)}(a)\vert}{2 (M+1)},\quad
\delta_2<\frac{(j+1)\vert f^{(j)}(b)\vert}{2 (M+1)}.\]
Finally, let $\alpha>0$ be the infimum of $\vert f\vert$ on the interval
$[a+\delta_1,b-\delta_2]$. The desired neighborhood of $f$ is defined by
requiring that $g\in\mathcal U$ if and only if:
\[\begin{split}
&\vert g^{(i)}(a)\vert>\frac{\vert f^{(i)}(a)\vert}2,\ 
\vert g^{(j)}(b)\vert>\frac{\vert f^{(j)}(b)\vert}2,\\ \\ & 
\Vert g^{(i+1)}-f^{(i+1)}\Vert_\infty<1,\ \Vert g^{(j+1)}-f^{(j+1)}\Vert_\infty<1,\ 
\Vert g-f\Vert_\infty<\alpha.
\end{split}\]
To check that this choice of $\mathcal U$ works, let $g\in\mathcal U$
be chosen so that $g$ has a zero of order $i$ at $a$ and a zero of order $j$
at $b$. Using the $i$-th order Taylor polynomial of $g$ around $a$, we get:
\[g(t)=(t-a)^i\left(\frac{g^{(i)}(a)}{i!}+r(t)\right),\]
where $r(t)$ satisfies:
\[\vert r(t)\vert\le\sup_{[a,t]}\frac{\vert g^{(i+1)}\vert}{(i+1)!}\cdot (t-a)\le
\frac{M+1}{(i+1)!}\cdot (t-a).\]
By our choice of $\delta_1$, it follows that $g$ has no zeroes in $]a,a+\delta_1]$;
similarly, $g$ has no zeroes in $[b-\delta_2,b[$. 

From $\Vert g-f\Vert_\infty<\alpha$, it follows that $g$ has no zeroes
in $[a+\delta_1,b-\delta_2]$, which concludes the proof.
\end{proof}
Finally, for the construction of our analytic counterexample, we use the observations
above, and a simple density result which is contained in the following Lemma:
\begin{dense}\label{thm:dense}
Let $k\in\N$ and $a_i,b_i\in\R$, $i=0,\ldots,k$, be fixed. Consider the following
subsets of $C^k([a,b],\R)$:
\[
\begin{split}
&A=\Big\{f\in C^k([a,b],\R):f^{(i)}(a)=a_i,\ f^{(i)}(b)=b_i, \ \ i=0,\ldots,k\Big\},
\\
&B=\Big\{f\in A: f\ \text{is a polynomial}\ \Big\}
\end{split}
\]
Then, $B$ is dense in $A$ in the $C^k$-topology.\qed
\end{dense}
We can therefore build a real analytic curve $\ell:[-1,b]\mapsto\Lambda$
which arises from an admissible quadruple $(g,R,P,S)$, with $R$ real analytic.
By condition (3) the only $(P,S)$-focal instant occurs at $t=0$; the restriction
of $g$ to the $\ell(0)\cap L_0$ is zero by condition (4). It follows that
the focal index of the quadruple is zero. On the other hand, by condition (6),
the sign of ${\rm det}(\beta(t))$ changes from negative to positive as $t$ passes
through $0$; moreover, the trace of $\beta(t)$ is negative around $t=0$. By
Proposition~\ref{thm:method}, this implies that the Maslov index of the quadruple 
is $-1$.

\mysubsection{Instability of the focal index}
\label{sub:instabilityfocal}
Let $(g,R,P,S)$ be the quadruple constructed in the previous subsection.
A {\em small\/} perturbation of $(g,R,P,S)$  preserves the Maslov index, by
Theorem~\ref{thm:LimitMaslov}. However, we observe that the focal index may change,
by the following arguments.

If we identify $\beta$ with a curve in $\R^3$, then the $(P,S)$-focal instants
occur precisely
at the intersections of this curve with the double cone $xy-z^2=0$. 
Given one such intersection $\beta(t_0)$, the degeneracy of $g$ on 
${\rm Ker}(\beta(t_0))$ means that $\beta(t_0)\ne0$ belongs to one of the straight lines
$x=z=0$ and $y=z=0$.

Such condition is evidently unstable by small perturbations, and a 
quadruple obtained from $(g,R,P,S)$ by a small perturbation will generically
satisfy the hypotheses of Theorem~\ref{thm:MaslovFocal} and therefore, its
focal index will be equal to $-1$.

\mysubsection{Instability of focal points with signature zero}
\label{sub:instabilitypoints}
Let's assume that $n=2$ and that $g$ is symmetric bilinear form of signature $0$,
i.e., ${\rm det}(g)<0$. An instant $t_0$ such that $\beta(t_0)=0$ is a $(P,S)$-focal
instant of signature $0$. Again, it is fairly obvious that a small perturbation
of the curve $\beta$ may not intersect the double cone $xy-z^2$ around $t_0$,
which amounts to say that a $(P,S)$-focal instant with signature $0$ may 
{\em evaporate\/} by small perturbations of the quadruple $(g,R,P,S)$.
\end{section}



\end{document}